\documentclass[10pt]{amsart}
\usepackage{latexsym, amsmath}
\usepackage{amssymb}

\newcommand{\varespilon}{\varepsilon}

\renewcommand{\epsilon}{\varepsilon}
\newcommand{\newsection}[1]
{\section{#1}\setcounter{theorem}{0} \setcounter{equation}{0}
\par\noindent}

\newtheorem{theorem}{Theorem}

\newtheorem{lemma}[theorem]{Lemma}
\newtheorem{corr}[theorem]{Corollary}

\newtheorem{proposition}[theorem]{Proposition}
\newtheorem{deff}[theorem]{Definition}

\newcommand{\bth}{\begin{theorem}}
\newcommand{\ble}{\begin{lemma}}
\newcommand{\bcor}{\begin{corr}}

\newcommand{\bdeff}{\begin{deff}}

\newcommand{\bprop}{\begin{proposition}}
\newcommand{\ele}{\end{lemma}}
\newcommand{\ecor}{\end{corr}}
\newcommand{\edeff}{\end{deff}}

\newcommand{\eprop}{\end{proposition}}

\newcommand{\cd}{\, \cdot\, }

\renewcommand{\Pi}{\varPi}

\renewcommand{\epsilon}{\varepsilon}

\newcommand{\Rplus}{{\Bbb R}_+}

\newcommand{\tidle}{\tilde}

\newcommand{\R}{{\Bbb R}}
\newcommand{\ext}{{\R^3\backslash\mathcal{K}}}

\begin{document}

\title[Quasilinear wave equations satisfying the null condition]
{Global existence of quasilinear, nonrelativistic wave equations
satisfying the null condition}
\thanks{The first and third authors were supported in part by the NSF}

\author{Jason Metcalfe}
\address{School of Mathematics, Georgia Institute of Technology,
Atlanta, GA  30332-0160}
\author{Makoto Nakamura}
\address{Graduate School of Information Sciences, Tohoku
University, Sendai 980-8579, Japan}
\author{Christopher D. Sogge}
\address{Department of Mathematics,  Johns Hopkins University,
Baltimore, MD 21218}

\maketitle

\newsection{Introduction}
The purpose of this paper is to provide a proof of global existence of
solutions to general quasilinear, multiple speed systems of wave
equations satisfying the null condition.  The techniques presented are
sufficient to handle both Minkowski wave equations and Dirichlet-wave
equations in the exterior of certain compact obstacles.

For the latter case, fix a smooth, compact obstacle
$\mathcal{K}\subset \R^3$.  We, then, wish to examine the quasilinear
system
\begin{equation}\label{1.1}
\begin{cases}
\Box u = F(u,du,d^2 u), \quad(t,x)\in \R_+\times\ext\\
u(t,\cd)|_{\partial\mathcal{K}}=0\\
u(0,\cd)=f, \quad \partial_t u(0,\cd)=g.
\end{cases}
\end{equation}
Here
\begin{equation}\label{1.2}
\Box=(\Box_{c_1},\Box_{c_2},\dots,\Box_{c_D})
\end{equation}
denotes a vector-valued multiple speed d'Alembertian where
$$\Box_{c_I}=\partial_t^2-c_I^2\Delta$$
and
$\Delta=\partial_1^2+\partial_2^2+\partial_3^2$ is the standard
Laplacian.  For clarity, we will assume that we are in the
nonrelativistic case.  That is, we assume that the wave speeds $c_I$
are positive and distinct.  Straightforward modifications can be made
to allow various components to have the same speed.  For convenience,
we will take $c_0=0$ and 
\begin{equation}\label{order}
0=c_0<c_1<c_2<\dots <c_D
\end{equation}
throughout.

We now describe our conditions on the nonlinearity $F$.  First of all,
$F$ is assumed to be linear in $d^2u$. 
$F$ is also required to vanish to second order.  That is,
$$\partial^\alpha F(0,0,0)=0,\quad |\alpha|\le 1.$$
Additionally, we assume
$$\partial_u^2 F(0,0,0)=0.$$
Thus, $F$ may be decomposed as
$$F(u,du,d^2u)=B(du)+Q(du,d^2u)+R(u,du,d^2u)+P(u,du)$$
where, for $1\le I\le D$, 
\begin{equation}\label{1.3}
B^I(du)=\sum_{\substack{1\le J,K\le D\\0\le j,k\le 3}}
A^{I,jk}_{JK}\partial_ju^J\partial_k u^K,
\end{equation}
\begin{equation}\label{1.4}
Q^I(du,d^2u)=\sum_{\substack{1\le J,K\le D\\ 0\le j,k,l\le 3}}
B_{K,l}^{IJ,jk} \partial_l u^K \partial_j\partial_k u^J,
\end{equation}
\begin{equation}\label{1.5}
R^I(u,du,d^2u)=\sum_{\substack{1\le J\le D\\ 0\le j,k\le 3}} 
C^{IJ,jk}(u,du) \partial_j\partial_k u^J
\end{equation}
with $C^{IJ,jk}(u,du)=O(|u|^2+|du|^2)$,
and $P(u,du)=O(|u|^3 + |du|^3)$ near $(u,du)=0$.  Here and throughout, we use the
notation $x_0=t$ and $\partial_0=\partial_t$ when convenient.
Additionally, $du=u'=\nabla_{t,x} u$ denotes the space-time
gradient.  The constants $B^{IJ,jk}_{K,l}$ are
real, as are the $C^{IJ,jk}(u,du)$ terms.
Moreover, the quasilinear
terms are assumed to satisfy the symmetry conditions
\begin{equation}\label{1.6}
B^{IJ,jk}_{K,l}=B^{JI,jk}_{K,l}=B^{IJ,kj}_{K,l},
\end{equation}
\begin{equation}\label{1.7}
C^{IJ,jk}(u,du)=C^{IJ,kj}(u,du)=C^{JI,jk}(u,du).
\end{equation}
In order to establish global existence, we require that the quadratic
terms satisfy the following null condition:
\begin{equation}\label{1.8}
\sum_{0\le j,k\le 3} A^{J,jk}_{JJ}\xi_j\xi_k =0,\quad \text{whenever }
\frac{\xi_0^2}{c_J^2} - \xi_1^2-\xi_2^2-\xi_3^2=0, \quad J=1,2,\dots, D,
\end{equation}
\begin{equation}\label{1.9}
\sum_{0\le j,k,l\le 3} B^{JJ,jk}_{J,l}\xi_j\xi_k\xi_l =0,\quad \text{whenever }
\frac{\xi_0^2}{c_J^2} - \xi_1^2-\xi_2^2-\xi_3^2=0, \quad J=1,2,\dots, D.
\end{equation}
This null condition guarantees that the self-interaction of each wave
family is nonresonant and is the natural one for systems of
quasilinear wave equations with multiple speeds.  It
is equivalent to the requirement that no plane wave solution of the
system is genuinely nonlinear.  This follows from an observation of
John and Shatah, and we refer the reader to John \cite{John} (p. 23)
and Agemi-Yokoyama \cite{AY}.  Additionally, in the setting of
elasticity, Tahvilday-Zadeh \cite{TZ} (see also Sideris \cite{Si}) 
 observed that \eqref{1.8}, 
\eqref{1.9} removed the physically unrealistic restrictions on the
growth of the stored energy imposed by the null conditions used, for
example, in \cite{MS}, \cite{Si2}, and \cite{So2}.  While general global
existence of solutions to \eqref{1.1} is only known (even in the
Minkowski setting) under the assumption of \eqref{1.8}, \eqref{1.9},
recent works of Lindblad-Rodnianski \cite{LR1, LR2} suggest that a weak form of the
null condition may be sufficient.

We now wish to describe our assumptions on the obstacle
$\mathcal{K}\subset \R^3$.  As mentioned above, we assume that
$\mathcal{K}$ is smooth and compact, but not necessarily connected.
By shifting and scaling, we may take
$$0\in \mathcal{K}\subset \{|x|<1\}$$
with no loss of generality.  The only additional assumption is that
there is exponential decay of local energy.  Specifically, if $u$ is a
solution to the homogeneous wave equation
$$\begin{cases}
\Box u=0\\
u(t,\cd)|_{\partial\mathcal{K}}=0\end{cases}$$
and the Cauchy data $u(0,\cd)$, $\partial_tu(0,\cd)$ are supported in
$\{|x|<4\}$, then we assume that there are constants $c,C>0$ so that
\begin{equation}\label{1.10}
\Bigl(\int_{\{x\in \ext\, :\, |x|<4\}} |u'(t,x)|^2\:dx\Bigr)^{1/2}\le
C e^{-ct} \sum_{|\alpha|\le 1} \|\partial^\alpha_x u'(0,\cd)\|_2.
\end{equation}

If the obstacle is nontrapping, a stronger version of \eqref{1.10}
holds with $|\alpha|=0$ (no loss of derivative).  See, e.g.,
Morawetz-Ralston-Strauss \cite{MRS}.  In the presence of trapped rays,
Ralston \cite{ralston} observed that this stronger version could not
hold, and Ikawa
\cite{Ikawa1, Ikawa2} showed that \eqref{1.10} holds for certain
finite unions of convex obstacles.  

In order to solve \eqref{1.1}, we must also require that the data
satisfies certain compatibility conditions.  Briefly, if we let
$J_ku=\{\partial^\alpha u \, : \, 0\le |\alpha|\le k\}$ and fix $m$,
we can write $\partial_t^k u(0,\cd)=\psi_k(J_k f,J_{k-1}g)$, $0\le
k\le m$ for any formal $H^m$ solution of \eqref{1.1}.  Here, $\psi_k$
is called a compatibility function and depends on $F$, $J_k f$, and
$J_{k-1}g$.  The compatibility condition for \eqref{1.1} with
$(f,g)\in H^m\times H^{m-1}$ states that the $\psi_k$ vanish on
$\partial\mathcal{K}$ when $0\le k\le m-1$.  Additionally, we say that
$(f,g)\in C^\infty$ satisfy the compatibility condition to infinite
order if this holds for all $m$. See, e.g., \cite{KSS} for a more 
detailed description of the compatibility conditions.

We can now state our main result.
\begin{theorem}\label{theorem1.1}
Let $\mathcal{K}$ be a fixed compact obstacle with smooth boundary
satisfying \eqref{1.10}.  Assume that $F(u,du,d^2u)$ and $\Box$ are as
above and that $(f,g)\in C^\infty (\ext)$ satisfy the compatibility
conditions to infinite order.  Then, there is an $\varepsilon_0>0$ and
an integer $N>0$ so that for all $\varepsilon < \varepsilon_0$, if
\begin{equation}\label{1.11}
\sum_{|\alpha|\le N}\|\langle x\rangle^{|\alpha|}\partial^\alpha_x
f\|_2 + \sum_{|\alpha|\le N-1} \|\langle
x\rangle^{1+|\alpha|}\partial^\alpha_x g\|_2\le \varepsilon,
\end{equation}
then \eqref{1.1} has a unique global solution $u\in C^{\infty}([0,\infty)\times\ext)$.
\end{theorem}

As mentioned above, we will also handle the Minkowski case.  Assuming
that $F$ and $\Box$ are as above, we show that solutions of
\begin{equation}\label{1.12}
\begin{cases}
\Box u = F(u,du,d^2u),\quad (t,x)\in \R_+\times\R^3\\
u(0,\cd)=f,\quad \partial_t u(0,\cd)=g
\end{cases}
\end{equation}
exist globally for small data.  Specifically, we will prove
\begin{theorem}\label{theorem1.2}
Assume that $F$ and $\Box$ are as above.  Then, there are constants
$\varepsilon_0, N>0$ so if $f,g$ are smooth functions satisfying
\begin{equation}\label{1.13}
\sum_{|\alpha|\le N}\|\langle x\rangle^{|\alpha|}\partial^\alpha_x
f\|_2 + \sum_{|\alpha|\le N-1} \|\langle
x\rangle^{1+|\alpha|}\partial^\alpha_x g\|_2\le \varepsilon,
\end{equation}
for all $\varepsilon < \varepsilon_0$, then the system \eqref{1.12} has a
unique global solution $u\in C^\infty([0,\infty)\times\R^3)$.
\end{theorem}

We note that during preparation of this paper it was discovered that
Theorem \ref{theorem1.2} was proven independently by Katayama \cite{K1} using
different techniques.  Additionally, in \cite{K2},
Katayama explored the possibility of allowing $F$ to contain certain
terms of the form $u^J \partial u^K$ if you assume the null condition
of \cite{Si2}, \cite{So2} rather than \eqref{1.8}, \eqref{1.9}.  The
obstacle result, Theorem \ref{theorem1.1}, is new.

By allowing general higher order terms, Theorem \ref{theorem1.2} 
extends the previously known results on multiple speed wave equations
due to Sideris-Tu \cite{Si3}, Agemi-Yokoyama \cite{AY},
Kubota-Yokoyama \cite{KY}, and Katayama \cite{K3}.    In a
similar way, Theorem \ref{theorem1.1} extends the previous result of
the authors \cite{MNS}.  

In studying both the Minkowski setting and the exterior domain, 
we will be using modifications of the method of
commuting vector fields due to Klainerman \cite{knull}. 
We will
restrict to the class of vector fields $\Gamma=\{Z,L\}$ that seem
``admissible'' for boundary value problems and studies of multiple
speed wave equations.  Here, $Z$ denotes the generators of space-time
translations and spatial rotations
\begin{equation}\label{1.14}
Z=\{\partial_i, x_j\partial_k - x_k\partial_j, \quad 0\le i\le 3, 1\le
j<k\le 3\}
\end{equation}
and $L$ is the scaling vector field
\begin{equation}\label{1.15}
L=t\partial_t + r\partial_r.
\end{equation}
Additionally, we will write $r=|x|$ and
\begin{equation}\label{1.16}
\Omega_{jk}=x_j\partial_k - x_k \partial_j
\end{equation}
for the generators of spatial rotations.  The generators of the
Lorentz rotations, $x_i\partial_t + t\partial_i$ when $c_I=1$, 
have an associated
speed and have unbounded normal components on the boundary of our compact
obstacle, and thus seem ill-suited to the problems in question.
Katayama \cite{K1, K2} has shown that these hyperbolic rotations can
be used in a limited fashion in the study of multiple speed wave
equations, but we do not require those techniques here.

The most significant new difficulty in this case versus the one
considered in \cite{MNS} is
the cubic terms not involving derivatives.  Those involving
derivatives can generally be handled using energy methods.  In
the approaches of Christodoulou \cite{christ} and Klainerman
\cite{knull}, such terms not involving derivatives 
were handled with a certain adapted energy
inequality that resembles, e.g., the work of Morawetz \cite{M}.  This
method relies on the use of the Lorentz rotations, and it is not clear
how to adapt it to the current setting.

The new argument that we utilize uses an analog of a 
pointwise estimate that was established by Kubota-Yokoyama \cite{KY}.  When combined
with the pointwise estimates of Keel-Smith-Sogge \cite{KSS3} and sharp
Huygens' principle, we are able to establish low regularity decay of
our solution $u$.  These improved estimates allow us to handle the
cubic terms without derivatives discussed in the previous paragraph.  
In \cite{MNS}, using only the estimates of
\cite{KSS3}, the authors were only able to get such decay for the gradient of
the solution $u'$.  

As in Keel-Smith-Sogge \cite{KSS2, KSS3}, we will utilize a class of
weighted $L^2_tL^2_x$-estimates where the weight is a negative power
of $\langle x\rangle = \langle r\rangle = \sqrt{1+r^2}$.  Such
estimates permit us to use the $O(\langle x\rangle^{-1})$ decay that
is obtained from Sobolev inequalities rather than the more standard
$O(t^{-1})$ decay which is difficult to prove without the use of the
Lorentz rotations.  Additionally, such estimates allow us to handle
the boundary terms that arise in the energy estimates of {\em
nonlinear} wave equations if we no longer have the convenient
assumption of star-shapedness on the obstacle.  This was one of the
main innovations of Metcalfe-Sogge \cite{MS}.

As in our previous work \cite{MNS}, we will require a class of
weighted Sobolev estimates.  The weights involve powers of $r$ and
$\langle t-r\rangle$.  In the Minkowski setting, these estimates are
originally due to Klainerman-Sideris \cite{KS} and Hidano-Yokoyama
\cite{HY}.

This paper is organized as follows.  In the next section, we gather
our preliminary estimates that will be needed to show global existence
in Minkowski space.  In particular, we collect the pointwise estimates
of Keel-Smith-Sogge \cite{KSS3} and Kubota-Yokoyama \cite{KY}.  
In Section 3, we prove Theorem \ref{theorem1.2}.  In
Section 4, we gather the estimates that we will require to prove
Theorem \ref{theorem1.1}.  Finally, in Sections 5-7, we prove our main
theorem, Theorem \ref{theorem1.1}.


\newsection{Preliminary estimates in Minkowski space}
In this section we gather the estimates for the free wave equation
that we will require in order to prove global existence.

\subsection{Energy estimates}
We begin with the standard energy estimates for perturbed wave
equations
\begin{equation}\label{2.1}
(\Box_\gamma u)^I=(\partial_t^2-c_I^2\Delta)u^I+\sum_{K=1}^D
\sum_{0\le j,k\le 3} \gamma^{IK,jk}\partial_j\partial_k u^K = G^I,
\quad I=1,\dots, D
\end{equation}
satisfying the symmetry condition
\begin{equation}\label{2.2}
\gamma^{IK,jk}=\gamma^{IK,kj}=\gamma^{KI,jk}, \quad 0\le j,k\le
3,\quad 1\le I,K\le D.
\end{equation}
As is standard, we let
$e_0=\sum_{I=1}^D e_0^I$ be the associated energy form where
\begin{multline}\label{2.3}
e_0^I(u,t)=(\partial_0u^I)^2 + \sum_{k=1}^3 c_I^2(\partial_k u^I)^2 +
2\sum_{J=1}^D\sum_{k=0}^3 \gamma^{IJ,0k}\partial_0 u^I \partial_k
u^J\\
-\sum_{J=1}^D \sum_{0\le j,k\le 3} \gamma^{IJ,jk}\partial_j
u^I\partial_k u^J.
\end{multline}
If we assume that
$$
\max_{\substack{1\le I,K\le D\\ 0\le j,k\le 3}}\|\gamma^{IK,jk}\|_\infty
$$
is sufficiently small, then it follows that
\begin{equation}\label{2.4}
\frac{1}{2}\sum_{1\le I\le D}\min(1,c_I^2)|\nabla_{t,x}u|^2\le
e_0(u)\le 2\sum_{1\le I\le D}\max(1,c_I^2)|\nabla_{t,x}u|^2.
\end{equation}
If we set $E(u,t)^2=\int_{\R^3} e_0(u,t)\:dx$ to be the associated
energy, then we have the energy inequality
\begin{multline}\label{2.5}
\sum_{|\alpha|\le M}\partial_t E(\Gamma^\alpha u,t)\le
C\sum_{|\alpha|\le M}\|\Gamma^\alpha G(t,\cd)\|_2 + \sum_{|\alpha|\le
M} \|[\Box_\gamma,\Gamma^\alpha] u(t,\cd)\|_2 \\
+C\sum_{|\alpha|\le M}E(\Gamma^\alpha u,t)\sum_{\substack{0\le
j,k,l\le 3\\ 1\le I,K\le D}} \|\partial_l \gamma^{IK,jk}(t,\cd)\|_\infty.
\end{multline}

In addition to the energy estimate \eqref{2.5}, we will need the
following $L^2_tL^2_x$ estimate of Keel-Smith-Sogge \cite{KSS2}
(Proposition 2.1).
\begin{lemma}\label{lemma2.1}
Suppose that $u\in C^\infty(\R\times\R^3)$ vanishes for large $x$ for
every $t$.  Then, there is a uniform constant $C$ so that
\begin{equation}\label{2.6}
(\log(2+t))^{-1/2}\|\langle x\rangle^{-1/2}
u'\|_{L^2([0,t]\times\R^3)} \le C\|u'(0,\cd)\|_2 + C\int_0^t \|\Box u(s,\cd)\|_2\:ds.
\end{equation}
\end{lemma}

\subsection{Pointwise estimates}
In this section, we will gather the pointwise estimates that will be
needed in the sequel.  The estimates that are presented are variants
of those in Keel-Smith-Sogge \cite{KSS3}, Sogge \cite{So2}, and
Kubota-Yokoyama \cite{KY}.  The key innovation in our approach to
Theorem \ref{theorem1.2} is the
use of both of these pointwise estimates and sharp Huygens' principle
to allow us to get good pointwise bounds for $u$ (not just $u'$ as in
\cite{MNS}).  This pointwise bound allows us to handle the higher
order terms without having to strengthen the null condition (as in
\cite{KY}).

In our first estimate, we will concentrate on the scalar wave equation
$\Box=(\partial_t^2-\Delta)$.  The transition to vector valued,
multiple speed wave equations is straightforward.

\begin{lemma}\label{lemma2.2}
Let $u$ be the solution of $\Box u(t,x)=F(t,x)$ with initial data
$u(0,\cd)=f$, $\partial_t u(0,\cd)=g$ for $(t,x)\in\R_+\times\R^3$.
Then,
\begin{multline}\label{2.7}
(1+t+|x|)|u(t,x)|\le C\sum_{|\alpha|\le 4} \|\langle
x\rangle^{|\alpha|}\partial^\alpha f\|_2 + C\sum_{|\alpha|\le 3}
\|\langle x\rangle^{1+|\alpha|} \partial^\alpha g\|_2\\
+C\sum_{\substack{\mu+|\alpha|\le 3\\ \mu\le 1}}\int_0^t \int_{\R^3}
|L^\mu Z^\alpha F(s,y)|\:\frac{dy\:ds}{\langle y\rangle}.
\end{multline}
\end{lemma}

\noindent{\em Proof of Lemma \ref{lemma2.2}:} For vanishing Cauchy
data, \eqref{2.7} can be found in Keel-Smith-Sogge \cite{KSS3} and
Sogge \cite{So2}.  Thus, it will suffice to show the estimate for
$\cos (t\sqrt{-\Delta}) f$ and $(\sin (t\sqrt{-\Delta}) /
\sqrt{-\Delta})g$.  The proof is similar to that in \cite{KSS3} for
the inhomogeneous case.  If we assume that $F=0$ above, we will show
\begin{multline}\label{2.8}
(1+t+|x|)|u(t,x)|\le C\sum_{\substack{|\alpha|+\mu\le 3\\\mu\le 1}}
\int_{\R^3} \Bigl(|(r\partial_r)^\mu Z^\alpha \nabla f|\:\frac{dy}{\langle
y\rangle} + \int_{\R^3} |(r\partial_r)^\mu Z^\alpha
f|\:\frac{dy}{\langle y\rangle^2}\\
+\int_{\R^3} |(r\partial_r)^\mu Z^\alpha g|\:\frac{dy}{\langle y\rangle}\Bigr)
\end{multline}
Our desired estimate \eqref{2.7} follows, then, via the Schwarz
inequality.

Let us first consider $(\sin (t\sqrt{-\Delta})/\sqrt{-\Delta})g$.
Using the positivity of the fundamental solution for the wave
equation, we have
\begin{equation}\label{2.9}
\begin{split}
|x|\Bigl|\frac{\sin (t\sqrt{-\Delta})}{\sqrt{-\Delta}} g\Bigr| &=
 \frac{t\,|x|}{4\pi} \Bigl|\int_{|\theta|=1} g(x-t\theta)\:d\sigma(\theta)\Bigr|\\
&\le \frac{1}{2}\int_{|t-|x||}^{t+|x|} \|sg(s\cd)\|_{L^\infty_\theta(S^2)}\:ds.
\end{split}
\end{equation}
By the embedding $H^{2,1}_\theta\hookrightarrow L^\infty_\theta$ it
follows that
\begin{equation}\label{2.10}
|x| \Bigl|\frac{\sin (t\sqrt{-\Delta})}{\sqrt{-\Delta}} g\Bigr| \le
 C\sum_{|\alpha|\le 2} \int_{|t-|x||\le |y|\le t+|x|} |\Omega^\alpha
 g(y)|\:\frac{dy}{|y|}.
\end{equation}

For $t\ge 10|x|$, apply the relation $sg(s\theta)=-\int_s^\infty
\partial_\tau (\tau g(\tau\theta))\:d\tau$ to \eqref{2.9} to see that
\begin{equation}\label{2.11}
\begin{split}
t\Bigl|\frac{\sin (t\sqrt{-\Delta})}{\sqrt{-\Delta}} g\Bigr| &\le
C\frac{t}{|x|}\int_{|t-|x||}^{t+|x|} \int_s^\infty \frac{1}{\tau} \|\tau\partial_\tau
(\tau g(\tau\theta))\|_{L^\infty_\theta(S^2)} \:d\tau\:ds\\
&\le C\int_{|t-|x||}^\infty \|\tau g(\tau\cd)\|_{L^\infty_\theta(S^2)} +
\|\tau(\tau\partial_\tau) g(\tau\cd)\|_{L^\infty_\theta(S^2)}
\:d\tau\\
&\le C\sum_{\mu\le 1,|\alpha|\le 2} \int_{|t-|x||\le |y|}
|(|y|\partial_{|y|})^\mu \Omega^\alpha g(y)|\:\frac{dy}{|y|}.
\end{split}
\end{equation}

By \eqref{2.10} and \eqref{2.11}, we obtain
\begin{equation}\label{2.12}
(t+|x|)\Bigl|\frac{\sin (t\sqrt{-\Delta})}{\sqrt{-\Delta}} g\Bigr| \le C\sum_{\mu\le
1, |\alpha|\le 2} \int_{|t-|x||\le |y|} |(|y|\partial_{|y|})^\mu
\Omega^\alpha g(y)|\:\frac{dy}{|y|}.
\end{equation}

We now wish to show that
\begin{equation}\label{2.13}
(1+t+|x|)\Bigl|\frac{\sin (t\sqrt{-\Delta})}{\sqrt{-\Delta}}g\Bigr|\le C
\sum_{\substack{|\alpha|+\mu\le 3\\\mu\le 1}} \int_{\R^3}
|(|y|\partial_{|y|})^\mu Z^\alpha g(y)|\:\frac{dy}{|y|}.
\end{equation}
For $t+|x|\ge 1$, \eqref{2.13} clearly follows from \eqref{2.12}.  For
$t+|x|\le 1$, let $\chi$ denote a smooth function with $\chi(x)\equiv 1$ for $|x|\le
1$ and $\chi(x)\equiv 0$ for $|x|>2$, and let $v$ be the solution to
the shifted wave equation
\begin{equation}\label{2.14}
\Box v(t,x)=0,\quad v(0,\cd)=0,\quad \partial_t v(0,x)=(\chi g)(x_1-10,x_2,x_3).
\end{equation}
By finite propagation, we have that $\frac{\sin
t\sqrt{-\Delta}}{\sqrt{-\Delta}} g = v(t,x_1+10,x_2,x_3)$ for
$t+|x|\le 1$, and \eqref{2.13} follows by applying \eqref{2.12} to
$v$.

Finally, we turn to the task of showing that our desired result
\begin{equation}\label{2.15}
(1+t+|x|)\Bigl|\frac{\sin (t\sqrt{-\Delta})}{\sqrt{-\Delta}}g\Bigr| \le
C\sum_{\substack{|\alpha|+\mu \le 3\\\mu\le 1}} \int_{\R^3}
|(|y|\partial_{|y|})^\mu Z^\alpha g(y)|\:\frac{dy}{\langle y\rangle}
\end{equation}
follows from \eqref{2.13}.  For $\chi$ as above, write $g=\chi g +
(1-\chi) g$.  When $g$ is replaced by $(1-\chi) g$, \eqref{2.15}
follows directly from \eqref{2.13}.  When $g$ is replaced by $\chi g$,
we instead apply \eqref{2.13} to the shifted function $v$.  It is this
use of the shifted function that introduces the translations on the
right sides of \eqref{2.13} and \eqref{2.15}.

Next, we consider $\cos (t\sqrt{-\Delta})f$.  We have
\begin{equation}\label{2.16}\begin{split}
\Bigl|\cos (t\sqrt{-\Delta}) f\Bigr| &= \Bigl|\partial_t
\Bigl(\frac{t}{4\pi}\int_{|\theta|=1}f(x-t\theta)\:d\sigma(\theta)\Bigr)\Bigr|\\
&\le \Bigl|\frac{\sin t\sqrt{-\Delta}}{t\sqrt{-\Delta}}f\Bigr|+\Bigl|\frac{\sin
t\sqrt{-\Delta}}{\sqrt{-\Delta}}|\nabla f|\Bigr|.
\end{split}
\end{equation}
By \eqref{2.15}, the $|\nabla f|$ part is bounded by the right side of
\eqref{2.8}.  For the first part, repeating the arguments of
\eqref{2.9} and \eqref{2.11}, we have
\begin{equation}\label{2.17}
\begin{split}
(t+|x|)\Bigl|\frac{\sin (t\sqrt{-\Delta})}{t\sqrt{-\Delta}}f\Bigr|& \le
\frac{t+|x|}{2t|x|}\int_{|t-|x||}^{t+|x|}
\|sf(s\cd)\|_{L^\infty_\theta(S^2)}\:ds\\
&\le \frac{t+|x|}{2t|x|}(t+|x|-|t-|x||)\int_{|t-|x||}^\infty
\|\partial_\tau(\tau f(\tau\cd))\|_{L^\infty_\theta(S^2)}\:d\tau\\
&\le C\sum_{\mu\le 1,|\alpha|\le 2} \int_{|y|\ge |t-|x||}
|(|y|\partial_{|y|})^\mu \Omega^\alpha f(y)|\:\frac{dy}{|y|^2}.
\end{split}
\end{equation}
Using the shifted function as in \eqref{2.13} and \eqref{2.15}, it
follows that
\begin{equation}\label{2.18}
(1+t+|x|)\Bigl|\frac{\sin (t\sqrt{-\Delta})}{t\sqrt{-\Delta}}f\Bigr|\le C
\sum_{\substack{|\alpha|+\mu\le 3\\\mu\le 1}} \int
|(|y|\partial_{|y|})^\mu Z^\alpha f(y)|\:\frac{dy}{\langle y\rangle^2}
\end{equation}
as desired.\qed

We now wish to explore the version of the pointwise estimate of
Kubota-Yokoyama \cite{KY} that we will use.  
We define the ``neighborhoods'' of the characteristic cones
$r=|x|=c_I t$ for $\Box_{c_I}$.  That is, with the $c_I$ as in
\eqref{order}, 
set
\begin{equation}\label{2.19}
\Lambda_I=\{(t,|x|)\in [1,\infty)\times [1,\infty)\, : \, |r-c_I t|\le
\delta t\}
\end{equation}
where $\delta=\frac{1}{3}\min_{1\le I\le D} (c_I-c_{I-1})$ and
$I=1,2,\dots, D$.  Note that
for $(t,x)\not\in \Lambda_I$, $|c_I t-|x||\approx t+|x|$.
Additionally, define
\begin{equation}\label{2.20}
z(s,\lambda)=\begin{cases} (1+|\lambda-c_J s|),&
\text{ for } (s,\lambda)\in\Lambda_J,\quad J=1,2,\dots,D\\
(1+\lambda),& \text{ otherwise.}\end{cases}
\end{equation}

With this notation, we then have
\begin{lemma}\label{lemma2.3}
Let $I=1,2,\dots, D$, and assume that $G^I(t,x)$ is a continuous
function of $(t,x)\in \R_+\times\R^3$.  Let $w^I$ be the solution of
$(\partial_t^2-c_I^2\Delta)w^I=G^I$ with vanishing Cauchy data at time $t=0$.
Then,
\begin{multline}\label{2.21}
(1+r+t)\Bigl(1+\log \frac{1+r+c_It}{1+|r-c_It|}\Bigr)^{-1}
|w^I(t,x)|\\\le C
\sup_{(s,y)\in D^I(t,r)}|y|(1+s+|y|)^{1+\mu}z^{1-\mu}(s,|y|) |G^I(s,y)|
\end{multline}
for any $\mu>0$ and 
\begin{equation}\label{2.22}
D^I(t,r)=\{(s,y)\in \R\times\R^3\,:\, 0\le s< t,\, |r-c_I(t-s)|\le
|y|\le r+c_I(t-s)\}.
\end{equation}
\end{lemma}

The above estimate is due to Kubota-Yokoyama \cite{KY} (Theorem 3.4).
If we combine \eqref{2.7} and \eqref{2.21} and use the fact that
$[\Box, Z]=0$ and $[\Box, L]=2\Box$, we get our main pointwise
estimates. 
\begin{theorem}\label{theorem2.4a}
Let $I=1,2,\dots, D$, and assume that $F^I(t,x)$, $G^I(t,x)$ are
smooth functions of $(t,x)\in\R_+\times\R^3$.  Let $w^I$ be the
solution of $(\partial_t^2-c_I^2\Delta)w^I=F^I+G^I$.  Then, there is a
uniform constant $\tilde{C}_1>0$ so that
\begin{multline}\label{2.23a}
(1+r+t)|\Gamma^\alpha w^I(t,x)|\le \tilde{C}_1 \sum_{|\beta|\le 4}
\|\langle x\rangle^{|\beta|} (\partial^\beta \Gamma^\alpha
w^I)(0,\cd)\|_2
\\+\tilde{C}_1 \sum_{|\beta|\le 3} \|\langle x\rangle^{1+|\beta|}
(\partial^\beta \Gamma^\alpha \partial_t w^I)(0,\cd)\|_2 + \tilde{C}_1
\sum_{|\beta|\le |\alpha|+3} \int_0^t \int |\Gamma^\beta
F^I(s,y)|\:\frac{dy\:ds}{\langle y\rangle}
\\+\tilde{C}_1 \Bigl(1+\log\frac{1+r+c_It}{1+|r-c_It|}\Bigr)
\sum_{|\beta|\le |\alpha|}\sup_{(s,y)\in D^I(t,r)}
|y|(1+s+|y|)^{1+\mu} z^{1-\mu}(s,|y|)|\Gamma^\beta G^I(s,y)|
\end{multline}
for any multiindex $\alpha$, $\mu>0$, and $D^I$ as in \eqref{2.22}.
\end{theorem}
Using strong Huygens' principle, we can establish the following
variant of the previous theorem.
\begin{theorem}\label{theorem2.4}
Fix $I=1,2,\dots, D$, and assume that $F^I(t,x), G^I(t,x)$ are
smooth functions of $(t,x)\in\R_+\times\R^3$.  Moreover, assume
that $F^I(t,x)$ is supported in $\Lambda_J$ for some $J\neq I$.
Let $w^I$ be the
solution of $(\partial_t^2-c_I^2\Delta)w^I=F^I+G^I$.
Then, there are uniform
constants $c, C_1>0$ depending on the wavespeeds $c_I,c_J$ so that
\begin{multline}\label{2.23}
(1+r+t)|\Gamma^\alpha w^I(t,x)|\le 
C_1\sum_{|\beta|\le 4}\|\langle x\rangle^{|\beta|}(\partial^\beta\Gamma^\alpha w^I)(0,\cd)\|_2 \\+
C_1\sum_{|\beta|\le 3}\|\langle x\rangle^{1+|\beta|}
(\partial^\beta\Gamma^\alpha \partial_t w^I)(0,\cd)\|_2 
+C_1 \sum_{|\beta|\le |\alpha|+2} \sup_{0\le s\le t} \int
|\Gamma^\beta F^I(s,y)|\:dy
\\+
 C_1 \sum_{|\beta|\le |\alpha|+3} \int_{\max(0,c|c_It-|x||-1)}^t
\int_{|y|\approx s}
|\Gamma^\beta F^I(s,y)|\: \frac{dy\:ds}{\langle y\rangle} \\+ C_1\Bigl(1+\log
\frac{1+r+c_It}{1+|r-c_It|}\Bigr) \sum_{|\beta|\le |\alpha|}\sup_{(s,y)\in D^I(t,r)}
|y|(1+s+|y|)^{1+\mu}z^{1-\mu}(s,|y|)|\Gamma^\beta G^I(s,y)|
\end{multline}
for any multiindex $\alpha$,  $\mu>0$ and $D^I$ as in \eqref{2.22}.
\end{theorem}

\noindent Here, and throughout, $|y|\approx s$ is used to denote that there is a
positive constant $\tilde{c}$ so that $(1/\tilde{c})|y|\le s\le \tilde{c}|y|$.

\noindent{\em Proof of Theorem \ref{theorem2.4}:}  
By \eqref{2.21}, we may take $G^I\equiv 0$ without restricting
generality.  We then note that there is a constant $c$ so that the
intersection of the backward light cone through $(t,x)$ with speed
$c_I$, $\{c_I(t-s)=|x-y|\}$, and $\Lambda_J$ is contained in
$[c|c_It-|x||,t]\times \{|y|\approx s\}$.  With this in mind, we fix a
smooth cutoff function $\rho$ so that $\rho(s)\equiv 1$ for $s\ge
c|c_It-|x||$ and $\rho(s)\equiv 0$ for $s\le c|c_It-|x||-1$.  Notice
that by strong Huygens' principle, we have $\Gamma^\alpha
w^I(t,x)=\Gamma^\alpha \tilde{w}$ where $\tilde{w}$ is the solution to
$$\Box_{c_I}\Gamma^\alpha \tilde{w}(s,y)=\rho(s)\Gamma^\alpha
F^I(s,y)+\rho(s)[\Box_{c_I}, \Gamma^\alpha] F^I(s,y)$$
and $\Gamma^\alpha \tilde{w}$ has the same Cauchy data as
$\Gamma^\alpha w$.

The result now follows from an application of \eqref{2.7} 
to $\Gamma^\alpha \tilde{w}$.  So long as the scaling
vector field $L$ in the third term on the right of \eqref{2.7} does
not hit $\rho$, the bound \eqref{2.23} follows and the third term on
the right is unnecessary.  If the $L$ in \eqref{2.7} is applied to
$\rho$, we get an additional term which is bounded by
$$C\sum_{|\beta|\le|\alpha|+2} \int_{\max(0,c|c_It-|x||-1)}^{c|c_It-|x||} 
\int_{|y|\approx s} s|\rho'(s)||\Gamma^\beta
F(s,y)|\:\frac{dy\:ds}{\langle y\rangle}.$$ 
Since $|y|\approx s$ and the time integral is taken over an interval
of length at most one, this term is easily seen to be dominated by the
third term in \eqref{2.23} which completes the proof. \qed

\subsection{Null form estimates and Sobolev-type estimates}
In this section, we gather our bounds on the null forms and some
weighted Sobolev-type estimates.  The first of these is the null form
estimate.  See, e.g., \cite{Si3}, \cite{So2}.

\begin{lemma}\label{lemma2.5}
Suppose that the quadratic parts of the nonlinearity $Q(du,d^2u)$,
$B(du)$ satisfy the null conditions \eqref{1.8} and \eqref{1.9}.
Then,
\begin{equation}\label{2.24}
\Bigl|\sum_{0\le j,k,l\le 3} B^{KK,jk}_{K,l}\partial_l 
u\partial_j\partial_k v\Bigr|\le C\langle
 r\rangle^{-1} (|\Gamma u||\partial^2 v|+|\partial u||\partial \Gamma
 v|)+ C\frac{\langle c_K t-r\rangle}{\langle t+r\rangle} |\partial
 u||\partial^2 v|,
\end{equation}
and
\begin{equation}\label{2.25}
\Bigl|\sum_{0\le j,k\le 3}A^{K,jk}_{KK} \partial_j 
u\partial_k v\Bigr|\le C\langle r\rangle^{-1}
 (|\Gamma u||\partial v|+|\partial u||\Gamma v|)+ C\frac{\langle
 c_Kt-r\rangle}{\langle t+r\rangle}|\partial u||\partial v|.
\end{equation}
\end{lemma}

For the Sobolev-type results, we begin with
\begin{lemma}\label{lemma2.6}
Suppose that $h\in C^\infty(\R^3)$.  Then, for $R>1$,
\begin{equation}\label{2.26}
\|h\|_{L^\infty(R/2<|x|<R)} \le C R^{-1}\sum_{|\alpha|+|\beta|\le 2}
\|\Omega^\alpha \partial_x^\beta h\|_{L^2(R/4<|x|<2R)}.
\end{equation}
\end{lemma}
\noindent This has become a rather standard result.  See Klainerman \cite{K}.  
A proof can also be found, e.g., in \cite{KSS2}.

Additionally, we have the following space-time weighted Sobolev
results.
\begin{lemma}\label{lemma2.7}
Let $u\in C^\infty_0(\R_+\times\R^3)$.  Then,
\begin{equation}\label{2.27}
\langle r\rangle^{1/2}|u(t,x)|\le C\sum_{|\alpha|\le 1}\|Z^\alpha u'(t,\cd)\|_2,
\end{equation}
\begin{equation}\label{2.28}
\|\langle c_I t-r\rangle 
\partial^2 u(t,\cd)\|_2\le C\sum_{|\beta|\le 1}\|\Gamma^\beta
u'(t,\cd)\|_2 + C\|\langle t+r\rangle \Box_{c_I} u(t,\cd)\|_2,
\end{equation}
\begin{equation}\label{2.29}
\langle r\rangle^{1/2}\langle c_I t-r\rangle | u'(t,x)|\le
C\sum_{|\beta|\le 1} \|Z^\beta u'(t,\cd)\|_2 +
C\sum_{|\beta|\le 1} \|\langle c_I t-r\rangle Z^\beta
\partial^2 u(t,\cd)\|_2,
\end{equation}
\begin{equation}\label{2.30}
\langle r\rangle \langle c_I t-r\rangle^{1/2}|u'(t,x)|\le
C\sum_{|\beta|\le 2} \|Z^\beta u'(t,\cd)\|_2 +
C\sum_{|\beta|\le 1} \|\langle c_I t-r\rangle Z^\beta
\partial^2 u(t,\cd)\|_2.
\end{equation}
\end{lemma}

The estimates \eqref{2.27} and \eqref{2.30} are shown in Sideris \cite{Si} 
(Proposition 3.3).
\eqref{2.28} is due to Klainerman-Sideris \cite{KS} (Lemma 2.3 and Lemma 3.1).
\eqref{2.29} is from Hidano-Yokoyama \cite{HY} (Lemma 4.1) and follows 
from \eqref{2.27}.


Lastly, by interpolating between \eqref{2.29} and
\eqref{2.30}, it is easy to see that
\begin{multline}\label{2.31}
\langle r\rangle^{1/2+\mu}\langle c_I t-r\rangle^{1-\mu} |\Gamma^\alpha
u'(t,x)|\le C \sum_{|\beta|\le |\alpha|+2}\|\Gamma^\beta u'(t,\cd)\|_2
\\ + C\sum_{|\beta|\le |\alpha|+1}\|\langle c_I t-r\rangle \Gamma^\beta
\partial^2 u(t,\cd)\|_2
\end{multline}
for any $0\le \mu\le 1/2$.

\newsection{Global existence in Minkowski space}
Here we prove Theorem \ref{theorem1.2}.  We will take $N=71$ in
\eqref{1.13}.  This, however, is not optimal.

To proceed, we shall require a standard local existence theorem.
\begin{theorem}\label{theorem3.1}
Let $f\in H^{71}(\R^3)$ and $g\in H^{70}(\R^3)$.  Then, there is a
$T>0$ dependent on the norm of the data so that the initial value
problem \eqref{1.12} has a $C^2$ solution satisfying
\begin{equation}\label{3.1}
u\in L^\infty([0,T];H^{71}(\R^3))\cap C^{0,1}([0,T];H^{70}(\R^3)).
\end{equation}
The supremum of all such $T$ is equal to the supremum of all $T$ such
that the initial value problem has a $C^2$ solution with
$\partial^\alpha u$ bounded for all $|\alpha|\le 2$.
\end{theorem}

This result is a multi-speed analog of Theorem 6.4.11 in \cite{H}
(which is stated only for scalar wave equations).  Since the proof is
based only on energy inequalities, the same argument yields Theorem
\ref{theorem3.1} provided we assume the symmetry conditions
\eqref{1.6} and \eqref{1.7}.

We are now ready to set up our continuity argument.  If $\varepsilon$
is as above, we will assume that we have a solution of our
equation \eqref{1.12} for $0\le t\le T$ satisfying the following:
\begin{align}
\sum_{|\alpha|\le 50}\|\Gamma^\alpha u'(t,\cd)\|_2 &\le A_0\varepsilon
\label{3.2}\\
(1+t+|x|)\sum_{|\alpha|\le 40}|\Gamma^\alpha u^I(t,x)|&\le
A_1\varepsilon\Bigl(1+\log\frac{1+t+|x|}{1+|c_It-|x||}\Bigr)
\label{3.3}\\
(1+t+|x|)\sum_{|\alpha|\le 60} |\Gamma^\alpha u^I(t,x)|&\le
A_2\varepsilon (1+t)^{1/10}
\log(2+t)\Bigl(1+\log\frac{1+t+|x|}{1+|c_It-|x||}\Bigr)\label{3.4}\\
(1+t+|x|)\sum_{|\alpha|\le 39}|\Gamma^\alpha u'(t,x)|&\le
B_1\varepsilon\label{3.5}\\
\sum_{|\alpha|\le 70}\|\Gamma^\alpha u'(t,\cd)\|_2&\le B_2\varepsilon
(1+t)^{1/40}\label{3.6}\\
\sum_{|\alpha|\le 65}\|\langle x\rangle^{-1/2} \Gamma^\alpha
u'\|_{L^2(S_t)}&\le B_3 \varepsilon (1+t)^{1/20}(\log(2+t))^{1/2}.\label{3.7}
\end{align}
Here $S_t$ denotes the time strip $[0,t]\times\R^3$.

By \eqref{1.13}, we have the estimate
\begin{multline*}
\sum_{I=1}^D \sum_{|\alpha|\le 67}(1+C_1+\tilde{C}_1)\{
\sum_{|\beta|\le 4}\|\langle x\rangle^{|\beta|}
(\partial^\beta\Gamma^\alpha u^I)(0,x)\|_2
\\+ \sum_{|\beta|\le 3}\|\langle x\rangle^{1+|\beta|}
(\partial^\beta\Gamma^\alpha \partial_t u^I)(0,x)\|_2
\} \le C_2\varepsilon 
\end{multline*}
for some constant $C_2>0$.
Here $\tidle{C_1}$ and $C_1$ are the constants 
occurring in \eqref{2.23a} and \eqref{2.23} respectively.  
In our estimates above, we choose $A_0=A_1=A_2=A\ge 10 \max(1,C_2)$.

We shall then prove that for
$\varepsilon$ sufficiently small,
\begin{enumerate}
\item[i.)] \eqref{3.2} holds with $A_0$ replaced by $A_0/2$.
\item[ii.)] \eqref{3.3}, \eqref{3.4}
hold with $A_1,A_2$ replaced by $A_1/2, A_2/2$ respectively.
\item[iii.)] \eqref{3.2}-\eqref{3.4} imply \eqref{3.5}-\eqref{3.7} for
a suitable choice of constants $B_1,B_2,B_3$.
\end{enumerate}
We will prove items (i.)-(iii.) in the next three subsections
respectively.

Before we begin with the proof of (i.), we will set up some
preliminary results under the assumption of \eqref{3.2}-\eqref{3.7}.  Let us first prove
\begin{equation}\label{3.8}
\sum_{|\alpha|\le 58} \langle r\rangle^{1/2+\mu}\langle
c_I t-r\rangle^{1-\mu}|\Gamma^\alpha \partial u^I (t,x)|\le C\varepsilon
(1+t)^{1/40},\quad 0\le \mu\le 1/2.
\end{equation}
Indeed, by \eqref{2.31} and \eqref{2.28}, we have that the left side
of \eqref{3.8} is controlled by
$$C\sum_{|\alpha|\le  60} \|\Gamma^\alpha u'(t,\cd)\|_2 +
C\sum_{|\alpha|\le 59} \|\langle t+r\rangle \Gamma^\alpha \Box_{c_I}
u^I(t,\cd)\|_2.$$
By \eqref{3.6}, the first term is controlled by the right side of
\eqref{3.8}.  Thus, it remains to show
\begin{equation}\label{3.9}
\sum_{|\alpha|\le 59} \|\langle t+r\rangle \Gamma^\alpha \Box_{c_I}
u^I(t,\cd)\|_2 \le C\varepsilon (1+t)^{1/40}.
\end{equation}
By our definition of $\Box u$, we have that the left side of \eqref{3.9}
is bounded by
\begin{multline*}
C\sum_{|\alpha|\le 30} \|\langle t+r\rangle \Gamma^\alpha
u'(t,\cd)\|_\infty \sum_{|\alpha|\le 60} \|\Gamma^\alpha u'(t,\cd)\|_2
\\+C\sum_{1\le J,K\le D} \Bigl\|\langle t+r\rangle \sum_{|\alpha|\le
31} |\Gamma^\alpha u^J| \sum_{|\alpha|\le 31} |\Gamma^\alpha u^K|
\sum_{|\alpha|\le 59} |\Gamma^\alpha u|\Bigr\|_2
\\+C\sum_{1\le J,K\le D} \Bigl\|\langle t+r\rangle \sum_{|\alpha|\le
31} |\Gamma^\alpha u^J|\sum_{|\alpha|\le 31}|\Gamma^\alpha u^K|
\sum_{|\alpha|\le 60} |\Gamma^\alpha \partial u|\Bigr\|_2.
\end{multline*}
By \eqref{3.5} and \eqref{3.6}, we see that the first term is
controlled by $C\varepsilon^2 (1+t)^{1/40}$ as desired.  For the second
term, we apply \eqref{3.3} to see that we have the bound
$$
C\varepsilon^2
\Bigl\|\Bigl(1+\log\frac{1+t+|x|}{1+|c_Jt-|x||}\Bigr)\Bigl(1+\log\frac{1+t+|x|}{1+|c_Kt-|x||}\Bigr)
(1+t+|x|)^{-1} \sum_{|\alpha|\le 59} |\Gamma^\alpha u(t,\cd)|\Bigr\|_2.
$$
We, then, see that this is $O(\varepsilon^3)$
using \eqref{3.4}.  The bound for the third term follows similarly
from applications of \eqref{3.3} and \eqref{3.6}.

If we argued similarly, using \eqref{3.2} instead of \eqref{3.6}, it
follows that
\begin{equation}\label{3.10}
\sum_{|\alpha|\le 48} \langle r\rangle^{1/2+\mu} \langle
c_I t-r\rangle^{1-\mu} |\Gamma^\alpha \partial u^I(t,x)|\le C\varepsilon, \quad 0\le
\mu\le 1/2,
\end{equation}
and
\begin{equation}\label{3.11}
\sum_{|\alpha|\le 49} \|\langle c_I t-r\rangle \Gamma^\alpha \partial^2
u^I(t,\cd)\|_2\le C\varepsilon.
\end{equation}
Indeed, the latter follows from \eqref{2.28} and the proof of
\eqref{3.9} where, as mentioned above, we use the lossless estimate
\eqref{3.2} rather than \eqref{3.6}.

\subsection{Proof of (i.):}  In this section, we will show that
\eqref{3.2}-\eqref{3.7} allow you to prove \eqref{3.2} with $A_0$ 
replaced by $A_0/2$.  By the standard energy inequality (see, e.g.,
\cite{S}), the square of the left side of \eqref{3.2} is controlled by
\begin{equation}\label{3.12}
\sum_{|\alpha|\le 50} \|
\Gamma^\alpha u'(0,\cd)\|^2_2 +
\sum_{|\alpha|\le 50} \int_0^t \int \Bigl|\langle
\partial_0 \Gamma^\alpha u,\Box \Gamma^\alpha u\rangle\Bigr|\:dy\:ds.
\end{equation}
It follows from \eqref{1.13} and our choice of $A_0$ that the first term
is controlled by $(A_0/10)^2 \varepsilon^2$.  Thus, it will suffice to show that
\begin{equation}\label{3.13}
\sum_{|\alpha|\le 50}\int_0^t \int \Bigl|\langle \partial_0
\Gamma^\alpha u, \Box \Gamma^\alpha u\rangle\Bigr|\:dy\:ds \le C\varepsilon^3.
\end{equation}

The left side of \eqref{3.13} is dominated by
\begin{multline}\label{3.14}
C\int_0^t \int_{\R^3} \sum_{K=1}^D \sum_{|\alpha|\le 50} |\partial_0
\Gamma^\alpha u^K| \sum_{|\alpha|+|\beta|\le 50} \Bigl|\sum_{0\le
j,k,l\le 3} \tilde{B}^{KK,jk}_{K,l} \partial_l \Gamma^\alpha u^K
\partial_j\partial_k \Gamma^\beta u^K\Bigr|\:dy\:ds
\\+C \int_0^t\int_{\R^3} \sum_{K=1}^D \sum_{|\alpha|\le 50}
|\partial_0 \Gamma^\alpha u^K|\sum_{|\alpha|+|\beta|\le 50}
\Bigl|\sum_{0\le j,k,l\le 3} \tilde{A}^{K,jk}_{KK} \partial_j
\Gamma^\alpha u^K \partial_k \Gamma^\beta u^K\Bigr|\:dy\:ds
\\+ C\int_0^t \int_{\R^3} \sum_{\substack{1\le I,J,K\le D\\ (I,K)\neq
(K,J)}} \sum_{|\alpha|\le 50} |\partial \Gamma^\alpha u^K|
\sum_{|\alpha|\le 50} |\partial\Gamma^\alpha u^I|\sum_{|\alpha|\le
51}|\partial \Gamma^\alpha u^J|\:dy\:ds
\\+ C\int_0^t \int_{\R^3} \sum_{|\alpha|\le 50} |\partial_0
\Gamma^\alpha u| \Bigl(\sum_{|\alpha|\le 31} |\Gamma^\alpha u|\Bigr)^2
\sum_{|\alpha|\le 52} |\Gamma^\alpha u|\:dy\:ds.
\end{multline}
Due to constants that are introduced when $L^\nu Z^\alpha$ commutes
with $\partial_{j,k,l}$, the coefficients $A^{K,jk}_{KK}$,
$B^{KK,jk}_{K,l}$ become new constants $\tilde{A}^{K,jk}_{KK}$,
$\tilde{B}^{KK,jk}_{K,l}$.  It is known, however, that $\Gamma$ preserves the
null forms.  That is, since the original constants satisfy
\eqref{1.8} and \eqref{1.9}, so do the new ones
$\tilde{A}^{K,jk}_{KK}$ and $\tilde{B}^{KK,jk}_{K,l}$.  See, e.g.,
Sideris-Tu \cite{Si3} (Lemma 4.1).

The first three terms are handled as in \cite{MNS}.  Let us begin with
the null terms (i.e., the first two terms in \eqref{3.14}).  By
\eqref{2.24} and \eqref{2.25}, these terms are dominated by
\begin{multline}\label{3.15}
C\int_0^t \int_{\R^3} \sum_{|\alpha|\le 50} |\Gamma^\alpha u'|
\sum_{|\alpha|\le 51} |\Gamma^\alpha u| \sum_{|\alpha|\le 51}
|\Gamma^\alpha u'|\:\frac{dy\:ds}{\langle y\rangle}
\\+C\int_0^t\int_{\R^3} \sum_{K=1}^D \frac{\langle c_K
s-r\rangle}{\langle s+r\rangle} \Bigl(\sum_{|\alpha|\le 51}
|\Gamma^\alpha \partial u^K|\Bigr)^3\:dy\:ds 
\end{multline}

In order to handle the contribution by the first term
of \eqref{3.15}, notice that by \eqref{3.4}
$$\sum_{|\alpha|\le 51} |\Gamma^\alpha u(s,y)|\le C\varepsilon \langle
s+|y|\rangle^{-9/10+}.$$
Thus, the first term in
\eqref{3.15} has a contribution to \eqref{3.14} which is dominated by
\begin{equation}\label{3.16}
C\varepsilon \int_0^t \langle s\rangle^{-9/10+}
\sum_{|\alpha|\le 51} \|\langle y\rangle^{-1/2} \Gamma^\alpha u'(s,\cd)\|^2_2\:ds
\end{equation}
by the Schwarz inequality.  By \eqref{3.7}, it follows that this
contribution is $O(\varepsilon^3)$.

In order to show that the second term in \eqref{3.15} satisfies a
similar bound, we apply \eqref{3.8} with $\mu=0$ and the Schwarz inequality 
 to see that it is controlled by
\begin{multline}\label{3.17}
C\varepsilon\int_0^t (1+s)^{1/40} \int_{\R^3} \frac{1}{\langle
r\rangle^{1/2} \langle s+r\rangle} \sum_{|\alpha|\le 51}
|\Gamma^\alpha \partial u|^2\:dy\:ds \\
\le C \int_0^t \langle s\rangle^{-19/40} \sum_{|\alpha|\le 51}
\|\langle y\rangle^{-1/2}\, \Gamma^\alpha u'(s,\cd)\|^2_2\:ds.
\end{multline}
It then follows from \eqref{3.7} that this term also has an
$O(\varepsilon^3)$ contribution to \eqref{3.14}.

We now wish to show that the multi-speed terms
\begin{equation}\label{3.18}
\int_0^t \int_{\R^3} \sum_{|\alpha|\le 50}|\partial \Gamma^\alpha
u^K|\sum_{|\alpha|\le 50} |\partial \Gamma^\alpha
u^I|\sum_{|\alpha|\le 51} |\partial \Gamma^\alpha u^J|\:dy\:ds
\end{equation}
with $(I,K)\neq (K,J)$ have an $O(\varepsilon^3)$ contribution to
\eqref{3.14}.  For simplicity, let us assume that $I\neq K$, $I=J$.  A
symmetric argument will yield the same bound for the remaining cases.
If we set $\delta < |c_I-c_K|/3$, it follows that $\{|y|\in
[(c_I-\delta) s, (c_I+\delta) s]\}\cap \{|y|\in [(c_K-\delta) s,
(c_K+\delta)s]\}=\emptyset$.  Thus, it will suffice to show the bound
when the spatial integral is taken over the complements of each of
these sets separately.  We will show the bound over $\{|y|\not\in
[(c_K-\delta) s, (c_K+\delta)s]\}$.  The same argument will
symmetrically yield the bound over the other set.

If we apply \eqref{3.8} with $\mu=0$, we see that over the indicated
set, \eqref{3.18} is bounded by
\begin{multline}\label{3.19}
C\varepsilon \int_0^t \int_{\{|y|\not\in
[(c_K-\delta)s,(c_K+\delta)s]\}} \langle s+r\rangle^{-39/40} \langle
r\rangle^{-1/2} \sum_{|\alpha|\le 51} |\partial \Gamma^\alpha
u^I|^2\:dy\:ds
\\\le C\varepsilon \int_0^t \langle s\rangle^{-19/40}
\sum_{|\alpha|\le 51} \|\langle y\rangle^{-1/2} \Gamma^\alpha
u'(s,\cd)\|^2_2 \:ds.
\end{multline}
Thus, it again follows from \eqref{3.7} that this term is $O(\varepsilon^3)$.

Finally, it remains to bound the contribution to \eqref{3.14} by the
cubic terms (the fourth term in \eqref{3.14}).  If we apply
\eqref{3.3} and \eqref{3.4}, it is clear that this term is dominated
by
$$C\varepsilon^3 \int_0^t \int
\frac{(\log(1+s+|y|))^4}{(1+s+|y|)^{29/10}}\sum_{|\alpha|\le 50}
|\partial_0 \Gamma^\alpha u|\:dy\:ds.$$
By Schwarz inequality and \eqref{3.6}, we see that this term is
$O(\varepsilon^4)$ which completes the proof of \eqref{3.13}.

\subsection{Proof of (ii.):} 
In this section, we wish to show that our pointwise estimates
\eqref{3.3} and \eqref{3.4} hold with $A_1,A_2$ replaced by $A_1/2,
A_2/2$ respectively.  Let us begin with \eqref{3.3}.

Fix a smooth cutoff function $\eta_J$ satisfying $\eta_J(s)\equiv 1$, 
$s\in [(c_J+(\delta/2))^{-1}, (c_J-(\delta/2))^{-1}]$ where, as in \eqref{2.19},
$\delta=(1/3)\min_{I}(c_I-c_{I-1})$, and $\eta_J(s)\equiv 0$, $s\not\in
[(c_J+\delta)^{-1}, (c_J-\delta)^{-1}]$.  We also set $\beta$ to be a
smooth function satisfying $\beta(x)\equiv 1$, $|x|<1$ and
$\beta(x)\equiv 0$, $|x|\ge 2$.  Then, let $\rho_J(x,t)=(1-\beta)(x)
\eta_J(|x|^{-1}t)$.  By construction when $|x|\ge 2$, $\rho_J$ is
identically $1$ in a conic neighborhood of $\{c_Jt=|x|\}$
and is supported on $\Lambda_J$.

We then set
\begin{equation}\label{3.20}
\tilde{F}^I=\sum_{\substack{1\le J\le D\\ J\neq I}} \sum_{0\le
j,k,l\le 3} B^{IJ,jk}_{J,l}\rho_J \partial_l u^J \partial_j\partial_k
u^J + \sum_{\substack{1\le J\le D\\J\neq I}} \sum_{0\le j,k\le 3}
A^{I,jk}_{JJ} \rho_J \partial_j u^J\partial_k u^J
\end{equation}
and $\tilde{G}^I=F^I-\tilde{F}^I$.  By \eqref{2.23} and our choice of $C_2$, we have that
the left side of \eqref{3.3} is dominated by
\begin{multline}\label{3.21}
C_2\varepsilon 
+C\Bigl(1+\log\frac{1+r+c_It}{1+|r-c_It|}\Bigr)\sum_{|\beta|\le 40}
\sup_{(s,y)\in D^I(t,r)} |y|(1+s+|y|)^{1+\mu}
z^{1-\mu}(s,|y|)|\Gamma^\beta \tilde{G}^I(s,y)|
\\+C\sum_{|\beta|\le 43} \int_{\max(0,c|c_It-|x||-1)}^t \int_{|y|\approx s} |\Gamma^\beta
\tilde{F}^I(s,y)|\:\frac{dy\:ds}{\langle y\rangle}
+C\sum_{|\beta|\le 42} \sup_{0\le s\le t} \int |\Gamma^\beta \tilde{F}^I(s,y)|\:dy.
\end{multline}
By construction, we have $C_2\varepsilon \le (A_1/10)\varepsilon$.

We now turn to the second to last term in \eqref{3.21}.  
Since
$|y|\approx s$ on the support of $\rho_J$, it follows that this term is
controlled by
\begin{multline*}
C\sum_{\substack{1\le J\le D\\J\neq I}}\int_{\max(0,c|c_It-|x||-1)}^t \frac{1}{1+s}\int_{|y|\approx s}
\sum_{|\beta|\le 44} |\Gamma^\beta \partial u^J|^2\:dy\:ds
\\\le C\Bigl(1+\log\frac{1+t}{1+|c_It-|x||}\Bigr)\sup_{0\le s\le t}
\sum_{|\beta|\le 44} \|\Gamma^\beta u'(s,\cd)\|_2^2.
\end{multline*}
The correct bound for the right side then follows from \eqref{3.2}.
If we apply the Schwarz inequality, it follows that the last term in
\eqref{3.21} is dominated by
$$C\sum_{|\beta|\le 43} \sup_{0\le s\le t} \|\Gamma^\beta u'(s,\cd)\|^2_2.$$ 
Thus, by \eqref{3.2}, 
we get the desired bound for the $\tilde{F}^I$
terms in \eqref{3.21}.

It remains to examine the $\tilde{G}^I$ term in \eqref{3.21}.  
The proof of \eqref{3.3} will be complete if we can show that
\begin{equation}\label{3.22}
\sum_{|\beta|\le 40}
\sup_{(s,y)\in D^I(t,r)} |y|(1+s+|y|)^{1+\mu}
z^{1-\mu}(s,|y|)|\Gamma^\beta \tilde{G}^I(s,y)|\le C\varepsilon^2.
\end{equation}

When
$\tilde{G}^I$ is replaced by the null forms
$$\sum_{0\le j,k\le 3} A^{I,jk}_{II}\partial_j u^I \partial_k u^I
+ \sum_{0\le j,k,l\le 3} B^{II,jk}_{I,l}\partial_l u^I
\partial_j\partial_k u^I,$$
we apply \eqref{2.24} and \eqref{2.25} to bound this term by
\begin{multline}\label{3.23}
C\sup_{(s,y)\in D^I(t,r)} (1+s+|y|)^{1+\mu} z^{1-\mu}(s,|y|)
\sum_{|\beta|\le 41}|\Gamma^\beta
u^I|\sum_{|\beta|\le 41}|\Gamma^\beta \partial u^I|\\
+C\sup_{(s,y)\in D^I(t,r)} |y|(1+s+|y|)^{\mu} z^{1-\mu}(s,|y|)
  \langle c_Is-|y|\rangle\sum_{|\beta|\le
21} |\Gamma^\beta \partial u^I| \sum_{|\beta|\le 41}|\Gamma^\beta
\partial u^I|.
\end{multline}
For the first term in \eqref{3.23}, if we apply \eqref{3.4}, we see
that it is controlled by
$$C\varepsilon \sup_{(s,y)\in D^I(t,r)} (1+s+|y|)^{1/10+\mu+} z^{1-\mu}(s,|y|)
\sum_{|\beta|\le 41}|\Gamma^\beta \partial u^I|.
$$
It follows, then, that this is $O(\varepsilon^2)$
by \eqref{3.10}.  Indeed, if $(s,|y|)\in\Lambda_I$, then $s\approx
|y|$ and $z(s,|y|)=\langle c_Is-|y|\rangle$.  For $(s,|y|)\not\in
\Lambda_I$, it follows that $s+|y|\approx |c_Is-|y||$ and
$z^{1-\mu}(s,|y|)\le \langle y\rangle^{1-\mu}$.  Similarly, by
\eqref{3.10}, it follows that the second term in \eqref{3.23} is
bounded by
$$C\varepsilon\sup_{(s,y)\in D^I(t,r)} |y|^{1/2} (1+s+|y|)^\mu
z^{1-\mu}(s,|y|)\sum_{|\beta|\le 41} |\Gamma^\beta \partial u^I|.$$
By \eqref{3.10} and the same considerations as above, this is in turn
$O(\varepsilon^2)$ as desired.

When we replace $\tilde{G}^I$ by
\begin{equation}\label{3.24}
\sum_{\substack{1\le J\le D\\ J\neq I}} \Bigl(\sum_{0\le
j,k,l\le 3} B^{IJ,jk}_{J,l}(1-\rho_J) \partial_l u^J \partial_j\partial_k
u^J + 
\sum_{0\le j,k\le 3}
A^{I,jk}_{JJ} (1-\rho_J) \partial_j u^J\partial_k u^J\Bigr)
\end{equation}
in the left side of \eqref{3.22}, we see that it is bounded by
$$C\sup_{(s,y)\in \text{supp}(1-\rho_J)}|y|
(1+s+|y|)^{1+\mu}z^{1-\mu} (s,|y|)
\sum_{|\beta|\le 41} |\Gamma^\beta \partial u^J|^2.$$
Since $\langle c_Js-|y|\rangle \gtrsim \langle s+|y|\rangle \gtrsim 
z(s,|y|)$ for $(s,|y|)$ in the support of $(1-\rho_J)$, 
it follows easily from \eqref{3.10}
with $\mu=0$ that this term is $O(\varepsilon^2)$ as desired.

Next, we shall examine \eqref{3.22} with $\tilde{G}^I$ replaced 
by the multi-speed terms
$$
\sum_{\substack{1\le J,K\le D\\ K\neq J}} 
\sum_{|\alpha|\le 41} \partial\Gamma^\alpha u^K\sum_{|\alpha|\le
41}\partial \Gamma^\alpha u^J.
$$
Suppose that $(s,|y|)\in \Lambda_J$.  Since $J\neq K$, we have
$|c_Ks-|y||\gtrsim (s+|y|)$.  Thus, if we apply \eqref{3.10} to the $u^K$
piece (with $\mu=0$), we see that the left side of \eqref{3.22} is
controlled by
$$C\varepsilon \sup_{(s,|y|)\in\Lambda_J} |y|^{1/2}(1+s+|y|)^\mu
(1+|c_Js-|y||)^{1-\mu}\sum_{|\beta|\le 41} |\Gamma^\beta \partial u^J|.$$
Since $|y|\approx s$ on $\Lambda_J$, we see that this term is also
$O(\varepsilon^2)$ by another application of \eqref{3.10}.  A symmetric
argument can be used when $(s,|y|)\in\Lambda_K$.  If $(s,|y|)\not\in
\Lambda_J\cup\Lambda_K$, then $|c_Js-|y||,|c_Ks-|y||\approx (s+|y|)$
and the bound follows from two applications of \eqref{3.10} with
$\mu=0$.

Finally, we are left with proving \eqref{3.22} when $\tilde{G}^I$ is
replaced by $R^I+P^I$.  In this case, the right side of \eqref{3.22}
is bounded by
\begin{multline}\label{3.25}
C\sum_{1\le J,K,L\le D} \sup_{(s,y)\in
D^I(t,r)}|y|(1+s+|y|)^{1+\mu}z^{1-\mu}(s,|y|)\\\times \sum_{|\beta|\le
22}|\Gamma^\beta u^J| \sum_{|\beta|\le 22} |\Gamma^\beta u^K|
\sum_{|\beta|\le 40}|\Gamma^\beta u^L|
\\+C \sum_{1\le J,K,L\le D}\sup_{(s,y)\in
D^I(t,r)}|y|(1+s+|y|)^{1+\mu}z^{1-\mu}(s,|y|)\\\times \sum_{|\beta|\le
22}|\Gamma^\beta u^J| \sum_{|\beta|\le 22} |\Gamma^\beta u^K|
\sum_{|\beta|\le 41}|\Gamma^\beta \partial u^L|.
\end{multline}
By the inductive hypothesis \eqref{3.3}, the first term in
\eqref{3.25} is controlled by
$$C\varepsilon^3 \sup_{(s,y)\in D^I(t,r)}
\frac{z^{1-\mu}(s,|y|)}{(1+s+|y|)^{1-\mu}}
\Bigl(1+\log\frac{1+s+|y|}{z(s,|y|)}\Bigr)^3.$$
Since $(\log x)^3/x^{1-\mu}$ is bounded for $x\ge 1$ and $\mu<1$, it
follows that the first term in \eqref{3.25} is $O(\varepsilon^3)$.
For the second term in \eqref{3.25}, if we apply \eqref{3.10}, we see
that it is bounded by
$$C\varepsilon \sum_{1\le J,K\le D} \sup_{(s,y)\in D^I(t,r)}
|y|^{1/2-\mu} (1+s+|y|)^{1+\mu} \sum_{|\beta|\le 22}|\Gamma^\beta u^J|
\sum_{|\beta|\le 22} |\Gamma^\beta u^K|.$$ 
It then follows easily via \eqref{3.3} that this term is also
$O(\varepsilon^3)$ as desired.  This completes the proof of
\eqref{3.22}, and thus, also \eqref{3.3}.

We now wish to prove that \eqref{3.4} can be obtained with $A_2$
replaced by $A_2/2$.  Here, we apply \eqref{2.23a} with $F^I$ replaced
by $B(du)+Q(du,d^2u)$ and $G^I$ replaced by $R(u,du,d^2u)+P(u,du)$ to
see that the left side of \eqref{3.4} is bounded by
\begin{multline}\label{3.26}
C_2\varepsilon+C\sum_{|\beta|\le 63} \int_0^t \int |\Gamma^\beta
[B(du)+Q(du,d^2u)](s,y)|\:\frac{dy\:ds}{\langle y\rangle}
\\+
C\Bigl(1+\log\frac{1+r+c_It}{1+|r-c_It|}\Bigr)\sum_{|\beta|\le 60}
\sup_{(s,y)\in D^I(t,r)} |y|(1+s+|y|)^{1+\mu}
z^{1-\mu}(s,|y|)\\\times|\Gamma^\beta [R(u,du,d^2u)+P(u,du)](s,y)|.
\end{multline}
By our choice of $A_2$, it follows that the first term in
\eqref{3.26} is controlled by $(A_2/10)\varepsilon$.  To complete the
proof of (ii.), it will suffice to show that the last two terms in
\eqref{3.26} are bounded by $C\varepsilon^2 (1+t)^{1/10}
\log(2+t)\Bigl(1+\log\frac{1+t+|x|}{1+|c_It-|x||}\Bigr)$.

Since $B(du)$ and $Q(du,d^2u)$ are quadratic, this is relatively easy
for the second term.  In fact, this term is bounded by
$$C\int_0^t \int \sum_{|\alpha|\le 64} |\Gamma^\alpha \partial
u(s,y)|^2\:\frac{dy\:ds}{\langle y\rangle}.$$
Since this is controlled by the square of the left side of
\eqref{3.7}, the desired bound follows immediately.

To complete the proof of (ii.), it suffices to show that
\begin{multline}\label{3.27}
\sup_{(s,y)\in D^I} |y|(1+s+|y|)^{1+\mu}
z^{1-\mu}(s,|y|)
\sum_{|\beta|\le 60} |\Gamma^\beta [R(u,du,d^2u)+P(u,du)](s,y)|
\\\le C
\varepsilon^3 (1+t)^{1/10}\log(2+t).
\end{multline}
The left side of \eqref{3.27} is controlled by
\begin{multline}\label{3.27a}
C\sup_{(s,y)\in D^I}|y|(1+s+|y|)^{1+\mu} z^{1-\mu}(s,|y|)
\Bigl(\sum_{|\beta|\le 32} |\Gamma^\beta u|\Bigr)^2 \sum_{|\beta|\le
60} |\Gamma^\beta u|
\\+C\sup_{(s,y)\in D^I} |y|(1+s+|y|)^{1+\mu} z^{1-\mu}(s,|y|)
\Bigl(\sum_{|\beta|\le 32} |\Gamma^\beta u|\Bigr)^2 \sum_{|\beta|\le
61} |\Gamma^\beta u'|.
\end{multline}
By \eqref{3.3} and \eqref{3.4}, we see that the first term is
dominated by
$$C\varepsilon^3 \sup_{(s,y)\in D^I}
\Bigl(\frac{z(s,|y|)}{1+s+|y|}\Bigr)^{1-\mu}
\Bigl(1+\log\frac{1+s+|y|}{z(s,|y|)}\Bigr)^3 (1+s)^{1/10} \log(2+s).$$
As above, since $(\log x)^3/x^{1-\mu}$ is bounded for $x>1$ and $\mu$
small, we easily obtain the desired bound.  For the second term in
\eqref{3.27a}, applying \eqref{2.26} and \eqref{3.6} we see that it is
dominated by
$$C\varepsilon (1+t)^{1/40} \sup_{(s,y)} (1+s+|y|)^{1+\mu}
z^{1-\mu}(s,|y|) \Bigl(\sum_{|\beta|\le 32} |\Gamma^\beta u|\Bigr)^2.$$
Applying \eqref{3.3} yields the desired bound \eqref{3.27} and finishes the proof
of (ii.).

\subsection{Proof of (iii.):}  In this section, we finish the
continuity argument, and thus the proof of Theorem \ref{theorem1.2},
by showing that \eqref{3.5}-\eqref{3.7} follow from
\eqref{3.2}-\eqref{3.4}.  

We begin with \eqref{3.5}.  Outside of $\Lambda_I$, $\log
\frac{1+t+|x|}{1+|c_It-|x||}$ is $O(1)$, and \eqref{3.5} follows
directly from \eqref{3.3}.  Within $\Lambda_I$, we have
$t\approx |x|$, and \eqref{3.5} follows from \eqref{2.26} and \eqref{3.2}.

Next, we want to show that the higher order energy bound \eqref{3.6}
holds.  We will apply \eqref{2.5} with
\begin{equation}\label{3.28}
\gamma^{IJ,jk}=-\sum_{\substack{1\le K\le D\\0\le l\le 3}}
B^{IJ,jk}_{K,l}\partial_l u^K -  C^{IJ,jk}(u,u')
\end{equation}
and
\begin{equation}\label{3.29}
G^I=B^I(du)+P^I(u,du).
\end{equation}
In order to prove \eqref{3.6}, by \eqref{2.4}, \eqref{3.2}, and an
induction argument, it will suffice to prove the following.
\begin{lemma}\label{lemma3.2}
Assume that \eqref{3.2}-\eqref{3.5} hold and $M\le 70$.  Additionally,
suppose that 
\begin{equation}\label{3.30}
\sum_{|\alpha|\le M-1} E(\Gamma^\alpha u,t)\le C\varepsilon
(1+t)^{C\varepsilon + \sigma}
\end{equation}
with $\sigma>0$.  Then, there is a constant $C'$ so that
\begin{equation}\label{3.31}
\sum_{|\alpha|\le M} E(\Gamma^\alpha u,t)\le C'\varepsilon (1+t)^{C'\varepsilon+C'\sigma}.
\end{equation}
\end{lemma}

\noindent{\em Proof of Lemma \ref{lemma3.2}:} 
Since 
\begin{equation}\label{3.32}
\sum_{|\alpha|\le M}|[\Box_\gamma, \Gamma^\alpha]u|\le
C\sum_{|\alpha|\le M-1} |\Gamma^\alpha \Box u| +
C\sum_{\substack{|\alpha|+|\beta|\le M \\ |\beta|\le M-1}}
|\Gamma^\alpha \gamma \Gamma^\beta \partial^2 u|
\end{equation}
and since \eqref{3.3} and \eqref{3.5} imply that
\begin{equation}\label{3.33}
\sum_{|\alpha|\le N} |\Gamma^\alpha \gamma|\le \frac{C\varepsilon}{1+t}
\end{equation}
for $N\le 39$, it follows from \eqref{2.5} that
\begin{multline}\label{3.34}
\sum_{|\alpha|\le M} \partial_t E(\Gamma^\alpha u,t) \le
C\sum_{|\alpha|\le M} \|\Gamma^\alpha B(du)(t,\cd)\|_2 + C\sum_{|\alpha|\le
M}\|\Gamma^\alpha P(u,du)(t,\cd)\|_2 \\+ C\sum_{|\alpha|\le M-1}
\|\Gamma^\alpha Q(du,d^2u)(t,\cd)\|_2 + C\sum_{|\alpha|\le M-1}
\|\Gamma^\alpha R(u,du,d^2u)(t,\cd)\|_2 \\+ 
C\sum_{\substack{|\alpha|+|\beta|\le M \\ |\beta|\le M-1}}
\|\Gamma^\alpha \gamma\Gamma^\beta\partial^2u\|_2
+\frac{C\varepsilon}{1+t}
\sum_{|\alpha|\le M} E(\Gamma^\alpha u, t).
\end{multline}

Note that it
follows from \eqref{3.5} that
\begin{multline}\label{3.35}
\sum_{|\alpha|\le M}\|\Gamma^\alpha B^I(du)(t,\cd)\|_2 +
\sum_{|\alpha|\le M-1} \|\Gamma^\alpha Q^I(du,d^2u)(t,\cd)\|_2\\\le
\frac{C\varepsilon}{1+t} \sum_{|\alpha|\le M} E(\Gamma^\alpha u,t).
\end{multline}
Additionally, by \eqref{3.3}, we have
\begin{multline*}
\sum_{|\alpha|\le M}\|\Gamma^\alpha P^I(u,du)(t,\cd)\|_2
+\sum_{|\alpha|\le M-1}\|\Gamma^\alpha R^I(u,du,d^2u)(t,\cd)\|_2\\
\le
C\varepsilon^2 \sum_{|\alpha|\le M, |\beta|\le 1}
\Bigl\|\Bigl(\frac{1+\log (1+t)}{1+t+|x|}\Bigr)^2 \Gamma^\alpha
\partial^\beta u(t,\cd)\Bigr\|_2.
\end{multline*}
Since the coefficients of $\Gamma$ are $O(1+t+|x|)$, it follows from
\eqref{3.3} that this is
\begin{multline}\label{3.36}
\le C\varepsilon^3
(1+t)^{-3/2+} + C\varepsilon^2
\frac{(1+\log(1+t))^2}{1+t} \sum_{|\alpha|\le M-1} E(\Gamma^\alpha u, t)
\\+ \frac{C \varepsilon^2 (1+\log(1+t))^2}{(1+t)^{2}} \sum_{|\alpha|\le M}
E(\Gamma^\alpha u,t).
\end{multline}
The first term on the right
side corresponds to the case $|\alpha|=|\beta|=0$ on the right side of
the previous equation.  Similarly, the second term is for the case
$|\beta|=0$, and the last term bounds the case $|\alpha|,|\beta|\neq
0$.   By a similar argument, the
fifth term on the right of \eqref{3.34} is also controlled by the
right sides of \eqref{3.35} and \eqref{3.36}.  Thus, we see that
\begin{multline}\label{3.37}
\sum_{|\alpha|\le M}\partial_t E(\Gamma^\alpha u,t) \le
\frac{C\varepsilon}{1+t} \sum_{|\alpha|\le M} E(\Gamma^\alpha u,t)\\ +
C\varepsilon^3 (1+t)^{-3/2+\delta} + \frac{C\varepsilon^2
(1+\log(1+t))^2}{1+t} \sum_{|\alpha|\le M-1} E(\Gamma^\alpha u,t).
\end{multline}

Integrating both sides in $t$, applying the smallness assumption on
the data \eqref{1.13} and the inductive hypothesis
\eqref{3.30}, and using Gronwall's inequality yields \eqref{3.31} as
desired.\qed

We are, thus, left with the task of showing \eqref{3.7}.  Applying
\eqref{2.6} with $u$ replaced by $\Gamma^\alpha u$, we see that the
left side of \eqref{3.7} is controlled by 
\begin{multline}\label{3.38}
C (\log (2+t))^{1/2} \Bigl(\sum_{|\alpha|\le 66} \|\Gamma^\alpha f\|_2 
+ \sum_{|\alpha|\le 65} \|\Gamma^\alpha g\|_2 + \sum_{|\alpha|\le 65} \int_0^t
\|\Gamma^\alpha \Box u(s,\cd)\|_2\:ds\Bigr).\end{multline}
By \eqref{1.13}, the first two terms satisfy the desired bound.  
Since
\begin{multline}\label{3.39}
\sum_{|\alpha|\le 65} \|\Gamma^\alpha \Box u(s,\cd)\|_2 \le
C\Bigl\|\sum_{|\alpha|\le 33} |\Gamma^\alpha u'| \sum_{|\alpha|\le
66}|\Gamma^\alpha u'|\Bigr\|_2 \\+ \Bigl\|\Bigl(\sum_{|\alpha|\le
33}|\Gamma^\alpha u|\Bigr)^2 \sum_{|\alpha|\le 67} |\Gamma^\alpha u|\Bigr\|_2,
\end{multline}
we may use \eqref{3.3},\eqref{3.5}, and the fact that the coefficients
of $\Gamma$ are $O(1+t+|x|)$ to see that the right side of
\eqref{3.39} is dominated by
\begin{multline}\label{3.40}
\frac{C\varepsilon}{1+s}\sum_{|\alpha|\le 66}\|\Gamma^\alpha
u'(s,\cd)\|_2 
+
C\varepsilon^2\frac{(1+\log(2+s))^2}{1+s}\sum_{|\alpha|\le 66}
\|\Gamma^\alpha u'(s,\cd)\|_2 
\\+ C\varepsilon^2 \frac{(1+\log(2+s))^2}{(1+s)^{3/2-}}\|(1+s+|\cd|)^{-1/2-}u(s,\cd)\|_2.
\end{multline}
Plugging \eqref{3.39} and \eqref{3.40} into \eqref{3.38}, we see that
the third term of \eqref{3.38} is bounded by the right side of
\eqref{3.7} by using
\eqref{3.3} and \eqref{3.6}.

This completes the proof of (iii.), and hence the proof of Theorem \ref{theorem1.2}.

\newsection{Preliminary estimates in the exterior domain}
In this section, we will collect the exterior domain analogs of the
estimates in Section 2.  Many of these estimates were previously
established in \cite{KSS3}, \cite{MNS}, and \cite{MS}.  The main new
item will be the use of the pointwise estimates found in the second
subsection.

\subsection{Energy estimates}  We begin by gathering the $L^2$
estimates that we will need in order to show global existence in the
exterior domain.  These estimates are from Metcalfe-Sogge \cite{MS}
(see also \cite{KSS3}),
and unless stated otherwise, their proofs can be found there.
Specifically, we will be concerned with solutions $u\in
C^\infty(\R_+\times\ext)$ of the Dirichlet-wave equation

\begin{equation}\label{4.1}
\begin{cases}
\Box_\gamma u=F\\
u|_{\partial\mathcal{K}}=0\\
u|_{t=0}=f,\quad \partial_t u|_{t=0}=g
\end{cases}
\end{equation}
with $\Box_\gamma$ as in \eqref{2.1}.
We shall assume that the $\gamma^{IJ,jk}$ satisfy the symmetry
conditions \eqref{2.2}
as well as the size condition
\begin{equation}\label{4.3}
\sum_{I,J=1}^D \sum_{j,k=0}^3 \|\gamma^{IJ,jk}(t,x)\|_{\infty}\le
\delta
\end{equation}
for $\delta$ sufficiently small (depending on the wave speeds).
 The energy estimate will involve bounds for the gradient of the
perturbation terms
$$\|\gamma'(t,\cd)\|_\infty = \sum_{I,J=1}^D\sum_{j,k,l=0}^3
\|\partial_l \gamma^{IJ,jk}(t,\cd)\|_\infty,$$ and the energy form
associated with $\Box_\gamma$, $e_0(u)=\sum_{I=1}^D e_0^I(u)$,
where $e_0^I(u)$ is given by \eqref{2.3}.

The most basic estimate will lead to a bound for
$$E_M(t)=E_M(u)(t)=\int \sum_{j=0}^M e_0(\partial^j_t
u)(t,x)\:dx.$$

\begin{lemma}\label{lemma4.1}
Fix $M=0,1,2,\dots$, and assume that the perturbation terms
$\gamma^{IJ,jk}$ satisfy \eqref{2.2} and \eqref{4.3}.  Suppose also that $u\in C^\infty$
solves \eqref{4.1} and for every $t$, $u(t,x)=0$ for large $x$.
Then there is an absolute constant $C$ so that
\begin{equation}\label{4.5}
\partial_t E^{1/2}_M(t)\le C\sum_{j=0}^M \|\Box_\gamma
\partial_t^ju(t,\cd)\|_2+C\|\gamma'(t,\cd)\|_\infty E^{1/2}_M(t).
\end{equation}
\end{lemma}

Before stating the next result, let us introduce some notation. If
$P=P(t,x,D_t,D_x)$ is a differential operator, we shall let
$$[P,\gamma^{kl}\partial_k\partial_l]u=\sum_{1\le I,J\le
D}\sum_{0\le k,l\le 3} |[P,\gamma^{IJ,kl}\partial_k\partial_l]
u^J|.$$

In order to generalize the above energy estimate to include the
more general vector fields $L, Z$, we will need to use a variant
of the scaling vector field $L$.  We fix a bump function $\eta\in
C^\infty(\R^3)$ with $\eta(x)=0$ for $x\in \mathcal{K}$ and
$\eta(x)=1$ for $|x|>1$.  Then, set $\tilde{L}=\eta(x)r\partial_r
+ t\partial_t$. Using this variant of the scaling vector field and
an elliptic regularity argument, one can establish

\begin{proposition}\label{proposition4.2}
Suppose that the constant in \eqref{4.3} is small.  Suppose
further that
\begin{equation}\label{4.6}
\|\gamma'(t,\cd)\|_\infty \le \delta/(1+t),
\end{equation}
and
\begin{multline}\label{4.7}
\sum_{\substack{j+\mu\le
N_0+\nu_0\\\mu\le\nu_0}}\left(\|\tilde{L}^\mu\partial^j_t\Box_\gamma
u(t,\cd)\|_2+\|[\tilde{L}^\mu\partial_t^j,\gamma^{kl}\partial_k\partial_l]u(t,\cd)\|_2\right)\\
\le \frac{\delta}{1+t}\sum_{\substack{j+\mu\le N_0+\nu_0\\
\mu\le\nu_0}}\|\tilde{L}^\mu\partial_t^j
u'(t,\cd)\|_2+H_{\nu_0,N_0}(t),
\end{multline}
where $N_0$ and $\nu_0$ are fixed.  Then
\begin{multline}\label{4.8}
\sum_{\substack{|\alpha|+\mu\le N_0+\nu_0\\
\mu\le\nu_0}}\|L^\mu\partial^\alpha u'(t,\cd)\|_2 \\ \le
C\sum_{\substack{|\alpha|+\mu\le N_0+\nu_0-1\\
\mu\le\nu_0}}\|L^\mu\partial^\alpha\Box u(t,\cd)\|_2 +
C(1+t)^{A\delta}\sum_{\substack{\mu+j\le N_0+\nu_0\\ \mu\le
\nu_0}}\left(\int e_0(\tilde{L}^\mu\partial_t^j
u)(0,x)\:dx\right)^{1/2}\\
+C(1+t)^{A\delta}\Bigl(\int_0^t \sum_{\substack{|\alpha|+\mu\le
N_0+\nu_0-1\\\mu\le\nu_0-1}}\|L^\mu\partial^\alpha \Box
u(s,\cd)\|_2\:ds +\int_0^t H_{\nu_0,N_0}(s)\:ds\Bigr)\\
+C(1+t)^{A\delta}\int_0^t \sum_{\substack{|\alpha|+\mu\le
N_0+\nu_0\\ \mu\le \nu_0-1}}\|L^\mu\partial^\alpha
u'(s,\cd)\|_{L^2(|x|<1)}\:ds,
\end{multline}
where the constants $C$ and $A$ are absolute constants.
\end{proposition}

In practice $H_{\nu_0,N_0}(t)$ will involve $L^2_x$ norms
of $|L^\mu \partial^\alpha u'|^2$ with $\mu+|\alpha|$ much smaller
than $N_0+\nu_0$, and so the integral involving $H_{\nu_0,N_0}$
can be dealt with using an inductive argument and the weighted
$L^2_tL^2_x$ estimates that will be presented at the end of this
subsection.

In proving our existence results for \eqref{1.1}, the key step
will be to obtain a priori $L^2$-estimates involving $L^\mu
Z^\alpha u'$.  Begin by setting
\begin{equation}\label{4.9}
Y_{N_0,\nu_0}(t)=\int \sum_{\substack{|\alpha|+\mu\le
N_0+\nu_0\\\mu\le\nu_0}}e_0(L^\mu Z^\alpha u)(t,x)\:dx.
\end{equation}
We, then, have the following proposition which shows how the
$L^\mu Z^\alpha u'$ estimates can be obtained from the ones
involving $L^\mu\partial^\alpha u'$.
\begin{proposition}\label{proposition4.3}
Suppose that the constant $\delta$ in \eqref{4.3} is small and that
\eqref{4.6} holds.  Then,
\begin{multline}\label{4.10}
\partial_t Y_{N_0,\nu_0}\le C Y^{1/2}_{N_0,\nu_0} \sum_{\substack{
|\alpha|+\mu\le N_0+\nu_0\\ \mu\le\nu_0}} \|\Box_\gamma L^\mu Z^\alpha
u(t,\cd)\|_2 + C \|\gamma'(t,\cd)\|_\infty Y_{N_0,\nu_0} \\
+C \sum_{\substack{|\alpha|+\mu\le N_0+\nu_0+1\\ \mu\le
\nu_0}} \|L^\mu \partial^\alpha u'(t,\cd)\|^2_{L^2(|x|<1)}.
\end{multline}
\end{proposition}

As in \cite{KSS2} and \cite{KSS3} we shall also require some
weighted $L^2_tL^2_x$ estimates.  They will be used, for example,
to control the local $L^2$ norms such as the last term in
\eqref{4.10}.  For convenience, for the remainder of this subsection,
allow $\Box=\partial_t^2-\Delta$ to denote the unit speed, scalar
d'Alembertian.  
The transition from the following estimates to
those involving \eqref{1.2} is straightforward.  Also, allow
$$S_T=\{[0,T]\times\ext\}$$
to denote the time strip of height $T$ in $\Rplus\times\ext$.

We, then, have the following proposition which is an exterior domain
analog of \eqref{2.6}.  
\begin{proposition}\label{proposition4.4}
Fix $N_0$ and $\nu_0$.  Suppose that $\mathcal{K}$ satisfies the
local exponential energy decay \eqref{1.10}.  Suppose also that $u\in
C^\infty$ satisfies $u(t,x)=0$, $t<0$.
Then
there is a constant $C=C_{N_0,\nu_0,\mathcal{K}}$ so that if $u$
vanishes for large $x$ at every fixed $t$
\begin{multline}\label{4.11}
(\log(2+T))^{-1/2}\sum_{\substack{|\alpha|+\mu\le
N_0+\nu_0\\\mu\le\nu_0}} \|\langle x\rangle^{-1/2}L^\mu \partial^\alpha
u'\|_{L^2(S_T)} 
\\\le 
C\int_0^T \sum_{\substack{|\alpha|+\mu\le
N_0+\nu_0+1\\\mu\le\nu_0}} \|\Box L^\mu \partial^\alpha
u(s,\cd)\|_2\:ds +C \sum_{\substack{|\alpha|+\mu\le
N_0+\nu_0\\\mu\le\nu_0}} \|\Box L^\mu\partial^\alpha
u\|_{L^2(S_T)}
\end{multline}
and
\begin{multline}\label{4.12}
(\log(2+T))^{-1/2}\sum_{\substack{|\alpha|+\mu\le
N_0+\nu_0\\\mu\le\nu_0}}\|\langle x\rangle^{-1/2}L^\mu Z^\alpha
u'\|_{L^2(S_T)}\\
\le 
C\int_0^T \sum_{\substack{|\alpha|+\mu\le
N_0+\nu_0+1\\\mu\le\nu_0}} \|\Box L^\mu Z^\alpha u(s,\cd)\|_2\:ds
+ C \sum_{\substack{|\alpha|+\mu\le N_0+\nu_0\\\mu\le\nu_0}}
\|\Box L^\mu Z^\alpha u\|_{L^2(S_T)}.
\end{multline}
\end{proposition}

We end this subsection with a couple of results that follow from the
local energy decay \eqref{1.10}.
\begin{lemma}\label{lemma4.5}
Suppose that \eqref{1.10} holds and that $\Box u(t,x)=0$ for $|x|>4$. 
Suppose also that $u(t,x)=0$ for $t\le 0$.  Then, if $N_0$
and $\nu_0$ are fixed and if $c>0$ is as in \eqref{1.10}, the following estimate holds : 
\begin{multline}\label{4.13}
\sum_{\substack{|\alpha|+\mu\le N_0+\nu_0\\\mu\le\nu_0}} \|L^\mu
\partial^\alpha u'(t,\cd)\|_{L^2(\{\ext\, : \,|x|<4\})}
\\\le C \sum_{\substack{|\alpha|+\mu\le
N_0+\nu_0-1\\\mu\le\nu_0}}\|L^\mu\partial^\alpha \Box
u(t,\cd)\|_2 
\\+C\int_0^t e^{-(c/2)(t-s)}\sum_{\substack{|\alpha|+\mu\le
N_0+\nu_0+1\\\mu\le\nu_0}}\|L^\mu \partial^\alpha \Box u(s,\cd)\|_2\:ds.
\end{multline}
\end{lemma}

To be able to handle the last term in \eqref{4.8}, we shall need
the following.
\begin{lemma}\label{lemma4.6}
Suppose that \eqref{1.10} holds, and suppose that $u\in C^\infty$
satisfies $u(t,x)=0$ for $t<0$.  Then, for
fixed $N_0$ and $\nu_0$ and $t>2$,
\begin{multline}\label{4.14}
\sum_{\substack{|\alpha|+\mu\le N_0+\nu_0\\\mu\le\nu_0}} \int_0^t
\|L^\mu \partial^\alpha u'(s,\cd)\|_{L^2(|x|<2)}\:ds \\
\le C \sum_{\substack{|\alpha|+\mu\le N_0+\nu_0+1\\\mu\le\nu_0}}
\int_0^t \left(\int_0^s \|L^\mu \partial^\alpha \Box
u(\tau,\cd)\|_{L^2(||x|-(s-\tau)|<10)}\:d\tau\right)\:ds.
\end{multline}
\end{lemma}

\subsection{Pointwise estimates}  Here, we will describe the various
pointwise estimates that we shall require.  These include variants of
those of Keel-Smith-Sogge \cite{KSS3} and Metcalfe-Sogge \cite{MS}
and exterior domain analogs of the estimates of Kubota-Yokoyama
\cite{KY}.

Let us begin with the former.  We will need analogs of the pointwise estimates
of \cite{KSS3} and \cite{MS} 
that allow Cauchy data that vanishes
in a neighborhood of the obstacle.  That is, we will estimate
solutions of the scalar wave equation with boundary
$(\partial_t^2-\Delta)w(t,x)=F(t,x)$.
Additionally, we will require that $w(0,x)=\partial_t w(0,x)=0$ if $|x|\le 6$,
and $F(t,x)=0$ if $|x|\le 6$ and $0\le t\le 1$.  With these
assumptions, we can greatly reduce the technical details involving the
compatibility conditions.  In the sequel, we will reduce our study to
this case.
Assuming, as we do throughout, that $\mathcal{K}\subset\{x\in\R^3\,
:\, |x|<1\}$, we have
\begin{theorem}\label{theorem4.7}
Suppose that the local energy decay bounds \eqref{1.10} hold for
$\mathcal{K}$.  Additionally, assume that $w(t,x)=0$ for 
$x\in \partial \mathcal{K}$, $w(0,x)=\partial_t w(0,x)=0$
for $|x|\le 6$, and $F(t,x)=0$ if $0\le t\le 1$ and $|x|\le 6$.  
Then, if $|\alpha|=M$,
\begin{multline}\label{4.15}
(1+t+|x|)|L^\nu Z^\alpha w(t,x)|\le 
C\sum_{\substack{j+|\beta|+k\le \nu+M+8 \\ j\le 1}} \|\langle x\rangle^{j+|\beta|}\partial_x^\beta\partial_t^{k+j}w(0,x)\|_2 \\
+C\int_0^t \int_\ext \sum_{\substack{|\beta|+\mu\le M+\nu + 7\\\mu\le
\nu+1}} |L^\mu Z^\beta F(s,y)|\:\frac{dy\:ds}{|y|} \\
+C \int_0^t \sum_{\substack{|\beta|+\mu\le M+\nu+4\\ \mu\le \nu+1}}
\|L^\mu \partial^\beta F(s,\cd)\|_{L^2(\{x\in\ext\, :\, |x|<2\})}\:ds.
\end{multline}
\end{theorem}

\noindent{\em Proof of Theorem \ref{theorem4.7}:}
If $w$ has vanishing Cauchy data with $F(t,x)=0$ for $0\le t\le 1$ and $x\in \ext$,
\eqref{4.15} follows from Theorem 3.1 in \cite{MS}.
We, thus, may assume $F(t,x)=0$ for $0\le t<\infty$ and $|x|\le 6$ and
that the Cauchy data is as stated above.
The proof follows from the arguments of \cite{MS} for the
inhomogeneous case very closely.  We include a sketch of the proof for
completeness.

We first note that if we argue as in \cite{KSS3} (Lemma 4.2) we have
\begin{multline}\label{4.16}
(1+t+|x|)|L^\nu Z^\alpha w(t,x)|\le 
C\sum_{\substack{j+|\beta|+k\le M+\nu+4 \\ j\le 1}} \|\langle x\rangle^{j+|\beta|}\partial_x^\beta\partial_t^{k+j}w(0,x)\|_2 \\
+C\sum_{\substack{|\beta| +\mu \le M+\nu+3\\ \mu\le \nu+1 }} \int_0^t\int_\ext |L^\mu Z^\beta F(s,y)|\:\frac{dy\:ds}{|y|} \\
+C\sum_{\substack{|\beta|+\mu\le M+\nu+1 \\ \mu\le \nu}}\sup_{\substack{0\le s\le t \\ |y|\le 2} }(1+s)|L^\mu \partial^\beta w(s,y)|.
\end{multline}

\noindent While the arguments in \cite{KSS3} are given for vanishing
Cauchy data, straightforward modifications allow the current setting.

It remains to prove bound in the region $|x|<2$.  
We show 
\begin{multline}\label{4.17}
\sum_{\substack{|\beta|+\mu\le M+\nu+1\\ \mu\le \nu}}
\sup_{\substack{0\le s\le t\\ |y|\le 2}}(1+s)|L^\mu \partial^\beta w(s,y)|
\le C\sum_{\substack{j+|\alpha|+k\le M+\nu+8 \\ j\le 1}}\|\langle x\rangle^{j+|\alpha|}\partial_x^\alpha\partial_t^{j+k} w(0,x)\|_2 \\
+C\sum_{\substack{|\beta|+\mu\le M+\nu+7\\ \mu\le
\nu+1}}\int_0^t\int_\ext |L^\mu Z^\beta F(s,y)|\:\frac{dy\:ds}{|y|}.
\end{multline}
To see this, write $w=w_0+w_r$ where $w_0$ solves the boundaryless
wave equation $(\partial_t^2-\Delta)w_0=F$ with initial data
$w_0(0,\cd)=w(0,\cd)$ and $\partial_t w_0(0,\cd)=\partial_tw(0,\cd)$.  If we fix $\eta\in
C^\infty_0(\R^3)$ with $\eta(x)\equiv 1$ for $|x|<2$ and $\eta(x)\equiv
0$ for $|x|\ge 3$ and set $\tilde{w}=\eta w_0+w_r$, it follows that
$w=\tilde{w}$ for $|x|<2$.
Thus, it will suffice to show \eqref{4.17}
with $w$ replaced by $\tilde{w}$.  Notice that $\tilde{w}$ solves the Dirichlet-wave
equation
$$(\partial_t^2-\Delta)\tilde{w}=-2\nabla \eta\cd \nabla_x w_0-(\Delta
\eta)w_0$$
with vanishing initial data since the support of $\eta$ does not intersect the 
supports of $F$, $w(0,\cd)$ and $\partial_tw(0,\cd)$ and   
that this forcing term vanishes unless $2\le |x|\le 3$.

In order to complete the proof, we begin by noting the following
consequence of the Fundamental Theorem of Calculus:
$$\sup_{|y|\le 2}|(1+s)L^\mu \partial^\beta \tilde{w}(s,y)|\le 
C\sum_{j=0,1}\sup_{|y|\le 2}\int_0^s
|(\tau\partial_\tau)^j L^\mu \partial^\beta \tilde{w}(\tau,y)|\:d\tau.$$
Using Sobolev's lemma and the fact that the Dirichlet condition allows
us to control $\tilde{w}$ locally by $\tilde{w}'$, we see that 
the left hand side of \eqref{4.17} is bounded by 
\begin{multline*}
 C\sum_{\substack{|\beta|+\mu \le M+\nu+3\\ \mu\le \nu}}\sum_{j=0,1} 
\int_0^t\|(\tau d\tau)^jL^\mu\partial^\beta
\tilde{w}(\tau,y)\|_{L^2(\ext,|x|\le4)} d\tau 
\\ \le C\sum_{\substack{|\beta|+\mu\le M+\nu+3\\ \mu\le \nu+1}}
\int_0^t\|L^\mu\partial^\beta \tilde{w}'(\tau,y)\|_{L^2(\ext,|x|\le4)} d\tau
\end{multline*}
By \eqref{4.13}, it follows that the right side of the above estimate is
controlled by
$$C\int_0^t \sum_{\substack{|\beta|+\mu\le M+\nu+5\\\mu\le \nu+1}}
\|L^\mu \partial^\beta w_0(s,\cd)\|_{L^\infty(2\le |x|\le 3)}\:ds.$$

From \eqref{2.10}, \eqref{2.16}, and the fact that $1/t\le
1/|y|$ on the domain of integration in \eqref{2.10}, we have
\begin{multline}\label{4.18}
\|L^\mu \partial^\beta w_0(s,\cd)\|_{L^\infty(2\le |x|\le 3)} \le
C\sum_{|\alpha|\le 2} \int_{|s-|y||\le 4} |(\Omega^\alpha\nabla L^\mu \partial^\beta 
 w_0)(0,y)|\frac{dy}{|y|} 
\\
+C\sum_{|\alpha|\le 2} \int_{|s-|y||\le 4} |(\Omega^\alpha L^\mu
\partial^\beta w_0)(0,y)|\frac{dy}{|y|^2}
\\
+C\sum_{|\alpha|\le 2} \int_{|s-|y||\le 4} |(\Omega^\alpha \partial_t L^\mu \partial^\beta  w_0)(0,y)|\frac{dy}{|y|}
\\
+C\sum_{|\alpha|\le 2}\int_0^s\int_{|s-\tau-|y||\le 4}|\Omega^\alpha(\partial_t^2-\Delta)L^\mu\partial^\beta w_0(\tau,y)|\:\frac{dy\:d\tau}{|y|}.
\end{multline}
Since the sets $\Lambda_s=\{y\,:\, |s-|y||\le 4\}$ satisfy
$\Lambda_s\cap \Lambda_s'=\emptyset$ if $|s-s'|\ge 10$, if we sum over
$|\beta|+\mu\le M+\nu+5$, $\mu\le \nu+1$, and integrate over $s\in
[0,t]$, we conclude that the left side of \eqref{4.17} is controlled by
\begin{multline*}
C\sum_{k+|\beta|\le M+\nu+7} \int \langle x\rangle^{|\beta|}
|(\partial^\beta \partial_t^k \nabla w)(0,y)|\:\frac{dy}{|y|}
\\+ C\sum_{k+|\beta|\le M+\nu+7} \int \langle x\rangle^{|\beta|}
|(\partial^\beta \partial_t^k w)(0,y)|\:\frac{dy}{|y|^2}
\\+C \sum_{k+|\beta|\le M+\nu+7} \int \langle x\rangle^{|\beta|}
|(\partial^\beta \partial_t^k \partial_t w)(0,y)|\:\frac{dy}{|y|}
\\+C\sum_{\substack{|\beta|+\mu\le M+\nu+7 \\ \mu\le
\nu+1}}\int_0^t\int_\ext|L^\mu Z^\beta F(s,y)|\:\frac{dy\:ds}{|y|}.
\end{multline*}
Using the Schwarz inequality, \eqref{4.17}, and thus \eqref{4.15},
follows. \qed

For the remainder of the estimates in this section, it will suffice to
take $w$ to be a solution to the following Dirichlet-wave equation
with vanishing initial data.  
\begin{equation}\label{4.19}
\begin{cases}
(\partial_t^2 - c_I^2 \Delta)w(t,x)=F(t,x), \quad (t,x)\in
\R_+\times\ext\\
w(t,x)=0,\quad x\in\partial\mathcal{K}\\
w(t,x)=0,\quad t\le 0.
\end{cases}
\end{equation}
In the sequel, we will reduce showing
that \eqref{1.1} has a global solution to showing that an equivalent
system of nonlinear wave equations with vanishing data has a global
solution.  Since the previous theorem will suffice to make this
reduction, it is unnecessary to consider nonvanishing Cauchy data in
the subsequent estimates. 

We will need the following version of \eqref{4.15} that does not require a
loss of a scaling vector field on the right.  
\begin{theorem}\label{theorem4.8}
Suppose that the local energy decay bound \eqref{1.10} holds for
$\mathcal{K}$.  Suppose that $w$ is a solution to \eqref{4.19} and
$|\alpha|=M$.  Then,
\begin{multline}\label{4.20}
(1+|x|)|L^{\nu_0} Z^\alpha w(t,x)| \le \int_0^t \int_\ext
\sum_{\substack{|\beta|+\nu\le M+\nu_0+6\\\nu\le \nu_0}}|L^\nu Z^\beta
F(s,y)|\:\frac{dy\:ds}{|y|} 
\\+C\int_0^t \sum_{\substack{|\beta|+\nu\le M+\nu_0+3\\\nu\le\nu_0}}
\|L^\nu \partial^\beta F(s,\cd)\|_{L^2(\{x\in\ext\,:\,|x|<4\})}ds.
\end{multline}
\end{theorem}

Here, we refer the
reader to similar arguments in the previous articles of Keel-Smith-Sogge \cite{KSS3}
(Theorem 4.1), Metcalfe-Sogge \cite{MS} (Theorem 3.1), and the authors
\cite{MNS} (Lemma 3.3, Lemma 3.4).  Since we are only requiring decay
in $|x|$, the proof is based only on the Minkowski estimate
\begin{equation}\label{4.21}
\begin{split}
|x||w_0(t,x)|&\le C\int_0^t \int_{||x|-(t-s)|}^{|x|+(t-s)}
 \sup_{|\theta|=1} |\Box w_0(s,r\theta)|\:r\:dr\:ds\\
&\le C \int_0^t \int_{\{y\in\R^3\,:\, |y|\in
 [||x|-(t-s)|,|x|+(t-s)]\}} \sum_{|a|\le 2} |\Omega^a \Box w_0(s,y)|\:\frac{dy\:ds}{|y|}.
\end{split}
\end{equation}
We, thus, do not require the additional $L$ that appears on the right
side of the estimates in \cite{KSS3}, \cite{MS}, and \cite{MNS}.

Letting $\Lambda_I$ be the small conic neighborhood of the
characteristic cone $|x|=c_It$ for $\Box_{c_I}$ defined by
\eqref{2.19}, we also have the following estimate when the forcing
term is localized to such a region.  This is an analog of \eqref{2.23}
for the Dirichlet-wave equation.

\begin{theorem}\label{theorem4.9}
Let $w$ be a solution to \eqref{4.19}.  Suppose that $F(t,x)$ is
supported in some $\Lambda_J$ for $J\neq I$.  Then, there are constants $c,c',C>0$
depending on $c_I,c_J$ so that for $t>2$, $I,J=1,2,\dots, D$,
\begin{multline}\label{4.24}
(1+t+|x|)|L^\nu Z^\alpha w(t,x)|\\\le C \int_{\max(0,c|c_I t-|x||-1)}^t
\int_{|y|\approx s} \sum_{\substack{|\beta|+\mu\le
|\alpha|+\nu+3\\\mu\le\nu+1}}|L^\mu Z^\beta F(s,y)|\:\frac{dy\:
ds}{|y|}
\\
+C \sum_{\substack{|\beta|+\mu\le
|\alpha|+\nu+6\\\mu\le \nu}} \sup_{0\le s\le t} \int |L^\mu Z^\beta F(s,y)|\:dy
\\
+C\sup_{0\le s\le t}\sum_{\substack{|\beta|+\mu\le
|\alpha|+\nu+7\\\mu\le \nu+1}}\int_{c' s}^s \int_{|y|\approx \tau}
|L^\mu Z^\beta F(\tau,y)|\:\frac{dy\:d\tau}{|y|}
\\
+C \sup_{0\le s\le t} (1+s)\sum_{\substack{|\beta|+\mu\le
|\alpha|+\nu+3\\\mu\le\nu}} \|L^\mu \partial^\beta F(s,\cd)\|_{L^\infty(|x|<10)}.
\end{multline}
\end{theorem}

\noindent Here, as before, $|y|\approx s$ indicates that there is some
positive 
constant $c$
so that $\frac{1}{c} s \le |y|\le c s$.

We shall need an analog of Lemma \ref{lemma2.3},
the result of
Kubota-Yokoyama \cite{KY}, for Dirichlet-wave equations.  With $z$
as in \eqref{2.20}, we have
\begin{theorem}\label{theorem4.11}
Let $I=1,2,\dots D$, and let $w$ be a solution to \eqref{4.19}.  Then,
for any $\mu >0$,
\begin{multline}\label{4.32}
(1+t+r)\Bigl(1+\log\frac{1+t+r}{1+|c_I t-r|}\Bigr)^{-1} |L^{\nu_0} Z^\alpha
w(t,x)|\\
\le C \sup_{(s,y)} |y|(1+s+|y|)^{1+\mu}
z^{1-\mu}(s,|y|)\sum_{\substack{|\beta|+\nu \le
|\alpha|+\nu_0\\\nu\le\nu_0}}|L^\nu Z^\beta F(s,y)|\\
+C\sup_{(s,y)}|y|(1+s+|y|)^{1+\mu}z^{1-\mu}(s,|y|) \sum_{\substack{|\beta|+\nu\le
|\alpha|+\nu_0+3\\\nu\le\nu_0}} |L^\nu \partial^\beta \partial F(s,y)|.\\
\end{multline}
\end{theorem}

The proofs of Theorem \ref{theorem4.9} and Theorem \ref{theorem4.11}
are quite similar, and we will only provide the proof of the latter.  
In order to prove Theorem \ref{theorem4.9}, we need only
replace the applications of \eqref{2.21} by \eqref{2.23} which is
the appropriate free space analog of \eqref{4.24}.

\noindent{\em Proof of Theorem \ref{theorem4.11}:} 
We begin by claiming that
\begin{multline}\label{4.34}
(1+t+r)\Bigl(1+\log\frac{1+t+r}{1+|c_It-r|}\Bigr)^{-1} |L^\nu_0
Z^\alpha w(t,x)| 
\\\le C\sup_{(s,y)} |y|(1+s+|y|)^{1+\mu} z^{1-\mu}(s,|y|)
\sum_{\substack{|\beta|+\nu\le |\alpha|+\nu_0\\\nu\le\nu_0}} |L^\nu
Z^\beta F(s,y)|
\\+C\sup_{(s,y), |y|<2} (1+s)\sum_{\substack{|\beta|+\nu\le
|\alpha|+\nu_0+1\\\nu\le \nu_0}} |L^\nu \partial^\beta w(s,y)|.
\end{multline}
Indeed, over $|x|<2$, the left side is clearly bounded by the second
term on the right side since the coefficients of $Z$ are $O(1)$ on
this set.  To see the estimate on $|x|\ge 2$, we fix a cutoff function
$\rho\in C^\infty$ where $\rho(x)\equiv 0$ for $|x|<3/2$ and
$\rho(x)\equiv 1$ for $|x|>2$.  If we let $w_j$ denote the solutions
to the boundaryless wave equations
$(\partial_t^2-c_I^2\Delta)w_j=G_j$, $j=1,2$ where
$G_1=\rho(\partial_t^2-c_I^2\Delta)w$ and
$G_2=-2c_I^2\nabla\rho\cdot\nabla_x w-c_I^2(\Delta\rho)w$, we see that $w=w_1+w_2$.
Since $[\Box,Z]=0$ and $[\Box,L]=2\Box$, we can establish the bound
for the $w_1$ piece by applying \eqref{2.21}.  Arguing as in Lemma 4.2
of Keel-Smith-Sogge \cite{KSS3}, we see that the $w_2$ term is bounded
by the second term on the right side of \eqref{4.34}.

To finish the proof, it thus suffices to show
\begin{multline}\label{4.35}
\sup_{0\le s\le t} (1+s)\sum_{\substack{|\beta|+\nu\le
|\alpha|+\nu_0+1\\\nu\le\nu_0}} \|L^\nu \partial^\beta
w(s,\cd)\|_{L^\infty(|x|<2)}
\\\le C\sup_{(s,y)} |y|(1+s+|y|)^{1+\mu} z^{1-\mu}(s,|y|)
\sum_{\nu\le\nu_0} |L^\nu F(s,y)|
\\+C\sup_{(s,y)} |y|(1+s+|y|)^{1+\mu} z^{1-\mu}(s,|y|)
\sum_{\substack{|\beta|+\nu\le |\alpha|+\nu_0+3\\\nu\le \nu_0}} |L^\nu
\partial^\beta \partial F(s,y)|.
\end{multline}

When $F(s,y)=0$ for $|y|>10$, we can apply the following lemma, which
is essentially Lemma 3.3 from \cite{MNS}.
\begin{lemma}\label{lemma4.10}
Suppose that $w$ is as above.  Suppose further that
$(\partial_t^2-c_I^2\Delta)w(s,y)=F(s,y)=0$ if $|y|>10$.  Then,
\begin{equation}\label{4.30}
(1+t)\sup_{|x|<2} |L^\nu \partial^\alpha w(t,x)|\le C\sup_{0\le s\le
t} \sum_{\substack{|\beta|+\mu\le |\alpha|+\nu+2\\\mu\le\nu}}
(1+s)\|L^\mu \partial^\beta F(s,\cd)\|_2.
\end{equation}
\end{lemma}
Since $F$ is supported on $|y|<10$ and since $|y|$ is bounded below on
the complement of $\mathcal{K}$, it follows that this term is
controlled by the right side of \eqref{4.35}.

We also need an estimate for solutions whose forcing terms vanish near
the obstacle.  Assume now that $F(s,y)=0$ for $|y|<5$ and write
$w=w_0+w_r$
where $w_0$ solves the boundaryless wave equation
$(\partial_t^2-c_I^2\Delta)w_0=F$ with vanishing initial data.
Fixing $\eta\in C^\infty_0(\R^3)$ satisfying $\eta(x)\equiv 1$ for
$|x|<2$ and $\eta(x)\equiv 0$ for $|x|\ge 3$ and setting
$\tilde{w}=\eta w_0+w_r$, we see that $w=\tilde{w}$ on $|x|<2$.  Since
$\tilde{w}$ solves the Dirichlet-wave equation
$$(\partial_t^2-c_I^2\Delta)\tilde{w}=-2c_I\nabla \eta\cdot\nabla_x
w_0-c_I^2 (\Delta\eta)w_0=G$$
and $G$ vanishes unless $2\le |x|\le 3$, we may apply Lemma
\ref{lemma4.10} to see
\begin{multline*}
(1+t)\sup_{|x|<2}
\sum_{\substack{|\beta|+\nu \le |\alpha|+\nu_0+1\\\nu\le\nu_0}}
 |L^\nu \partial^\beta
w(t,x)|\le (1+t)\sup_{|x|<2}
\sum_{\substack{|\beta|+\nu \le |\alpha|+\nu_0+1\\\nu\le\nu_0}}
|L^\nu \partial^\beta \tilde{w}(t,x)|\\
\le C\sup_{0\le s\le t} \sum_{\substack{|\beta|+\nu \le
|\alpha|+\nu_0+3\\\nu\le\nu_0}}(1+s)
|L^\nu\partial^\beta w_0'(s,\cd)\|_{L^\infty(|x|<3)} \\
 + C\sup_{0\le s\le t} (1+s)\sum_{\nu\le\nu_0}
\|L^\nu w_0(s,\cd)\|_{L^\infty(|x|<3)}.
\end{multline*}
We thus see that \eqref{4.35} follows from an application of \eqref{2.21}.\qed

\subsection{Sobolev-type estimates}
In this subsection, we state the exterior domain analogs of Lemma
\ref{lemma2.7} that we will require.  The proofs of the relevant
extensions to the exterior domain can be found in \cite{MNS} (Lemma 4.2 and Lemma 4.3).

\begin{lemma}\label{lemma4.12}
Suppose that $u(t,x)\in C^\infty_0(\R\times\ext)$ vanishes for $x\in
\partial\mathcal{K}$.  Then, if $|\alpha|=M$ and $\nu$ are fixed
\begin{multline}\label{4.36}
\|\langle c_I t-r\rangle L^\nu Z^\alpha \partial^2 u(t,\cd)\|_2\le
C\sum_{\substack{|\beta|+\mu\le M+\nu+1\\\mu\le\nu+1}} \|L^\mu Z^\beta
u'(t,\cd)\|_2 
\\+C\sum_{\substack{|\beta|+\mu\le M+\nu\\\mu\le\nu}}\|\langle
t+r\rangle L^\mu Z^\beta (\partial_t^2-c_I^2 \Delta)u(t,\cd)\|_2 +
C(1+t)\sum_{\mu\le \nu}\|L^\mu u'(t,\cd)\|_{L^2(|x|<2)}.
\end{multline}
and
\begin{multline}\label{4.37}
r^{1/2+\theta}\langle c_I t-r\rangle^{1-\theta} |\partial L^\nu Z^\alpha
u(t,x)|\le C\sum_{\substack{|\beta|+\mu\le
M+\nu+2\\\mu\le\nu+1}}\|L^\mu Z^\beta u'(t,\cd)\|_2
\\+C\sum_{\substack{|\beta|+\mu\le M+\nu+1\\ \mu\le\nu}}\|\langle
t+r\rangle L^\mu Z^\beta (\partial_t^2-c_I^2 \Delta)u(t,\cd)\|_2 +
C(1+t)\sum_{\mu\le\nu} \|L^\mu u'(t,\cd)\|_{L^\infty(|x|<2)}
\end{multline}
for any $0\le\theta\le 1/2$.
\end{lemma}


\newsection{The continuity argument in the exterior domain}
In this section, we will prove the main result, Theorem
\ref{theorem1.1}.  We shall take $N=322$ in the smallness hypothesis
\eqref{1.11}.  This can be improved considerably, but here
we will take such a liberty in order to avoid unnecessary
technicalities.

Our global existence theorem will be based on the following local
existence result.
\begin{theorem}\label{theorem5.1}
Suppose that $f$ and $g$ are as in Theorem \ref{theorem1.1} with $N\ge
7$ in \eqref{1.11}.  Then, there is a $T>0$ so that the initial value
problem \eqref{1.1} with this initial data has a $C^2$ solution
satisfying
$$u\in L^\infty([0,T];H^N(\ext))\cap C^{0,1}([0,T];H^{N-1}(\ext)).$$
The supremum of such $T$ is equal to the supremum of all $T$ where the
initial value problem has a $C^2$ solution with $\partial^\alpha u$
bounded for all $|\alpha|\le 2$.  Also, one can take $T\ge 2$ if
$\|f\|_{H^N}+\|g\|_{H^{N-1}}$ is small enough.
\end{theorem}

This is essentially from Keel-Smith-Sogge \cite{KSS} 
(Theorem 9.4 and Lemma 9.6).  
These were only stated for diagonal single-speed
systems.  Since the proofs relied only on energy estimates, the
results extend to the current setting provided \eqref{1.6} and
\eqref{1.7} hold. 

Prior to setting up the continuity argument, it is convenient to
reduce to an equivalent system of nonlinear equations with vanishing
Cauchy data.  By doing so, we will avoid complications related to the
compatibility conditions.  We first
reduce to an equivalent system of nonlinear equations whose data
vanish in a neighborhood of the obstacle.  Initially, we note that if
$\varepsilon$ in \eqref{1.11} is sufficiently small, then there is a
constant $C$ so that
\begin{equation}\label{5.1}
\sup_{0\le t\le 2} \sum_{|\alpha|\le 322} \|\partial^\alpha
u(t,\cd)\|_{L^2(|x|\le 10)} \le C\varepsilon.
\end{equation}
This, again, follows from the local existence theory (see, e.g.,
\cite{KSS}).  On the other hand, over $\{t\in [0,2]\}\times \{|x|\ge
6\}$, by finite propagation speed, $u$ corresponds to a solution of the boundaryless wave equation
$\Box u = F(u,du,d^2u)$.  If we take $N=322$ in \eqref{1.13}, it is clear
that the analogs of \eqref{3.4} and \eqref{3.6} yield
\begin{equation}\label{5.2}
\sup_{0\le t\le 2}\sum_{|\alpha|+\mu\le 321}\|L^\mu Z^\alpha
u'(t,\cd)\|_{L^2(|x|\ge 6)} + \sup_{\substack{0\le t\le 2\\|x|\ge 6}}
(1+t+|x|)\sum_{|\alpha|+\mu\le 311} |L^\mu Z^\alpha u(t,x)| \le C\varepsilon.
\end{equation} 
Here we have used our assumption that $\mathcal{K}\subset\{|x|<1\}$.

We will use this local solution to set up our reduction.  First, we
fix a cutoff function $\eta\in C^\infty(\R\times\R^3)$ satisfying
$\eta(t,x)\equiv 1$ if $t\le 3/2$ and $|x|\le 6$, $\eta(t,\cd)\equiv 0$ for
$t>2$, and $\eta(\cd, x)\equiv 0$ for $|x|>8$.  If we set 
$$u_0=\eta u,$$
it follows that $\Box u_0 = \eta F(u,du,d^2u) + [\Box, \eta] u.$
Thus, $u$ solves \eqref{1.1} for $0<t<T$ if and only if $w=u-u_0$ solves
\begin{equation}\label{5.3}
\begin{cases}
\Box w=(1-\eta)F(u,du,d^2u) - [\Box, \eta] u\\
w|_{\partial\mathcal{K}}=0\\
w(0,\cd)=(1-\eta)(0,\cd)f\\
\partial_t w(0,\cd)=(1-\eta)(0,\cd) g -\eta_t(0,\cd)f 
\end{cases}
\end{equation}
for $0<t<T$.

We now fix a smooth cutoff function $\beta$ with $\beta(t)\equiv 1$
for $t\le 1$ and $\beta(t)\equiv 0$ for $t\ge 3/2$.
If we let $v$ be the solution of the linear equation
\begin{equation}\label{5.4}
\begin{cases}
\Box v=\beta(1-\eta)F(u,du,d^2u)-[\Box,\eta] u\\
v|_{\partial\mathcal{K}}=0\\
v(0,\cd)=(1-\eta)(0,\cd)f\\
\partial_t v(0,\cd)=(1-\eta)(0,\cd)g-\eta_t(0,\cd)f,
\end{cases}
\end{equation}
we will show that there is an absolute constant so that
\begin{multline}\label{5.5}
(1+t+|x|)\sum_{\mu+|\alpha|\le 302} |L^\mu Z^\alpha v(t,x)| +
\sum_{\mu+|\alpha|\le 300} \|L^\mu Z^\alpha v'(t,\cd)\|_2
\\+(\log(2+t))^{-1} \sum_{\mu+|\alpha|\le 298} \|\langle
x\rangle^{-1/2} L^\mu Z^\alpha v'\|_{L^2(S_t)}
\le C_2\varepsilon
\end{multline}
where, as above, $S_t=[0,t]\times\ext$ denotes the time strip of
height $t$.

Indeed, by \eqref{4.15}, the first term on the left side of \eqref{5.5}
is bounded by
\begin{multline}\label{5.6}
C\sum_{\substack{j+|\beta|+k\le 310 \\ j\le 1}}
\|\langle x\rangle^{j+|\beta|}\partial_x^\beta \partial_t^{k+j}v(0,\cd)\|_2\\
+C\int_0^t \int \sum_{|\alpha|+\mu\le 309} |L^\mu Z^\alpha
\left(\beta(s)(1-\eta)(s,y) F(u,du,d^2u)(s,y)\right)|\:\frac{dy\:ds}{|y|}
\\+C\int_0^t \int \sum_{|\beta|+\mu\le 309} |L^\mu Z^\beta
[\Box,\eta]u|\:\frac{dy\:ds}{|y|}
\\+C\int_0^t \sum_{|\beta|+\mu \le 306} \|L^\mu \partial^\beta
[\Box,\eta]u\|_{L^2(\{x\in\ext\,:\,|x|<2\})} \:ds.
\end{multline}
It follows from \eqref{1.11} that the first term in \eqref{5.6}
is $O(\varepsilon)$.  Since $[\Box,\eta]u$ vanishes unless $t\le
2$ and $|x|\le 8$, the last two terms in \eqref{5.6} are also
$O(\varepsilon)$ by \eqref{5.1}.  
Thus, it remains to study the second term in \eqref{5.6}.  
This term is bounded by
\begin{multline*}
C\int_0^{3/2} \int_{|y|\ge 6} \sum_{|\alpha|+\mu\le 310} |L^\mu Z^\alpha
u'(s,y)|^2\:\frac{dy\:ds}{|y|} \\+ C\int_0^{3/2} \int_{|y|\ge 6}
\sum_{|\alpha|+\mu\le 311} |L^\mu Z^\alpha u(s,y)|^3
\:\frac{dy\:ds}{|y|}.
\end{multline*}
This is also clearly $O(\varepsilon)$ by \eqref{5.2}.

For the second term on the left of \eqref{5.5}, we use the standard
energy integral method  (see, e.g., Sogge \cite{S},
p.12) to see that
\begin{multline*}
\partial_t \sum_{|\alpha|+\mu\le 300}\|L^\mu Z^\alpha v'(t,\cd)\|_2^2
\\\le C\Bigl(\sum_{|\alpha|+\mu\le 300} \|L^\mu Z^\alpha v'(t,\cd)\|_2
\Bigr)\Bigl(\sum_{|\alpha|+\mu \le 300} \|L^\mu Z^\alpha \Box
v(t,\cd)\|_2\Bigr)
\\+C\sum_{|\alpha|+\mu\le 300}
\Bigl|\int_{\partial\mathcal{K}}\partial_0 L^\mu Z^\alpha
v(t,\cd)\nabla L^\mu Z^\alpha v(t,\cd)\cd {\mathbf{n}}\:d\sigma\Bigr|,
\end{multline*}
where ${\mathbf {n}}$ is the outward normal at a given point on
$\partial\mathcal{K}$.  Since $\mathcal{K}\subset\{|x|<1\}$ and since
$\Box v=\beta(t)(1-\eta)\Box u -[\Box,\eta]u$, it follows that
\begin{multline}\label{5.7}
\sum_{\mu+|\alpha|\le 300} \|L^\mu Z^\alpha v'(t,\cd)\|_2^2 \le
C\sum_{|\alpha|+\mu\le 300} \|L^\mu Z^\alpha v'(0,\cd)\|_2^2
\\+C\Bigl(\int_0^t \sum_{|\alpha|+\mu\le 300} \|L^\mu Z^\alpha
\beta(s)(1-\eta)(s,\cd) F(u,du,d^2u)(s,\cd)\|_2\:ds\Bigr)^2
\\+C\Bigl(\int_0^t \sum_{|\alpha|+\mu\le 300} \|L^\mu Z^\alpha
(-[\Box,\eta]u)(s,y)\|_2\:ds\Bigr)^2 
\\+ C\int_0^t \sum_{|\alpha|+\mu \le 301} \|L^\mu \partial^\alpha v'(s,\cd)\|^2_{L^2(|x|<1)}\:ds.
\end{multline}
The first term is $O(\varepsilon)$ by \eqref{1.11}.  
Since $[\Box,\eta]u$ is compactly supported in both $t$ and $x$, the third 
term in the right of \eqref{5.7} is also $O(\varepsilon)$ by \eqref{5.1}. 
Using the bound that we just obtained for the first term in the left of \eqref{5.5}, it
follows that the last term in \eqref{5.7} also satisfies the desired
bound.  We are left with studying the second term in \eqref{5.7}.  
This is clearly controlled by
\begin{multline*}
C\Bigl(\int_0^{3/2} \sum_{|\alpha|+\mu\le 301} \||L^\mu Z^\alpha
u'(s,\cd)|^2\|_{L^2(|x|>6)}\:ds\Bigr)^2
\\+ C\Bigl(\int_0^{3/2} \sum_{|\alpha|+\mu\le 302} \||L^\mu Z^\alpha
u(s,\cd)|^3\|_{L^2(|x|>6)}\:ds\Bigr)^2.
\end{multline*}
These terms are also easily seen to be $O(\varepsilon)$ by
\eqref{5.2}, which establishes the estimate for the second term in \eqref{5.5}.

Finally, it remains to show that the third term on the left side of
\eqref{5.5} is $O(\varepsilon)$.  To do so, we first notice that by
\eqref{4.37} we have
\begin{multline}\label{5.8}
r\langle c_I t-r\rangle^{1/2} |\partial L^\mu Z^\alpha v^I(t,x)| \le
C\sum_{|\beta|+\nu\le 300} \|L^\nu
Z^\beta v'(t,\cd)\|_2 
\\+C\sum_{|\beta|+\nu\le 299}
\|\langle t+r\rangle L^\nu Z^\beta (\partial_t^2-c_I^2 \Delta)v^I(t,\cd)\|_2 +
C(1+t)\sum_{\nu\le 298} \|L^\nu v'(t,\cd)\|_{L^\infty(|x|<2)}.
\end{multline}
for $\mu+|\alpha|\le 298$.
The first and last term on the right side of \eqref{5.8} are clearly
$O(\varepsilon)$ by the bounds for the first two terms in the left
side of \eqref{5.5}.  Since $\Box v = \beta(1-\eta)\Box u-[\Box,\eta]u$, the
second term on the right of \eqref{5.8} is controlled by
\begin{multline*}
C\sum_{|\beta|+\nu\le 300} \sup_{0\le t\le 3/2}\|\langle r\rangle |L^\nu Z^\beta
u'(t,\cd)|^2\|_{L^2(|x|>6)} \\
+ C\sum_{|\beta|+\nu\le 301}\sup_{0\le t\le 3/2} 
\|\langle r\rangle|L^\nu Z^\beta u(t,\cd)|^3\|_{L^2(|x|> 6)} 
+C\sum_{|\alpha|\le 300}\sup_{0\le t\le 2}
\|\partial_{t,x}^\alpha u\|_{L^2(|x|\le 8)}.
\end{multline*}
This is also $O(\varepsilon)$ by \eqref{5.1} and \eqref{5.2}.  
Thus, we have
\begin{equation}\label{5.9}
\sum_{\mu+|\alpha|\le 298} 
r\langle c_I t-r\rangle^{1/2} |\partial L^\mu Z^\alpha v^I(t,x)|\le C\varepsilon.
\end{equation}

In order to use this to bound the last term on the left of
\eqref{5.5}, notice that we can write
\begin{multline}\label{5.10}
\sum_{\mu+|\alpha|\le 298} \|\langle
x\rangle^{-1/2} L^\mu Z^\alpha v'\|^2_{L^2(S_t)} \le
C\sum_{\mu+|\alpha|\le 298} \int_0^t \frac{1}{1+s} \|L^\mu Z^\alpha
v'(s,\cd)\|^2_{L^2(|x|\ge c_1s/2)} \:ds \\+ C\sum_{\mu+|\alpha|\le 298}
\int_0^t \|\langle x\rangle^{-1/2} L^\mu Z^\alpha
v'(s,\cd)\|^2_{L^2(|x|\le c_1 s/2)}\:ds.
\end{multline}
By the bound for the second term on the left side of \eqref{5.5}, the
first term in \eqref{5.10} is clearly controlled by $C\varepsilon^2 \log(2+t)$.
If we apply \eqref{5.9} to the second term in \eqref{5.10}, assuming
as in \S 4 that the wavespeeds satisfy $0<c_1<c_2<\dots<c_D$, we see
that it is controlled by
$$C\varepsilon^2 \int_0^t \frac{1}{1+s} \|\langle
x\rangle^{-3/2}\|_{L^2(|x|\le c_1 s/2)}^2\:ds.$$
This is easily seen to be bounded by $C\varepsilon^2
(\log(2+t))^2$, which completes the proof of \eqref{5.5}.

The bounds \eqref{5.5} will allow us in many instances to restrict our
study to $w-v$
which is the solution of
\begin{equation}\label{5.11}
\begin{cases}
\Box (w-v)=(1-\beta)(1-\eta)F(u,du,d^2u), \quad (t,x)\in
\R_+\times\ext\\
(w-v)(t,x)=0,\quad x\in {\partial\mathcal{K}}\\
(w-v)(t,x)=0,\quad t\le 0.
\end{cases}
\end{equation}
Here, as mentioned earlier, 
we have vanishing Cauchy data, which allows us to avoid technical
details involving the compatibility conditions.

Depending on the linear estimates we employ, at times we shall use
certain $L^2$ and $L^\infty$ bounds for $u$ while at other times we
shall use them for $w-v$ or $w$.  Since $u=(w-v)+v+u_0$ and $u_0,v$
satisfy the bounds \eqref{5.1}, \eqref{5.5} respectively, it will
always be the case that bounds for $w-v$ will imply those for $w$
which in turn imply the same bounds for $u$ and vice versa.

We are now ready to set up the continuity argument.  If
$\varepsilon>0$ is as above, we shall assume that we have a
solution of our equation \eqref{1.1} for $0\le t\le T$ satisfying the
following dispersive estimates
\begin{equation}\label{5.12}
(1+t+|x|)\sum_{|\alpha|\le 201} |Z^\alpha w^I(t,x)|\le A_0 \varepsilon
\Bigl(1+\log\frac{1+t+|x|}{1+|c_I t-|x||}\Bigr)
\end{equation}
\begin{equation}\label{5.13}
(1+t+|x|)\sum_{\substack{|\alpha|+\nu\le 190\\\nu\le M}} |L^\nu
Z^\alpha w^I(t,x)|\le A_1\varepsilon
(1+t)^{\tilde{b}_{M+1}\varepsilon}\Bigl(1+\log\frac{1+t+|x|}{1+|c_I t-|x||}\Bigr)
\end{equation}
\begin{equation}\label{5.14}
(1+t+|x|)\sum_{\substack{|\alpha|+\nu\le 255\\\nu\le M}}|L^\nu Z^\alpha
w^I(t,x)|\le A_2\varepsilon (1+t)^{b_{M+1}\varepsilon}
\Bigl(1+\log\frac{1+t+|x|}{1+|c_It-|x||}\Bigr)
\end{equation}
\begin{equation}\label{5.15}
(1+|x|)\sum_{\substack{|\alpha|+\nu\le 180\\\nu\le N}} |L^\nu Z^\alpha
w^I(t,x)|\le A_3 \varepsilon (1+t)^{\tilde{c}_N \varepsilon}
\Bigl(1+\log\frac{1+t+|x|}{1+|c_It-|x||}\Bigr)
\end{equation}
\begin{equation}\label{5.14b}
(1+|x|)\sum_{\substack{|\alpha|+\nu\le 255\\\nu\le N}} |L^\nu Z^\alpha
w^I(t,x)|\le A_4 \varepsilon (1+t)^{c'_N \varepsilon}
\Bigl(1+\log\frac{1+t+|x|}{1+|c_It-|x||}\Bigr)
\end{equation}
\begin{equation}\label{5.16}
(1+t+|x|)\sum_{|\alpha|\le 200} |Z^\alpha w'(t,x)|\le
B_1\varepsilon
\end{equation}
for $M=0,1,2$ and $N=0,1,2,3$, and the following energy estimates
\begin{equation}\label{5.17}
\sum_{\substack{|\alpha|+\nu\le 220\\\nu\le 1}}\|L^\nu Z^\alpha
w'(t,\cd)\|_2\le A_5\varepsilon
\end{equation}
\begin{equation}\label{5.18}
\sum_{|\alpha|\le 300} \|\partial^\alpha u'(t,\cd)\|_2 \le
B_2\varepsilon (1+t)^{\tilde{C}\varepsilon}
\end{equation}
\begin{equation}\label{5.19}
\sum_{\substack{|\alpha|+\nu\le 202\\\nu\le N}} \|L^\nu Z^\alpha
u'(t,\cd)\|_2 + \sum_{\substack{|\alpha|+\nu\le
201\\\nu\le N}} \|\langle x\rangle^{-1/2} L^\nu Z^\alpha
u'\|_{L^2(S_t)} \le B_3 \varepsilon(1+t)^{\tilde{a}_N\varepsilon}
\end{equation}
\begin{equation}\label{5.20}
\sum_{\substack{|\alpha|+\nu\le 297-8N\\\nu\le N}} \|L^\nu Z^\alpha
u'(t,\cd)\|_2 + \sum_{\substack{|\alpha|+\nu\le 295-8N\\\nu\le N}} \|\langle
x\rangle^{-1/2} L^\nu Z^\alpha u'\|_{L^2(S_t)} \le B_4 \varepsilon(1+t)^{a_N \varepsilon}.
\end{equation}
As before, the $L^2_x$ norms are taken over $\ext$, and the
weighted $L^2_tL^2_x$-norms are taken over $S_t=[0,t]\times\ext$.

In \eqref{5.18}, $\tilde{C}$ is independent of the losses
$\tilde{a}_M$, $\tilde{b}_M$, $\tilde{c}_M$, $a_M$, $b_M$, and
$c_M'$.  The other associated losses satisfy
\begin{equation}\label{losses}
\tilde{a}_M\ll \tilde{b}_M\ll \tilde{c}_M\ll a_M\ll b_M\ll
c'_M\ll \tilde{a}_{M+1}
\end{equation}
for $M=1,2,3$, and
$$\tilde{C}\ll \tilde{a}_0 =\tilde{c}_0 = a_0\ll c'_0\ll \tilde{a}_1.$$


It is worth noting that \eqref{5.14}, \eqref{5.16}, \eqref{5.17},
\eqref{5.18}, and \eqref{5.20} are the estimates that made up the
simpler argument in the preceding paper \cite{MNS}.  \eqref{5.12} is
the main new estimate required in order to handle the higher order
terms that do not involve derivatives.  The remaining estimates are
technical pieces that are needed (or convenient) to make the argument
work.

In the estimates \eqref{5.12}-\eqref{5.14b} and \eqref{5.17}, we take
$A_j=4C_2$ where $j=0,1,\dots,5$ and $C_2$ is the uniform constant
appearing in the bounds \eqref{5.5} for $v$.  If $\varespilon$ is
small, all of these estimates are valid for $T=2$ by Theorem
\ref{theorem5.1}.  With this in mind, we shall prove that for
$\varepsilon >0$ sufficiently small depending on $B_1,\dots,B_4$
\begin{enumerate}
\item [$(i.)$] \eqref{5.12}-\eqref{5.14b} and \eqref{5.17} are valid with
$A_j$ replaced by $A_j/2$;
\item [$(ii.)$] \eqref{5.16} and \eqref{5.18}-\eqref{5.20} are a
consequence of \eqref{5.12}-\eqref{5.14b} and \eqref{5.17} for suitable
constants $B_j$.
\end{enumerate}
By the local existence theorem, it will follow that a solution exists
for all $t>0$ if $\varepsilon>0$ is sufficiently small.  We now
explore $(i.)$ and $(ii.)$ in the next two sections respectively.

\newsection{Proof of $(i.)$}

In this section, we will show step $(i.)$ of the proof of Theorem
\ref{theorem1.1}.  Specifically, we must show that
\eqref{5.12}-\eqref{5.14b} and \eqref{5.17} hold with $A_j$ replaced by $A_j/2$
under the assumption of \eqref{5.12}-\eqref{5.20}.

\subsection{Preliminaries:}
We begin with some preliminary estimates 
that follow from \eqref{5.12}-\eqref{5.20}.

First, we shall prove that if $|\alpha|+\nu\le 270$, $\nu\le 2$, 
then there is a constant $\tilde{b}$ so that
\begin{equation}\label{6.1}
\langle
r\rangle^{1/2+\theta} \langle c_I t-r\rangle^{1-\theta} |L^\nu Z^\alpha
\partial u^I (t,x)|\le C\varepsilon (1+t)^{\tilde{b}\varepsilon}, \quad 0\le \theta \le 1/2
\end{equation}
for any $0\le \theta\le 1/2$.  Additionally,
\begin{equation}\label{6.2}
\|\langle t+r\rangle
L^\nu Z^\alpha \Box u(t,\cd)\|_2\le C\varepsilon (1+t)^{\tilde{b}\varepsilon}.
\end{equation}
for $|\alpha|+\nu\le 271$ and $\nu\le 2$.

By \eqref{4.37}, \eqref{6.1} follows from \eqref{5.13}, \eqref{5.20}, and
\eqref{6.2}.  It, thus, suffices to show \eqref{6.2}.  To do so, notice
that the left side can be controlled by
\begin{multline}\label{6.3}
C\Bigl\|\Bigl(\langle t+r\rangle \sum_{\substack{|\alpha|+\nu\le
190\\\nu\le 2}} |L^\nu Z^\alpha
u'(t,\cd)|\Bigr)\sum_{\substack{|\alpha|+\nu\le 272\\\nu\le 2}}|L^\nu
Z^\alpha u'(t,\cd)|\Bigr\|_2
\\+ C\Bigl\|\langle t+r\rangle \Bigl(\sum_{\substack{|\alpha|+\nu\le
190\\\nu\le 2}} |L^\nu Z^\alpha u(t,\cd)|\Bigr)^2 \sum_{\substack{|\alpha|+\nu
\le 273\\\nu\le 2}} |L^\nu Z^\alpha u(t,\cd)|\Bigr\|_2.
\end{multline}
For the first term, if we apply \eqref{5.13}, we establish the bound
$$C\varepsilon (1+t)^{\tilde{b}_{3}\varepsilon}
(1+\log(1+t))
\sum_{\substack{|\alpha|+\nu\le 272\\\nu\le 2}} \|L^\nu Z^\alpha
u'(t,\cd)\|_2.$$
The desired estimate for the first term, thus, follows from \eqref{5.20}.

For the second term in \eqref{6.3}, we will again apply \eqref{5.13}.
Since the coefficients of $\Gamma=\{L,Z\}$ are $O(t+r)$, it follows
that this term is controlled by
$$C\varepsilon^2 (\log(2+t))^2 (1+t)^{2\tilde{b}_{3}\varepsilon}
\Bigl[\sum_{\substack{|\alpha|+\nu\le 272\\\nu\le 2}}\|L^\nu Z^\alpha
u'(t,\cd)\|_2 + \|\langle t+r\rangle^{-1} u(t,\cd)\|_2\Bigr].$$
The desired bound then follows from \eqref{5.12} and \eqref{5.20}, thus
completing the proof of \eqref{6.2}.

We will argue similarly to show a lossless version of \eqref{6.1} and
\eqref{6.2} that does not involve the scaling vector field $L$.  In
particular, we shall prove, for $|\alpha|\le 218$ and any $0\le
\theta\le 1/2$,
\begin{equation}\label{6.4}
r^{(1/2)+\theta}\langle c_I t-r\rangle^{1-\theta} |Z^\alpha \partial u^I(t,\cd)|\le C\varepsilon,
\end{equation}
and for $|\alpha|\le 219$,
\begin{equation}\label{6.5}
\|\langle t+r\rangle Z^\alpha \Box u(t,\cd)\|_2\le C\varepsilon.
\end{equation}
As before, \eqref{6.4} follows from \eqref{6.5} by \eqref{4.37},
\eqref{5.16}, and \eqref{5.17}.  

In order to show \eqref{6.5}, we again expand the left side to get the
bound
\begin{multline}\label{6.6}
C\Bigl\|\langle t+r\rangle \sum_{|\alpha|\le 190} |Z^\alpha u'(t,\cd)|
\sum_{|\alpha|\le 220} |Z^\alpha u'(t,\cd)|\Bigr\|_2
\\+C \Bigl\|\langle t+r\rangle \Bigl(\sum_{|\alpha|\le 190}|Z^\alpha
u(t,\cd)|\Bigr)^2 \sum_{|\alpha|\le 221} |Z^\alpha u(t,\cd)|\Bigr\|_2.
\end{multline}
By \eqref{5.1}, \eqref{5.16} and \eqref{5.17}, the first term is
$O(\varepsilon)$ as desired.  
Applying \eqref{5.1} and \eqref{5.12} to
the second term in \eqref{6.6}, we see that it is dominated by
$$
C\varepsilon^2 (1+\log(1+t))^2 \sum_{|\alpha|\le 221} \|\langle t+|x|\rangle^{-1}
|Z^\alpha u(t,\cd)|\|_2.
$$

\noindent By \eqref{5.1} and \eqref{5.14}, it follows that this term
is $O(\varepsilon)$ if $\varepsilon>0$ is sufficiently small.  This
completes the proof of \eqref{6.5}.

Notice that \eqref{6.1} and \eqref{6.4} hold when $u$ is replaced by
$w-v$.  Indeed, since $|L^\mu Z^\alpha \Box (w-v)|\lesssim
\sum_{\substack{|\beta|+\nu\le|\alpha|+\mu\\\nu\le \mu}}|L^\nu
Z^\beta \Box u|$, \eqref{6.2} and \eqref{6.5} hold with $w-v$
substituted for $u$.  Thus, the appropriate versions of \eqref{6.1}
and \eqref{6.4} are consequences of \eqref{4.37}, \eqref{5.1},
\eqref{5.5}, \eqref{5.13}, \eqref{5.16}, \eqref{5.17} and \eqref{5.20}.

\subsection{Proof of \eqref{5.12}:}  
Assuming \eqref{5.12}-\eqref{5.20}, we must show that \eqref{5.12}
holds with $A_0$ replaced by $A_0/2$.  Since the better bounds
\eqref{5.5} hold for $v$, it will suffice to show
\begin{equation}\label{a.1}
(1+t+|x|)\sum_{|\alpha|\le 201} |Z^\alpha (w-v)^I(t,x)|\le
C\varepsilon^2 \Bigl(1+\log\frac{1+t+|x|}{1+|c_It-|x||}\Bigr).
\end{equation}

Fix a smooth cutoff function $\eta_J$ satisfying $\eta_J(s)\equiv 1$
for $s\in [(c_J+(\delta/2))^{-1}, (c_J-(\delta/2))^{-1}]$ with $\delta
= (1/3)\min (c_I-c_{I-1})$, and $\eta_J(s)\equiv 0$ for $s\not\in
 [(c_J+\delta)^{-1}, (c_J-\delta)^{-1}]$.  Then, set
$\rho_J(t,x)=\eta_J(|x|^{-1}t)$.  Since we may assume that
$0\in\mathcal{K}$, we have that $|x|$ is bounded below on the
complement of $\mathcal{K}$, and the function $\rho_J$ is smooth and
homogeneous of degree $0$ in $(t,x)$.  Clearly, $\rho_J$ is
identically one on a conic neighborhood of $\{|x|=c_Jt\}$, and its
support does not intersect any $\{|x|=c_It\}$ for $I\neq
J$.  Let
\begin{equation}\label{a.2}
\tilde{F}^I=\sum_{\substack{1\le J\le D\\J\neq I}} \rho_J \sum_{0\le
j,k,l\le 3} B^{IJ,jk}_{J,l}\partial_l u^J \partial_j\partial_k u^J +
\sum_{\substack{1\le J\le D\\J\neq I}} \rho_J \sum_{0\le j,k\le 3}
A^{I,jk}_{JJ}\partial_j u^J \partial_k u^J,
\end{equation}
and set $\tilde{G}^I=F^I-\tilde{F}^I$.

By \eqref{4.24} and \eqref{4.32}, we have that the left side of
\eqref{a.1} is bounded by
\begin{multline}\label{a.3}
C\int_{\max (0,c|c_It-|x||-1)}^t \int_{|y|\approx s}
\sum_{\substack{|\alpha|+\nu\le 204\\\nu\le 1}} |L^\nu Z^\alpha
\tilde{F}^I(s,y)|\:\frac{dy\:ds}{|y|} 
\\+C\sum_{|\alpha|\le 207} \sup_{0\le s\le t} \int |Z^\alpha
\tilde{F}^I(s,y)|\:dy
\\+C\sup_{0\le s\le t} \int_{c's}^s \int_{|y|\approx \tau}
\sum_{\substack{|\alpha|+\nu\le 208\\\nu\le 1}} |L^\nu Z^\alpha
\tilde{F}^I(\tau,y)|\:\frac{dy\:d\tau}{|y|}
\\+C\sup_{0\le s\le t} (1+s) \sum_{|\alpha|\le 204} \|\partial^\alpha
\tilde{F}^I(s,\cd)\|_{L^\infty (|x|<10)}
\\+C\Bigl(1+\log\frac{1+t+|x|}{1+|c_It-|x||}\Bigr)\sup_{(s,y)}
|y|(1+s+|y|)^{1+\mu}z^{1-\mu}(s,|y|) \sum_{|\alpha|\le 201} |Z^\alpha
\tilde{G}^I(s,y)|
\\+C\Bigl(1+\log\frac{1+t+|x|}{1+|c_It-|x||}\Bigr)\sup_{(s,y)}
|y|(1+s+|y|)^{1+\mu}z^{1-\mu}(s,|y|)\sum_{|\alpha|\le 204}
|\partial^\alpha \partial \tilde{G}^I(s,y)|.
\end{multline}
We need to show that each of these terms is bounded by
$C\varepsilon^2\Bigl(1+\log\frac{1+t+|x|}{1+|c_It-|x||}\Bigr)$.

For the first term in \eqref{a.3}, it follows immediately that we have
the bound
\begin{multline*}
C\int_{\max(0,c|c_It-|x||-1)}^t \frac{1}{1+s} \int_{|y|\approx s}
\sum_{\substack{|\alpha|+\nu\le 205\\\nu\le 1}} |L^\nu Z^\alpha
\partial u^J|^2\:dy\:ds
\\\le C\Bigl(1+\log\frac{1+t}{1+|c_It-|x||}\Bigr) \sup_{0\le s\le t}
\sum_{\substack{|\alpha|+\nu\le 205\\\nu\le 1}} \|L^\nu Z^\alpha u'(s,\cd)\|_2^2.
\end{multline*}
The desired bound follows from \eqref{5.1} and \eqref{5.17}.  The
third term in \eqref{a.3} can be handled quite similarly.  The second
term above is easily seen to be $O(\varepsilon^2)$ by the Schwarz
inequality, \eqref{5.1}, and \eqref{5.17}.  The fourth term above is
bounded by
$$\sup_{0\le s\le t} (1+t)\Bigl\|\sum_{|\alpha|\le 103}
|\partial^\alpha \partial u^J| \sum_{|\alpha|\le 205} |\partial^\alpha
\partial u^J|\Bigr\|_\infty.$$
If we apply \eqref{5.1} and \eqref{5.16}, this is controlled by 
$$C\varepsilon \sup_{0\le s\le t} \sum_{|\alpha|\le 205}
\|\partial^\alpha u'\|_\infty.$$
Thus, by Sobolev's lemma, \eqref{5.1}, and \eqref{5.17}, we see that
this term is also $O(\varepsilon^2)$.

It remains to show that
\begin{multline}\label{a.4}
\sup_{(s,y)} |y|(1+s+|y|)^{1+\mu}z^{1-\mu}(s,|y|) \sum_{|\alpha|\le
201} |Z^\alpha \tilde{G}^I(s,y)|
\\+\sup_{(s,y)}|y|(1+s+|y|)^{1+\mu}z^{1-\mu}(s,|y|)\sum_{|\alpha|\le
204} |\partial^\alpha \partial \tilde{G}^I(s,y)|
\end{multline}
is $O(\varepsilon^2)$.

When $\tilde{G}^I$ is replaced by the null forms
$$\sum_{0\le j,k\le 3} A^{I,jk}_{II}\partial_j u^I \partial_k u^I +
\sum_{0\le j,k,l\le 3} B^{II,jk}_{I,l}\partial_l u^I
\partial_j\partial_k u^I,$$
we can apply \eqref{2.24} and \eqref{2.25} to see that \eqref{a.4} is
controlled by 
\begin{multline}\label{a.5}
C\sup_{(s,y)} (1+s+|y|)^{1+\mu} z^{1-\mu}(s,|y|)
\sum_{\substack{|\alpha|+\nu\le 207\\\nu\le 1}} |L^\nu Z^\alpha
u^I(s,y)| \sum_{|\alpha|\le 206} |Z^\alpha \partial u^I(s,y)|
\\+C\sup_{(s,y)} |y|(1+s+|y|)^\mu z^{1-\mu}(s,|y|)\langle
c_Is-|y|\rangle \sum_{|\alpha|\le 103} |Z^\alpha \partial u^I(s,y)|
\sum_{|\alpha|\le 206} |Z^\alpha \partial u^I(s,y)|.
\end{multline}
For the first term, if we apply \eqref{5.1} and \eqref{5.14}, we get the bound
$$C\varepsilon \sup_{(s,y)}
(1+s+|y|)^{\mu+b_2\varepsilon+}z^{1-\mu}(s,|y|) \sum_{|\alpha|\le 206}
|Z^\alpha \partial u^I|.$$
If $\mu$ and $\varepsilon$ are sufficiently small, the desired
$O(\varepsilon^2)$ bound follows from \eqref{6.1}.  Indeed, if
$(s,|y|)\in\Lambda_I$, then $z(s,|y|)=\langle c_Is-|y|\rangle$ and
$|y|\gtrsim (1+s+|y|)$.  On the other hand, if
$(s,|y|)\not\in\Lambda_I$, then $\langle c_Is-|y|\rangle\gtrsim
\langle s+|y|\rangle$ and $\langle y\rangle^{1-\mu}\gtrsim z^{1-\mu}(s,|y|)$.

For the second term in \eqref{a.5}, we apply \eqref{6.4} to obtain the
bound
$$C\varepsilon \sup_{(s,y)} |y|^{1/2} (1+s+|y|)^\mu z^{1-\mu}(s,|y|)
\sum_{|\alpha|\le 206} |Z^\alpha \partial u^I|.$$
Using considerations as above, this is $O(\varepsilon^2)$ by a
subsequent application of \eqref{6.4}.

When $\tilde{G}^I$ is replaced by
$$\sum_{\substack{1\le J\le D\\J\neq I}} \sum_{0\le j,k,l\le 3}
(1-\rho_J) B^{IJ,jk}_{J,l} \partial_l u^J \partial_j\partial_k u^J
+ \sum_{\substack{1\le J\le D\\J\neq I}} \sum_{0\le j,k\le 3}
(1-\rho_J) A^{I,jk}_{JJ} \partial_j u^J \partial_k u^J,$$
\eqref{a.4} is dominated by
$$C\sum_{\substack{1\le J\le D\\J\neq I}} \sup_{(s,y)\not\in
\Lambda_J} |y|(1+s+|y|)^{1+\mu} z^{1-\mu}(s,|y|)
\Bigl(\sum_{|\alpha|\le 206} |Z^\alpha \partial u^J|\Bigr)^2.$$
Since $\langle c_Js-|y|\rangle^{1-\mu}\gtrsim \langle
s+|y|\rangle^{1-\mu} \gtrsim z^{1-\mu}(s,|y|)$ on the support of
$(1-\rho_J)$, these terms are $O(\varepsilon^2)$ by two applications
of \eqref{6.4} with $\theta=0$.

If $\tilde{G}^I$ in \eqref{a.4} is replaced by the remaining quadratic
terms
$$\sum_{\substack{1\le J,K\le D\\J\neq K}} \sum_{0\le j,k,l\le 3}
B^{IJ,jk}_{K,l} \partial_l u^K \partial_j\partial_k u^J +
\sum_{\substack{1\le J,K\le D\\J\neq K}} \sum_{0\le j,k\le 3}
A^{I,jk}_{JK} \partial_j u^J \partial_k u^K,$$
we see that it is bounded by
\begin{equation}\label{a.6}
C\sum_{\substack{1\le J,K\le D\\J\neq K}} \sup_{(s,y)}
|y|(1+s+|y|)^{1+\mu} z^{1-\mu}(s,|y|) \sum_{|\alpha|\le 206} |Z^\alpha
\partial u^J| \sum_{|\alpha|\le 206} |Z^\alpha \partial u^K|.
\end{equation}
If $(s,y)\not\in \Lambda_J\cup\Lambda_K$, then we can argue as in the
previous case to see that this is $O(\varepsilon^2)$.  Thus, let us
assume that $(s,y)\in\Lambda_J$, and hence $(s,y)\not\in\Lambda_K$.
The reverse case will follow symmetrically.  For such $(s,y)$, we have
$z^{1-\mu}(s,|y|)=\langle c_Js-|y|\rangle^{1-\mu}$.  Thus, by
\eqref{6.4}, we see that in this case \eqref{a.6} is controlled by
$$C\varepsilon \sum_{\substack{1\le K\le D\\K\neq J}}
\sup_{(s,y)\in\Lambda_J}
|y|^{1/2-\mu}(1+s+|y|)^{1+\mu}\sum_{|\alpha|\le 206} |Z^\alpha \partial
u^K|.$$
Since $\langle c_Ks-|y|\rangle\gtrsim \langle s+|y|\rangle$ and
$|y|\approx s$ on $\Lambda_J$, the desired $O(\varepsilon^2)$ bound
follows from \eqref{6.4}.

Finally, when $\tilde{G}^I$ is replaced by the cubic terms $R^I+P^I$,
\eqref{a.4} is dominated by
\begin{multline}\label{a.7}
C\sup_{(s,y)}|y|(1+s+|y|)^{1+\mu}z^{1-\mu}(s,|y|)
\Bigl(\sum_{|\alpha|\le 201} |Z^\alpha u(s,y)|\Bigr)^3
\\+C\sup_{(s,y)}|y|(1+s+|y|)^{1+\mu}z^{1-\mu}(s,|y|)
\Bigl(\sum_{|\alpha|\le 104} |Z^\alpha
u(s,y)|\Bigr)^2\sum_{|\alpha|\le 206} |Z^\alpha \partial u(s,y)|.
\end{multline}
By the inductive hypothesis \eqref{5.12}, the first term in \eqref{a.7}
is controlled by
$$C\varepsilon^3 \sup_{(s,y)}
\frac{z^{1-\mu}(s,|y|)}{(1+s+|y|)^{1-\mu}}
\Bigl(1+\log\frac{1+s+|y|}{z(s,|y|)}\Bigr)^3.$$
Since $(\log x)^3/x^{1-\mu}$ is bounded for $x\ge 1$ and $\mu < 1$, it
follows that this term is $O(\varepsilon^3)$.  For the second term
in\eqref{a.7}, by \eqref{6.4}, we have the bound
$$C\varepsilon \sup_{(s,y)} |y|^{1/2-\mu}
(1+s+|y|)^{1+\mu}\Bigl(\sum_{|\alpha|\le 104} |Z^\alpha
u(s,y)|\Bigr)^2.$$
This term is then easily seen to be $O(\varepsilon^3)$ by \eqref{5.1}
and \eqref{5.12}
which completes the proof.

\subsection{Proof of \eqref{5.13}:}
In this section, we show that if you assume \eqref{5.12}-\eqref{5.20},
then you can prove \eqref{5.13} with $A_1$ replaced by $A_1/2$.  By
the arguments in the previous section, this clearly holds when $M=0$.
As before, by \eqref{5.5}, it suffices to show
\begin{multline}\label{b.1}
\Bigl(1+\log\frac{1+t+|x|}{1+|c_It-|x||}\Bigr)^{-1}
(1+t+|x|)\sum_{\substack{|\alpha|+\nu\le 190\\\nu\le M}} |L^\nu
Z^\alpha (w-v)^I(t,x)|
\\\le C\varepsilon^2 (1+t)^{\tilde{b}_{M+1}\varepsilon}.
\end{multline}

Since $\Box (w-v)=(1-\beta)(1-\eta)\Box u =
(1-\beta)(1-\eta)(B+Q+R+P)$, by \eqref{4.15} and \eqref{4.32}, we see
that the left side of \eqref{b.1} is dominated by
\begin{multline}\label{b.2}
C\int_0^t \int \sum_{\substack{|\alpha|+\nu\le 197\\\nu\le M+1}}
|L^\nu Z^\alpha (B^I+Q^I)(s,y)|\:\frac{dy\:ds}{|y|}
\\+C\sup_{(s,y)}|y|(1+s+|y|)^{1+\mu}z^{1-\mu}(s,|y|)\sum_{\substack{|\alpha|+\nu\le
190\\\nu\le M}} |L^\nu Z^\alpha (R^I+P^I)(s,y)|
\\+C\sup_{(s,y)}|y|(1+s+|y|)^{1+\mu}z^{1-\mu}(s,|y|)\sum_{\substack{|\alpha|+\nu\le
193\\\nu \le M}} |L^\nu \partial^\alpha \partial (R^I+P^I)(s,y)|.
\end{multline}
Here we have used the fact that the last term in \eqref{4.15} is
controlled by the second term in the right of \eqref{4.15} using
Sobolev estimates and the fact that we
may assume $0\in \mathcal{K}$ without loss of generality.

By \eqref{1.3} and \eqref{1.4}, we have that the
first term in \eqref{b.2} is dominated by
$$C\sum_{\substack{|\alpha|+\nu\le 198\\\nu\le M+1}} \|\langle
x\rangle^{-1/2} L^\nu Z^\alpha u'\|^2_{L^2(S_t)}.$$
It, thus, follows from \eqref{5.19} that these terms are bounded by
$C\varepsilon^2 (1+t)^{2\tilde{a}_{M+1}\varepsilon}$ as desired.

The last two terms of \eqref{b.2} are controlled by a constant times
\begin{multline}\label{b.3}
\sup_{(s,y)}|y|(1+s+|y|)^{1+\mu} z^{1-\mu}(s,|y|)
\Bigl(\sum_{|\alpha|\le 190} |Z^\alpha u(s,y)|\Bigr)^2
\sum_{\substack{|\alpha|+\nu\le 190\\\nu\le M}} |L^\nu Z^\alpha
u(s,y)|
\\+\sup_{(s,y)} |y|(1+s+|y|)^{1+\mu} z^{1-\mu}(s,|y|)
\Bigl(\sum_{\substack{|\alpha|+\nu\le 190\\\nu\le M-1}} |L^\nu
Z^\alpha u(s,y)|\Bigr)^3
\\+\sup_{(s,y)} |y|(1+s+|y|)^{1+\mu} z^{1-\mu}(s,|y|)
\Bigl(\sum_{\substack{|\alpha|+\nu\le 190\\\nu\le M}} |L^\nu Z^\alpha
u(s,y)|\Bigr)^2 \sum_{\substack{|\alpha|+\nu\le 195\\\nu\le M}} |L^\nu
Z^\alpha \partial u(s,y)|.
\end{multline}
For the first term in \eqref{b.3}, we apply \eqref{5.1}, \eqref{5.12}, and
\eqref{5.13} to obtain the bound
$$C\varepsilon^3 (1+t)^{\tilde{b}_{M+1}\varepsilon} \sup_{(s,y)}
\frac{z^{1-\mu}(s,|y|)}{(1+s+|y|)^{1-\mu}}
\Bigl(1+\log\frac{1+s+|y|}{z(s,|y|)}\Bigr)^3 \le C\varepsilon^3
(1+t)^{\tilde{b}_{M+1}\varepsilon}.$$
If we apply \eqref{5.13} and argue similarly, it follows that the
second term is controlled by $C\varepsilon^3
(1+t)^{3\tilde{b}_M\varepsilon}$.
Finally,
for the third term in \eqref{b.3}, we first apply \eqref{6.1} to see
that it is controlled by
$$C\varepsilon \sup_{(s,|y|)} |y|^{1/2-\mu} (1+s+|y|)^{1+\mu}
(1+s)^{\tilde{b}\varepsilon} \Bigl(\sum_{\substack{|\alpha|+\nu\le
190\\\nu\le M}} |L^\nu Z^\alpha u(s,y)|\Bigr)^2.$$
It, thus, follows from the inductive hypothesis \eqref{5.13} that this
term is $O(\varepsilon^3)$ if $\varepsilon>0$ is sufficiently small,
which completes the proof of \eqref{5.13}.

\subsection{Proof of \eqref{5.14}:}  
In this section, by proving
\begin{equation}\label{c.1}
\Bigl(1+\log\frac{1+t+|x|}{1+|c_It-|x||}\Bigr)^{-1}
(1+t+|x|)\sum_{\substack{|\alpha|+\nu\le 255\\\nu\le M}} |L^\nu
Z^\alpha (w-v)^I(t,x)|\le C\varepsilon^2 (1+t)^{b_{M+1}\varepsilon},
\end{equation}
we show that \eqref{5.14} holds with $A_2$ replaced by $A_2/2$ for $M=0,1,2.$

Here, again, we apply \eqref{4.15} and \eqref{4.32} to see that the
left side of \eqref{c.1} is controlled by
\begin{multline}\label{c.2}
C\sum_{\substack{|\alpha|+\nu\le 262\\\nu\le M+1}} \int_0^t \int
|L^\nu Z^\alpha (Q^I+B^I)(s,y)|\:\frac{dy\:ds}{|y|}
\\+C\sup_{(s,y)} |y|(1+s+|y|)^{1+\mu} z^{1-\mu}(s,|y|)
\sum_{\substack{|\alpha|+\nu\le 255\\\nu\le M}} |L^\nu Z^\alpha
(R^I+P^I)(s,y)|
\\+C\sup_{(s,y)}|y|(1+s+|y|)^{1+\mu}z^{1-\mu}(s,|y|)\sum_{\substack{|\alpha|+\nu\le
258\\\nu\le M}} |L^\nu \partial^\alpha \partial (R^I+P^I)(s,y)|.
\end{multline}

The first term is controlled by
$$C\sum_{\substack{|\alpha|+\nu\le 263\\\nu\le M+1}} \|\langle
x\rangle^{-1/2} L^\nu Z^\alpha u'\|^2_{L^2(S_t)}\le C\varepsilon^2
(1+t)^{2a_{M+1}\varepsilon}$$
by \eqref{5.20}.  For the last two terms in \eqref{c.2}, which involve
the cubic nonlinearities, we have the bound
\begin{multline}\label{c.3}
C\sup_{(s,y)}|y|(1+s+|y|)^{1+\mu}z^{1-\mu}(s,|y|)\Bigl(\sum_{|\alpha|\le
131} |Z^\alpha u(s,y)|\Bigr)^2 \sum_{\substack{|\alpha|+\nu\le
255\\\nu\le M}} |L^\nu Z^\alpha u(s,y)|
\\+C\sup_{(s,y)}|y|(1+s+|y|)^{1+\mu}z^{1-\mu}(s,|y|)
\Bigl(\sum_{\substack{|\alpha| + \nu\le 131\\\nu\le M}} |L^\nu
Z^\alpha u(s,y)|\Bigr)^2 \sum_{\substack{|\alpha|+\nu\le 255\\\nu\le
M-1}} |L^\nu Z^\alpha u(s,y)|
\\+C\sup_{(s,y)} |y|(1+s+|y|)^{1+\mu} z^{1-\mu}(s,|y|)
\Bigl(\sum_{\substack{|\alpha|+\nu \le 131\\\nu\le M}} |L^\nu Z^\alpha
u(s,y)|\Bigr)^2 \sum_{\substack{|\alpha|+\nu\le 260\\\nu\le M}} |L^\nu
Z^\alpha u'(s,y)|.
\end{multline}

Using the fact that $(\log x)^3/x^{1-\mu}$ is bounded for $x\ge 1$ and
$\mu<1$, the first term is dominated by $C\varepsilon^3
(1+t)^{b_{M+1}\varepsilon}$ by \eqref{5.12} and \eqref{5.14}.
Similarly, using \eqref{5.13} and \eqref{5.14}, the second term is
controlled by $C\varepsilon^3
(1+t)^{2\tilde{b}_{M+1}\varepsilon+b_M\varepsilon}$.  Again arguing as
in the proof of \eqref{5.13}, the final term in \eqref{c.3} is easily
seen to be $O(\varepsilon^3)$ if $\varepsilon>0$ is sufficiently small
 using \eqref{6.1} and \eqref{5.13}.  This completes the proof of
\eqref{c.1}, and hence, that of \eqref{5.14}.

\subsection{Proof of \eqref{5.15}:}
In the proof of part $(ii.)$, we will require pointwise estimates
that allow up to three occurences of the scaling vector field $L$.
This is not the case for the previous estimates due to the loss of an
$L$ associated to \eqref{4.15}.  Here, we may argue as in the proofs
of the previous esimates (in particular, that of \eqref{5.13}) 
replacing \eqref{4.15} by \eqref{4.20}.

Clearly \eqref{5.15} holds when $N=0$ by \eqref{5.12}.  Thus, by
\eqref{5.5}, in order to show that \eqref{5.15} holds with $A_3$
replaced by $A_3/2$, it suffices to show
\begin{equation}\label{d.1}
\Bigl(1+\log\frac{1+t+|x|}{1+|c_It-|x||}\Bigr)^{-1}
(1+|x|)\sum_{\substack{|\alpha|+\nu\le 180\\\nu\le N}} |L^\nu Z^\alpha
(w-v)^I(t,x)|\le C\varepsilon^2 (1+t)^{\tilde{c}_N\varepsilon}
\end{equation}
for $N=1,2,3$.

Since $\Box (w-v)=(1-\beta)(1-\eta)\Box u =
(1-\beta)(1-\eta)(B+Q+R+P)$, by \eqref{4.20} and \eqref{4.32}, we see
that the left side of \eqref{d.1} is controlled by
\begin{multline}\label{d.2}
C\int_0^t\int \sum_{\substack{|\alpha|+\nu\le 186\\\nu\le N}} |L^\nu
Z^\alpha (B^I+Q^I)(s,y)|\:\frac{dy\:ds}{|y|}
\\+C\sup_{(s,y)} |y|(1+s+|y|)^{1+\mu}z^{1-\mu}(s,|y|)
\sum_{\substack{|\alpha|+\nu\le 180\\\nu\le N}} |L^\nu Z^\alpha
(R^I+P^I)(s,y)|
\\+C\sup_{(s,y)}|y|(1+s+|y|)^{1+\mu} z^{1-\mu}(s,|y|)
\sum_{\substack{|\alpha|+\nu\le 183\\\nu\le N}} |L^\nu \partial^\alpha
\partial (R^I+P^I)(s,y)|.
\end{multline}
As above, using \eqref{5.19}, the first term is bounded by
$$C\sum_{\substack{|\alpha|+\nu\le 187\\\nu\le N}} \|\langle
x\rangle^{-1/2} L^\nu Z^\alpha u'\|^2_{L^2(S_t)}\le C\varepsilon^2
(1+t)^{2\tilde{a}_N\varepsilon}.$$

The cubic terms require a little additional care.  To begin, we have that
the last two terms of \eqref{d.2} are controlled by
\begin{multline}\label{d.3}
C\sup_{(s,y)}|y|(1+s+|y|)^{1+\mu}z^{1-\mu}(s,|y|)\Bigl(\sum_{|\alpha|\le
186} |Z^\alpha u(s,y)|\Bigr)^2 \sum_{\substack{|\alpha|+\nu\le
180\\\nu\le N}} |L^\nu Z^\alpha u(s,y)|
\\+C\sup_{(s,y)}|y|(1+s+|y|)^{1+\mu}
z^{1-\mu}(s,|y|)\Bigl(\sum_{|\alpha|\le 186} |Z^\alpha u(s,y)|\Bigr)^2
\sum_{\substack{|\alpha|+\nu\le 188\\\nu\le N}} |L^\nu Z^\alpha
u'(s,y)|
\\+C\sup_{(s,y)}|y|(1+s+|y|)^{1+\mu}z^{1-\mu}(s,|y|)
\Bigl(\sum_{\substack{|\alpha|+\nu\le 189\\\nu\le N-1}} |L^\nu
Z^\alpha u(s,y)|\Bigr)^3.
\end{multline}

Applying \eqref{5.1} and \eqref{5.15} to the first term and using \eqref{2.26}
and \eqref{5.19} in the second, we see that the first two terms of
\eqref{d.3} are controlled by
$$C\varepsilon (1+t)^{\tilde{c}_N\varepsilon}
\sup_{(s,y)} (1+s+|y|)^{1+\mu}z^{1-\mu}(s,|y|)
\Bigl(1+\log\frac{1+s+|y|}{z(s,|y|)}\Bigr)
\Bigl(\sum_{|\alpha|\le 186} |Z^\alpha u(s,y)|\Bigr)^2$$
Here, we have used that \eqref{losses} gives $\tilde{a}_N\le
\tilde{c}_N$.  If we in turn apply \eqref{5.12}, we see that this is
bounded by the right side of \eqref{d.1} as desired.  When $N=1$, this
is sufficient to complete the proof.  When $N=2,3$, we must also
consider the last term in \eqref{d.3}.  The bound here, however,
follows quite simply from three applications of \eqref{5.13} (and
\eqref{5.1}).  Doing so, we see that this last term is controlled by
$C\varepsilon^3 (1+t)^{3\tilde{b}_N\varepsilon}$.  Since we may choose
$\tilde{c}_N> 3 \tilde{b}_N$ (see \eqref{losses}), 
this is sufficient to complete the proof
of \eqref{d.1}.

\subsection{Proof of \eqref{5.14b}:}
In this section, we will argue much as in the previous section to
establish the higher order pointwise estimate that permits three
occurences of $L$.  Here, we must establish \eqref{5.14b} with $A_4$
replaced by $A_4/2$.  This is accomplished by showing
\begin{equation}\label{e.1}
\Bigl(1+\log\frac{1+t+|x|}{1+|c_It-|x||}\Bigr)^{-1}(1+|x|)
\sum_{\substack{|\alpha|+\nu\le 255\\\nu\le N}} |L^\nu Z^\alpha
(w-v)^I(t,x)| \le C\varepsilon^2 (1+t)^{c_N'\varepsilon}
\end{equation}
for $N=0,1,\dots,4$ and using \eqref{5.5}.

Applying \eqref{4.20} and \eqref{4.32} and arguing as in the proof of
\eqref{5.15}, we see that the left side of \eqref{e.1} is dominated by
\begin{multline}\label{e.2}
C\sum_{\substack{|\alpha|+\nu\le 262\\\nu\le N}} \|\langle
x\rangle^{-1/2} L^\nu Z^\alpha u'\|^2_{L^2(S_t)}
\\+C\sup_{(s,y)}|y|(1+s+|y|)^{1+\mu}z^{1-\mu}(s,|y|)
\Bigl(\sum_{|\alpha|\le 132} |Z^\alpha u(s,y)|\Bigr)^2
\sum_{\substack{|\alpha|+\nu\le 255\\\nu\le N}} |L^\nu Z^\alpha
u(s,y)|
\\+C\sup_{(s,y)}|y|(1+s+|y|)^{1+\mu}z^{1-\mu}(s,|y|)
\Bigl(\sum_{|\alpha|\le 132} |Z^\alpha u(s,y)|\Bigr)^2
\sum_{\substack{|\alpha|+\nu \le 260\\\nu\le N}} |L^\nu Z^\alpha
u'(s,y)|
\\+C\sup_{(s,y)}|y|(1+s+|y|)^{1+\mu} z^{1-\mu}(s,|y|)
\sum_{\substack{|\alpha|+\nu\le 132\\\nu\le N}} |L^\nu Z^\alpha
u(s,y)|\Bigl(\sum_{\substack{|\alpha|+\nu\le 255\\\nu\le N-1}} |L^\nu
Z^\alpha u(s,y)|\Bigr)^2
\\+C\sup_{(s,y)}|y|(1+s+|y|)^{1+\mu}z^{1-\mu} (s,|y|)
\sum_{\substack{|\alpha|+ \nu\le 132\\\nu\le N}} |L^\nu Z^\alpha
u(s,y)| 
\\\times\sum_{\substack{|\alpha|+\nu\le 132\\\nu\le N-1}} |L^\nu
Z^\alpha u(s,y)| \sum_{\substack{|\alpha|+\nu\le 260\\\nu\le N-1}}
|L^\nu Z^\alpha u'(s,y)|.
\end{multline}

Choosing $c_N'> \max(2a_N, 2b_N+\tilde{c}_N)$ as we may, we see that this is
bounded by the right side of \eqref{e.1}.  Indeed, the bound for the
first term in \eqref{e.2} follows directly from \eqref{5.20}.  For the
second term, we apply \eqref{5.1}, \eqref{5.12}, and \eqref{5.14b} as
before.  This suffices to handle the $N=0$ case.  In order to complete
the proof for $N=1,2,3$, we  similarly, bound the third term using
applications of \eqref{2.26}, \eqref{5.1}, \eqref{5.12}, and
\eqref{5.20}.  To get control over the fourth term, we apply
\eqref{5.1}, \eqref{5.14}, and \eqref{5.15}.  Using \eqref{5.1},
\eqref{5.13}, \eqref{5.15}, and \eqref{6.1}, one can see that the last
term is $O(\varepsilon^3)$ for small $\varepsilon$ which completes the
proof of \eqref{e.1}.

\subsection{Proof of \eqref{5.17}:}   In order to complete the proof
of part $(i.)$, it remains to show that the low order, lossless energy
inequality \eqref{5.17} with $A_5$ replaced by $A_5/2$ follows from
\eqref{5.12}-\eqref{5.20}.  Since $v$ satisfies the better bound
\eqref{5.5}, it suffices to establish
\begin{equation}\label{f.1}
\sum_{\substack{|\alpha|+\nu\le 220\\\nu\le 1}} \|L^\nu Z^\alpha
(w-v)'(t,\cd)\|_2^2 \le C\varepsilon^3.
\end{equation}

By the standard energy integral method, 
we have that the left side of \eqref{f.1} is bounded by
\begin{multline*}
C\sum_{\substack{|\alpha|+\nu\le 220\\\nu\le 1}} \int_0^t \int_{\ext}
|\langle \partial_0 L^\nu Z^\alpha (w-v),\Box L^\nu Z^\alpha
(w-v)\rangle|\:dy\:ds
\\+C\sum_{\substack{|\alpha|+\nu\le 220\\\nu\le 1}}
\Bigl|\int_0^t\int_{\partial\mathcal{K}} \partial_0L^\nu Z^\alpha
(w-v)\, \nabla_x L^\nu Z^\alpha (w-v)\cd \mathbf{n} \:d\sigma\:ds\Bigr|
\end{multline*}
where $\mathbf{n}$ is the outward normal at a given point on
$\mathcal{K}$ and $\langle \cd,\cd\rangle$ is the standard Euclidean
inner product on $\R^D$.  Since $\mathcal{K}\subset \{|x|<1\}$ and
since the coefficients of $Z$ are $O(1)$ on $\partial\mathcal{K}$, it
follows that the last term in controlled by
$$C\int_0^t \int_{\{x\in\ext\,:\,|x|<1\}} \sum_{\substack{|\alpha| +
\nu\le 221\\\nu\le 1}} |L^\nu \partial^\alpha
(w-v)'(s,y)|^2\:dy\:ds.$$
Additionally, by the commutation properties of $\Gamma$ with $\Box$
and the fact that $\Box(w-v)=(1-\beta)(1-\eta)\Box u$, we see that the
left side of \eqref{f.1} is dominated by
\begin{multline*}
C\int_0^t \int_\ext \sum_{\substack{|\alpha_j|+\nu_j\le 220\\\nu_j\le
1; j=1,2}} \Bigl|\langle\partial_0 L^{\nu_1}Z^{\alpha_1}
(w-v),L^{\nu_2}Z^{\alpha_2} F(u,du,d^2u)\rangle\Bigr|\:dy\:ds
\\+C\int_0^t\int_{\{x\in\ext\,:\,|x|<1\}}
\sum_{\substack{|\alpha|+\nu\le 221\\\nu\le 1}} |L^\nu \partial^\alpha
(w-v)'(s,y)|^2 \:dy\:ds.
\end{multline*}

If we expand using the definition of $F(u,du,d^2u)$, the preceding
equation is controlled by
\begin{multline}\label{f.2}
C\int_0^t\int_\ext \sum_{K=1}^D\sum_{\substack{|\alpha|+\nu\le
220\\\nu\le 1}} \Bigl|\partial_0L^\nu Z^\alpha (w-v)^K\Bigr|
\\\times \sum_{\substack{|\alpha|+\nu\le 220\\\nu\le 1}}
\sum_{\substack{|\beta|+\mu \le 220\\\mu\le 1}} \Bigl|\sum_{0\le
j,k,l\le 3} \tilde{B}^{KK,jk}_{K,l}\partial_l L^\nu Z^\alpha u^K
 \partial_j\partial_k L^\mu Z^\beta u^K\Bigr|\:dy\:ds
\\+C\int_0^t\int_\ext \sum_{K=1}^D \sum_{\substack{|\alpha|+\nu\le
220\\\nu\le 1}} \Bigl|\partial_0 L^\nu Z^\alpha (w-v)^K\Bigr|
\sum_{\substack{|\alpha|+ \nu\le 220\\\nu\le 1}}
\sum_{\substack{|\beta|+\mu\le 220\\\mu\le 1}} \Bigl|\sum_{0\le j,k\le
3} \tilde{A}^{K,jk}_{KK} \partial_j L^\nu Z^\alpha u^K 
\\\times \partial_k L^\mu Z^\beta u^K\Bigr|\:dy\:ds
\\+C\int_0^t\int_\ext \sum_{\substack{1\le I,J,K\le D\\(I,K)\neq
(K,J)}} \sum_{\substack{|\alpha|+\nu\le 220\\\nu\le 1}} |L^\nu
Z^\alpha \partial (w-v)^K|\sum_{\substack{|\alpha|+ \nu\le 220\\\nu\le
1}} |L^\nu Z^\alpha \partial u^I| \\\times \sum_{\substack{ |\alpha|+\nu \le
221\\\nu \le 1}} |L^\nu Z^\alpha\partial u^J|\:dy\:ds
\\+C\int_0^t \int_\ext \sum_{\substack{|\alpha|+\nu\le 220\\\nu\le 1}}
|L^\nu Z^\alpha \partial (w-v)|\Bigl(\sum_{\substack{|\alpha|+\nu\le
222\\\nu\le 1}} |L^\nu Z^\alpha u|\Bigr)^3 \:dy\:ds
\\+C\int_0^t\int_{\{x\in\ext\,:\,|x|<1\}}
\sum_{\substack{|\alpha|+\nu\le 221\\\nu\le 1}} |L^\nu \partial^\alpha
(w-v)'(s,y)|^2\:dy\:ds.
\end{multline}
By Lemma 4.1 of Sideris-Tu \cite{Si3}, the constants
$\tilde{A}^{K,jk}_{KK}$ and $\tilde{B}^{KK,jk}_{K,l}$ satisfy
\eqref{1.8} and \eqref{1.9}.

The first two terms in \eqref{f.2} satisfy the bounds of Lemma
\ref{lemma2.5}.  The third term involves quadratic interactions
between waves of different speeds, and the fourth term is the cumulative
effect of the nonlinearities of higher order.  The arguments to bound
the first three terms and the final term follow from those in
\cite{MNS}.  For completeness, we sketch the argument.

Let us first handle the null terms.  By \eqref{2.24} and \eqref{2.25},
the first two terms in \eqref{f.2} are controlled by
\begin{multline}\label{f.3}
C\int_0^t\int_\ext \sum_{\substack{|\alpha|+\nu\le 221\\\nu\le 2}}
|L^\nu Z^\alpha u|\sum_{\substack{|\alpha|+\nu\le 221\\\nu\le 2}}
|L^\nu Z^\alpha u'| \sum_{\substack{|\alpha|+\nu\le 220\\\nu\le 1}}
|L^\nu Z^\alpha (w-v)'|\:\frac{dy\:ds}{|y|}
\\+C\int_0^t\int_\ext \sum_{K=1}^D \frac{\langle
c_Ks-r\rangle}{\langle s+r\rangle} \sum_{\substack{|\alpha|+\nu\le
220\\\nu\le 1}} |L^\nu Z^\alpha
\partial (w-v)^K|\Bigl(\sum_{\substack{|\alpha|+\nu\le 221\\\nu\le 1}} |L^\nu
Z^\alpha u'|\Bigr)^2 \:dy\:ds.
\end{multline}

To handle the first term of \eqref{f.3}, notice
that by \eqref{5.1}, \eqref{5.5}, and \eqref{5.14}, we have
$$\sum_{\substack{|\alpha|+\nu\le 221\\\nu\le 2}} |L^\nu Z^\alpha
u(s,y)|\le C\varepsilon\langle
s+|y|\rangle^{-1+b_3\varepsilon}\log(2+s),$$
which means that the first term of \eqref{f.3} has a contribution to
\eqref{f.2} which is dominated by
\begin{multline*}
C\varepsilon \int_0^t \frac{\log (2+s)}{\langle
s\rangle^{1-b_3\varepsilon}}\sum_{\substack{|\alpha|+\nu\le
221\\\nu\le 2}} \|\langle y\rangle^{-1/2} L^\nu Z^\alpha u'(s,\cd)\|_2
\\\times\sum_{\substack{|\alpha|+\nu\le 220\\\nu\le 1}} \|\langle
y\rangle^{-1/2} L^\nu Z^\alpha (w-v)'(s,\cd)\|_2\:ds
\end{multline*}
by the Schwarz inequality.  Thus, if we subsequently apply the Schwarz
inequality, \eqref{5.1}, \eqref{5.5}, and \eqref{5.20}, we see that
this contribution is $O(\varepsilon^3)$ for $\varepsilon>0$ small.

We now want to show that the second term of \eqref{f.3} satisfies a
similar bound.  If we apply \eqref{6.1}, we see that the second term
of \eqref{f.3} is controlled by
$$C\varepsilon \int_0^t (1+s)^{-(1/2)+\tilde{b}\varepsilon} \int_\ext
\frac{1}{r^{1/2}\langle s+r\rangle^{1/2}}
\sum_{\substack{|\alpha|+\nu\le 221\\\nu\le 1}} |L^\nu Z^\alpha
u'|^2\:dy\:ds.$$
For $\varepsilon$ sufficiently small, it follows similarly that this
term is $O(\varepsilon^3)$ by the $L^2_tL^2_x$ estimates of
\eqref{5.20}.

For the multi-speed terms (i.e. the third term in \eqref{f.2}), let us
for simplicity assume that $I\neq K$, $I=J$.  A symmetric argument
will yield the same bound for the remaining cases.  If we set $\delta
< |c_I-c_K|/3+$, it follows that $\{|y|\in
[(c_I-\delta)s,(c_I+\delta)s]\} \cap
\{[(c_K-\delta)s,(c_K+\delta)s]\}=\emptyset$.  Thus, it will suffice
to show the bound when the spatial integral is taken over the
complements of each of these sets separately.  We will show the bound
over $\{|y|\not\in [(c_K-\delta)s,(c_K+\delta)s]\}$.  A symmetric
argument will yield the bound over the other set.

If we apply \eqref{6.1}, we see that over $\{|y|\not\in
[(c_K-\delta)s,(c_K+\delta) s]\}$ the third term in \eqref{f.2} is
bounded by
$$C\varepsilon \int_0^t \frac{1}{\langle
s\rangle^{(1/2)-\tilde{b}\varepsilon}} \int_{\{|y|\not\in
[(c_K-\delta)s,(c_K+\delta)s]\}} \frac{1}{r}
\sum_{\substack{|\alpha|+\nu\le 221\\\nu\le 1}} |L^\nu Z^\alpha
\partial u^I|^2\:dy\:ds.$$
Arguing as above, it is easy to see that these multiple speed
quadratic terms are also $O(\varepsilon^3)$ by \eqref{5.20}.

Next, we need to show that the last term in \eqref{f.2} enjoys an
$O(\varepsilon^4)$ contribution.  This is clear, however, since this
term is bounded by
$$\int_0^t \sum_{\substack{|\alpha|+\nu\le 221\\\nu\le 1}} \|L^\nu
\partial^\alpha (w-v)'(s,\cd)\|^2_\infty\:ds,$$
and an application of \eqref{c.1} yields the desired bound.

In order to finish the proof of \eqref{f.1}, and hence that of part
$(i.)$, it remains to bound the cubic terms in \eqref{f.2}.  By
\eqref{5.1} and 
\eqref{5.14}, it follows that this fourth term of \eqref{f.2} is
controlled by
$$C\varepsilon^3 \int_0^t \int_\ext \frac{(\log(2+s))^3
(1+s)^{3b_2\varepsilon}}{(1+s+|y|)^3} \sum_{\substack{|\alpha|+\nu\le
220\\\nu\le 1}} |L^\nu Z^\alpha \partial(w-v)|\:dy\:ds.$$
By the Schwarz inequality, \eqref{5.5}, and \eqref{5.17}, this last
term is also $O(\varepsilon^4)$ if $\varepsilon$ is sufficiently
small.  Thus, we have shown \eqref{f.1} and have finished the proof of $(i.)$.

\newsection{Proof of $(ii.)$}
We now begin part $(ii.)$ of the continuity argument.  In particular, we
need to show that \eqref{5.16}, \eqref{5.18}, \eqref{5.19}, and \eqref{5.20} follow
from \eqref{5.12}-\eqref{5.14b} and \eqref{5.17}.  This will complete
the proof of Theorem \ref{theorem1.1}.

\subsection{Proof of \eqref{5.16}:}   In this subsection, we briefly
prove that \eqref{5.16} holds.  Indeed, if $(t,x)\not\in\Lambda_I$, it
follows that $\langle c_It-|x|\rangle\gtrsim \langle t+|x|\rangle$.
Thus, away from the associated light cone, we have that
$1+\log\frac{1+t+|x|}{1+|c_It-|x|}$ is $O(1)$.  Hence, by
\eqref{5.12}, we have
$$(1+t+|x|)\sum_{|\alpha|\le 200} |Z^\alpha \partial w^I(t,x)|\le
C\varepsilon$$
when $(t,x)\not\in\Lambda_I$.  For $(t,|x|)\in \Lambda_I$, it follows
that $|x|\approx t$.  Thus, by \eqref{2.26}, we have
$$(1+t+|x|)\sum_{|\alpha|\le 200} |Z^\alpha \partial w^I(t,x)|\le
C\sum_{|\alpha|\le 202} \|Z^\alpha w'(t,\cd)\|_2$$
provided $(t,x)\in\Lambda_I$.  Since the right side of this equation
is $O(\varepsilon)$ by \eqref{5.17}, we have established \eqref{5.16}
as desired.

\subsection{Proof of \eqref{5.18}:}
The next step is to show that we have the higher order, lossy energy
estimates and mixed norm estimates when the scaling vector field does
not occur.  Here, we modify the arguments of \cite{MS} to allow the
general higher order terms in the nonlinearity.

In the notation of \S 2, we have $\Box_\gamma u = B(du)+P(u,du)$ where
$B(du)+P(u,du)$ is the semilinear part of the nonlinearity and
$$\gamma^{IJ,jk}=\gamma^{IJ,jk}(u,du)=-\sum_{\substack{1\le K\le
D\\0\le l\le 3}} B^{IJ,jk}_{K,l}\partial_l u^K - C^{IJ,jk}(u,du).$$
Also, note that by \eqref{5.12} and \eqref{5.16}
\begin{equation}\label{B.1}
\sum_{|\alpha|\le 1} \|\partial^\alpha \gamma (s,\cd)\|_\infty\le
\frac{C\varepsilon}{1+s}.
\end{equation}

We begin by proving \eqref{5.18}.  To do so, we first estimate the
energy of $\partial^j_t u$ for $j\le M\le 300$.  
Notice that by \eqref{4.5} and \eqref{B.1}, we have
\begin{equation}\label{B.2}
\partial_t E_M^{1/2}(u)(t)\le C\sum_{j\le M} \|\Box_{\gamma}
\partial_t^j u(t,\cd)\|_2 + \frac{C\varepsilon}{1+t} E_M^{1/2}(u)(t).
\end{equation}
Note also that for $M=1,2,\dots$
\begin{multline*}
\sum_{j\le M} |\Box_\gamma \partial_t^j u|\le C\Bigl(\sum_{j\le
M}|\partial_t^j u'|+\sum_{j\le M-1} |\partial^j_t\partial^2 u|\Bigr)
\sum_{|\alpha|\le 200} |\partial^\alpha u'|
\\+C\Bigl(|u|+\sum_{j\le M} |\partial_t^j u'| +\sum_{j\le M-1}
|\partial^j_t \partial^2 u|\Bigr)\Bigl(\sum_{|\alpha|\le 200}
|\partial^\alpha u|\Bigr)^2.
\end{multline*}
Using \eqref{5.1}, \eqref{5.12}, and \eqref{5.16}, it follows that
this is bounded by
$$\frac{C\varepsilon}{1+t}\Bigl(\sum_{j\le M}|\partial^j_t
u'|+\sum_{j\le M-1} |\partial^j_t \partial^2 u|\Bigr) +
\frac{C\varepsilon}{(1+t+|x|)^{3-}}.$$
If we use elliptic regularity and repeat this argument, we get
\begin{align*}
\sum_{j\le M-1} \|\partial_t^j \partial^2 u(t,\cd)\|_2 &\le
C\sum_{j\le M} \|\partial_t^j u'(t,\cd)\|_2 + C\sum_{j\le M-1}
\|\partial_t^j \Box u(t,\cd)\|_2\\
&\le C\sum_{j\le M} \|\partial_t^j u'(t,\cd)\|_2 +
\frac{C\varepsilon}{1+t} \sum_{j\le M-1} \|\partial_t^j \partial^2
u(t,\cd)\|_2 + \frac{C\varepsilon^3}{(1+t)^{3/2-}}.
\end{align*}
If $\varepsilon$ is small, we can absorb the second term into the left
side of the preceding inequality.  Therefore, if we combine the last
two estimates, we conclude that
$$\sum_{j\le M} \|\Box_\gamma \partial^j_t u(t,\cd)\|_2 \le
\frac{C\varepsilon}{1+t} \sum_{j\le M}\|\partial^j_t u'(t,\cd)\|_2 +
\frac{C\varepsilon^3}{(1+t)^{3/2-}}.$$

If we use this in \eqref{B.2}, we get that for small $\varepsilon>0$
$$\partial_t E^{1/2}_M(u)(t)\le
\frac{C\varepsilon}{1+t}E_M^{1/2}(u)(t) +
\frac{C\varepsilon^3}{(1+t)^{3/2-}}$$
since $\frac{1}{2}E^{1/2}_M(u)(t)\le \sum_{j\le M} \|\partial^j_t
u'(t,\cd)\|_2 \le 2 E^{1/2}_M(u)(t)$ when $\varepsilon$ is small.
Since \eqref{1.11} implies that $E^{1/2}_{300}(u)(0)\le C\varepsilon$, 
Gronwall's inequality
yields
\begin{equation}\label{B.3}
\sum_{j\le 300} \|\partial_t^j u'(t,\cd)\|_2 \le
2E_{300}^{1/2}(u)(t)\le C\varepsilon (1+t)^{\tidle{C}\varepsilon}
\end{equation}
for some constant $\tilde{C}>0$.
By elliptic regularity, this leads to the bound \eqref{5.18} if
$\varepsilon>0$ is sufficiently small.

\subsection{Proof of the base case, $N=0$, of \eqref{5.19} and
\eqref{5.20}:}

We begin by showing
\begin{equation}\label{C.1}
(\log(2+t))^{-1/2}\sum_{|\alpha|\le 298} \|\langle x\rangle^{-1/2}
\partial^\alpha u'\|_{L^2(S_t)}\le C (1+t)^{\tilde{C}\varepsilon}
\end{equation}
where $\tilde{C}$ is the constant appearing in \eqref{5.18}.  
By \eqref{4.11}, \eqref{5.1}, and \eqref{5.5}, we have
\begin{equation}\label{C.2}
\begin{split}
(\log(&2+t))^{-1/2}\sum_{|\alpha|\le 298} \|\langle x\rangle^{-1/2}
\partial^\alpha u'\|_{L^2(S_t)}\\
&\le C\varepsilon (\log(2+t))^{1/2} +
(\log(2+t))^{-1/2}\sum_{|\alpha|\le 298} \|\langle
x\rangle^{-1/2}\partial^\alpha (w-v)'\|_{L^2(S_t)}\\
&\le C\varepsilon(\log(2+t))^{1/2} + C\sum_{|\alpha|\le 299}\int_0^t
\|\partial^\alpha \Box (w-v)(s,\cd)\|_2\:ds \\
&\quad\quad\quad\quad\quad\quad\quad\quad\quad\quad\qquad\qquad\qquad\qquad +C \sum_{|\alpha|\le 298}
\|\partial^\alpha \Box (w-v)\|_{L^2(S_t)}.
\end{split}
\end{equation}
Since $\partial^\alpha \Box(w-v)=\partial^\alpha (1-\beta)(1-\eta)\Box
u$, the right side is
$$\le C\varepsilon(\log(2+t))^{1/2} + C\sum_{|\alpha|\le 299}\int_0^t
\|\partial^\alpha \Box u(s,\cd)\|_2\:ds + C\sum_{|\alpha|\le 298}
\|\partial^\alpha \Box u\|_{L^2(S_t)}.$$

If we apply \eqref{5.1}, \eqref{5.12}, and \eqref{5.16} as in the proof of
\eqref{5.18}, it is easy to see that
$$\sum_{|\alpha|\le 299} \|\partial^\alpha \Box u(s,\cd)\|_2 \le
\frac{C\varepsilon}{1+s} \sum_{|\alpha|\le 300} \|\partial^\alpha
u'(s,\cd)\|_2 + \frac{C\varepsilon^3}{(1+s)^{3/2-}}.$$
If we plug this into the previous equation and apply \eqref{5.18},
\eqref{C.1} follows immediately.

We next wish to show
\begin{equation}\label{C.3}
\sum_{|\alpha|\le 297} \|Z^\alpha u'(t,\cd)\|_2 \le C\varepsilon (1+t)^{a_0\varepsilon}
\end{equation}
for some $a_0\ge \tilde{C}$.
In order to show this, we will argue inductively.  That is, for $M\le 297$, we will
assume that
\begin{equation}\label{C.4}
\sum_{|\alpha|\le M-1} \|Z^\alpha u'(t,\cd)\|_2 \le C\varepsilon (1+t)^{C\varepsilon},
\end{equation}
and we will use this to show
\begin{equation}\label{C.5}
\sum_{|\alpha|\le M} \|Z^\alpha u'(t,\cd)\|_2 \le C\varepsilon (1+t)^{C'\varepsilon+\sigma}
\end{equation}
where $\sigma>0$ can be chosen arbitrarily small.  Notice that the
base case follows trivially from \eqref{5.17}.

In order to control the left side of \eqref{C.5}, we use
\eqref{4.10}.  To do so, we must estimate the first term in its right
side.  
We have
\begin{multline}\label{C.6}
\sum_{|\alpha|\le M}\|\Box_\gamma Z^\alpha u(t,\cd)\|_2 
\le C\sum_{|\gamma|\le 200} \|Z^\gamma u'(t,\cd)\|_\infty
\sum_{|\alpha|\le M} \|Z^\alpha u'(t,\cd)\|_2
\\+C\sum_{|\beta|,|\gamma|\le 200, |\alpha|\le M} \|Z^\beta u(t,\cd)
Z^\gamma u(t,\cd) Z^\alpha u'(t,\cd)\|_2
\\+C\sum_{|\beta|,|\gamma|\le 200, |\alpha|\le M} \|Z^\beta u(t,\cd)
Z^\gamma u(t,\cd) Z^\alpha u(t,\cd)\|_2
\end{multline}
By \eqref{5.1}, \eqref{5.12}, and \eqref{5.16}, the first two terms
are controlled by
$$\frac{C\varepsilon}{1+t}\sum_{|\alpha|\le M} \|Z^\alpha
u'(t,\cd)\|_2
\le \frac{C\varepsilon}{1+t} Y_{M,0}^{1/2}(t)$$
where $Y_{M,0}(t)$ is as in \eqref{4.9}.
Since the coefficients of $Z$ are $O(|x|)$, we can apply \eqref{5.1}
and \eqref{5.12} to see that the last term in \eqref{C.6} is dominated
by
$$\frac{C\varepsilon^3}{(1+t)^{3/2-}}
+\frac{C\varepsilon^2(1+\log(1+t))^2}{1+t} \sum_{|\alpha|\le M-1}
\|Z^\alpha u'(t,\cd)\|_2.$$
The first term here corresponds to the case
$|\alpha|=|\beta|=|\gamma|=0$ in the last term of \eqref{C.6}.

Plugging these bounds for \eqref{C.6} into \eqref{4.10} and applying
\eqref{B.1} yields
\begin{multline}\label{C.7}
\partial_t Y_{M,0}(t)\le \frac{C\varepsilon}{1+t}Y_{M,0}(t)
\\+Y^{1/2}_{M,0}(t)\Bigl(\frac{C\varepsilon^3}{(1+t)^{3/2-}} +
\frac{C\varepsilon^2(1+\log(1+t))^2}{1+t} \sum_{|\alpha|\le M-1}
\|Z^\alpha u'(t,\cd)\|_2\Bigr) \\+ C\sum_{|\alpha|\le M+1}
\|\partial^\alpha u'(t,\cd)\|^2_{L^2(|x|<1)}
\end{multline}
if $\varepsilon$ is sufficiently small.  Therefore, by Gronwall's
inequality, \eqref{1.11}, and the inductive hypothesis \eqref{C.4}, 
we have
\begin{multline*}
Y_{M,0}(t)\le C(1+t)^{C\varepsilon} \Bigl[\varepsilon^2 +
\Bigl(\sup_{0\le s\le t} Y^{1/2}_{M,0}(s)\Bigr) \varepsilon^2
(1+t)^{C\varepsilon+\sigma} \\+ \sum_{|\alpha|\le M+1}
\|\langle x\rangle^{-1/2} \partial^\alpha u'\|^2_{L^2(S_t)}\Bigr].
\end{multline*}
If we apply \eqref{C.1} to the last term, we see that \eqref{C.5}
follows.  By induction, this yields \eqref{C.3}.

Using \eqref{4.12}, this in turn implies
\begin{equation}\label{C.8}
(\log(2+t))^{-1/2}\sum_{|\alpha|\le 295} 
\|\langle x\rangle^{-1/2} Z^\alpha u'\|_{L^2(S_t)} \le C\varepsilon(1+t)^{a_0\varepsilon}
\end{equation}
which completes the proof of $N=0$ cases of \eqref{5.19} and \eqref{5.20}.

\subsection{Proof of \eqref{5.19} and \eqref{5.20} for $N>0$:}
In order to complete the proof of Theorem \ref{theorem1.1}, we must
show that \eqref{5.19} and \eqref{5.20} hold for $N=1,2,3$.  To do so,
we argue inductively in $N$.  We fix an $N$ and 
assume that \eqref{5.20} holds with $N$ replaced by $N-1$.  It then
remains to show \eqref{5.19} and \eqref{5.20} for that $N$.

We begin with the task of showing \eqref{5.19}.  The first step will
be to show that
\begin{equation}\label{D.1}
\sum_{\substack{|\alpha|+\mu\le 205\\\mu\le N}} \|L^\mu
\partial^\alpha u'(t,\cd)\|_2\le C\varepsilon (1+t)^{A\varepsilon+\sigma}.
\end{equation}
For this, we shall want to use \eqref{4.8}.  We must first establish
an appropriate version of \eqref{4.7} for $N_0+\nu_0\le 205$,
$\nu_0\le N$.  For this, we note that for $M\le 205$,
\begin{multline*}
\sum_{\substack{j+\mu\le M\\\mu\le N}}\Bigl(|\tilde{L}^\mu
\partial_t^j \Box_\gamma u|+|[\tilde{L}^\mu
\partial^j_t,\Box-\Box_\gamma]  u|\Bigr)
\\\le C\Bigl(\sum_{j\le M-N} |\tilde{L}^N\partial_t^j \partial u| +
\sum_{j\le M-N-1} |\tilde{L}^N\partial^j_t \partial^2
u|\Bigr)\Bigl(\sum_{|\alpha|\le 200} |\partial^\alpha u'| +
\Bigl(\sum_{|\alpha|\le 200} |\partial^\alpha u|\Bigr)^2\Bigr)
\\+ C\sum_{\substack{|\alpha|+\mu\le M-200\\\mu\le N}} |L^\mu
\partial^\alpha u'| \sum_{\substack{|\alpha|+\mu\le M\\\mu\le N-1}}
|L^\mu \partial^\alpha u'|
\\+C\sum_{\substack{|\alpha|+\mu\le M\\\mu\le N-1}} |L^\mu
\partial^\alpha u'| \sum_{\substack{|\alpha|+\mu\le
\max(M/2,M-200)\\\mu\le N-1}} |L^\mu \partial^\alpha u'|
\\+ C\Bigl(\sum_{\substack{|\alpha|+\mu\le M+1\\\mu\le N-1}} |L^\mu
\partial^\alpha u|\Bigr)^2 \sum_{\substack{|\alpha|+\mu \le
180\\\mu\le N}} |L^\mu \partial^\alpha u|.
\end{multline*}

\noindent Using elliptic regularity, \eqref{2.26}, \eqref{5.1},
\eqref{5.12}, and \eqref{5.17}, we conclude that
\begin{multline*}
\sum_{\substack{j+\mu\le M\\\mu\le N}} \Bigl(\|\tilde{L}^\mu
\partial^j_t \Box_\gamma u(t,\cd)\|_2 + \|[\tilde{L}^\mu\partial_t^j,
\Box-\Box_\gamma] u(t,\cd)\|_2\Bigr)\le \frac{C\varepsilon}{1+t}
\sum_{\substack{j+\mu\le M\\\mu\le N}} \|\tilde{L}^\mu \partial^j_t
u'(t,\cd)\|_2
\\+C\sum_{\substack{|\alpha|+\mu\le M-200\\\mu\le N}} \|\langle
x\rangle^{-1/2} L^\mu \partial^\alpha u'(t,\cd)\|_2
\sum_{\substack{|\alpha|+\mu\le 207\\\mu\le N-1}} \|\langle
x\rangle^{-1/2} L^\mu Z^\alpha u'(t,\cd)\|_2
\\+C\sum_{\substack{|\alpha|+\mu\le \max(M,2+M/2)\\\mu\le N-1}}
\|\langle x\rangle^{-1/2} L^\mu Z^\alpha u'(t,\cd)\|_2^2
\\+C\Bigl\|
\Bigl(\sum_{\substack{|\alpha|+\mu\le M+1\\\mu\le N-1}} |L^\mu
\partial^\alpha u(t,\cd)|\Bigr)^2 \sum_{\substack{|\alpha|+\mu \le
180\\\mu\le N}} |L^\mu \partial^\alpha u(t,\cd)|\Bigr\|_2.
\end{multline*}
Based on this, \eqref{4.7} holds with $\delta=C\varepsilon$ and
\begin{multline*}
H_{N,M-N}(t)=C\sum_{\substack{|\alpha|+\mu\le M-200\\\mu\le N}}
\|\langle x\rangle^{-1/2} L^\mu \partial^\alpha u'(t,\cd)\|_2^2 \\+
C\sum_{\substack{|\alpha|+\mu\le 207\\\mu\le N-1}} \|\langle
x\rangle^{-1/2} L^\mu Z^\alpha u'(t,\cd)\|_2^2
\\+C\Bigl\|\Bigl(\sum_{\substack{|\alpha|+\mu\le M+1\\\mu\le N-1}}
|L^\mu \partial^\alpha u(t,\cd)|\Bigr)^2
\sum_{\substack{|\alpha|+\mu\le 180\\\mu\le N}} |L^\mu \partial^\alpha
u(t,\cd)| \Bigr\|_2
\end{multline*}
Since $M+1\le 206$, \eqref{5.1}, \eqref{5.14}, and \eqref{5.15} imply
that this last term is controlled by
$$C\varepsilon^3 (1+\log(1+t))^3 (1+t)^{2b_N+\tilde{c}_N}
\Bigl\|\frac{1}{(1+t+|x|)^2} \frac{1}{|x|}\Bigr\|_2 \le
\frac{C\varepsilon^3}{(1+t)^{1+}}$$
if $\varepsilon$ is sufficiently small.

Since the conditions on the data give $\int e_0(\tilde{L}^\nu
\partial^j_t u)(0,x)\:dx \le C\varepsilon^2$ if $j+\nu\le 300$, it
follows from \eqref{4.8} and the inductive hypothesis (\eqref{5.20}
with $N$ replaced by $N-1$) that for $M\le 205$
\begin{multline}\label{D.2}
\sum_{\substack{|\alpha|+\mu\le M\\\mu\le N}} \|L^\mu \partial^\alpha
u'(t,\cd)\|_2 \le C\varepsilon (1+t)^{C\varepsilon+\sigma} 
\\+C(1+t)^{C\varepsilon}\sum_{\substack{|\alpha|+\mu\le M-200\\\mu\le
N}} \|\langle x\rangle^{-1/2} L^\mu \partial^\alpha u'\|^2_{L^2(S_t)}
\\+C(1+t)^{C\varepsilon} \int_0^t \sum_{\substack{|\alpha|+\mu\le
M\\\mu\le N-1}} \|L^\mu \partial^\alpha u'(s,\cd)\|_{L^2(|x|<1)}\:ds
\end{multline}
for some constant $\sigma>2a_{N-1}\varepsilon$.

If we apply \eqref{4.14}, \eqref{5.1}, and \eqref{5.5}, we get that
the last integral is dominated by $C\varepsilon \log(2+t)$ plus
\begin{multline*}
\int_0^t \sum_{\substack{|\alpha|+\mu\le M\\\mu\le N-1}} \|L^\mu
\partial^\alpha (w-v)'(s,\cd)\|_{L^2(|x|<1)}\:ds 
\\\le C\sum_{\substack{|\alpha|+\mu\le M+1\\\mu\le N-1}} \int_0^t
\Bigl(\int_0^s \|L^\mu \partial^\alpha \Box
(w-v)(\tau,\cd)\|_{L^2(||x|-(s-\tau)|<10)} \:d\tau\Bigr)\:ds.
\end{multline*}
Since $\Box (w-v)=(1-\beta)(1-\eta)\Box u$, we conclude that this last
term is bounded by
\begin{equation}\label{D.3}
\sum_{\substack{|\alpha|+\mu\le M+1\\\mu\le N-1}} \int_0^t
\Bigl(\int_0^s \|L^\mu \partial^\alpha \Box
u(\tau,\cd)\|_{L^2(||x|-(s-\tau)|<10)} \:d\tau\Bigr)\:ds.
\end{equation}

As in \cite{MS}, when $\Box u$ is replaced by 
the quadratic terms $B(du)+Q(du,d^2u)$ in \eqref{D.3}, we see
from an application of \eqref{2.26} that the integrand is bounded by
$$\sum_{\substack{|\alpha|+\mu\le 209\\\mu\le N-1}} \|\langle
x\rangle^{-1/2} L^\mu Z^\alpha
u'(\tau,\cd)\|^2_{L^2(||x|-(s-\tau)|<20)}.$$
Since the sets $\{(\tau,x)\,:\,||x|-(j-\tau)|<20\}$, $j=0,1,2,\dots$
have finite overlap, we conclude that, in this case, the last integral
in \eqref{D.2} is bounded by
$$C\varepsilon \log(2+t)+C\sum_{\substack{|\alpha|+\mu\le 209\\\mu\le
N-1}} \|\langle x\rangle^{-1/2} L^\mu Z^\alpha u'\|^2_{L^2(S_t)} \le
C\varepsilon (1+t)^{2a_{N-1}\varepsilon}.$$
The last inequality follows from the inductive hypothesis.

When $\Box u$ is replaced by the higher order terms
$P(u,du)+R(u,du,d^2u)$, we see that \eqref{D.3} is bounded by
\begin{equation}\label{D.4}
C\int_0^t \Bigl(\int_0^s \Bigl\|\Bigl(\sum_{\substack{|\alpha|+\mu\le
208\\\mu\le N-2}} |L^\mu \partial^\alpha u(\tau,\cd)|\Bigr)^2
\sum_{\substack{|\alpha|+ \mu\le 208\\\mu\le N-1}} |L^\mu
\partial^\alpha
u(\tau,\cd)|\Bigr\|_{L^2(||x|-(s-\tau)|<10)}\:d\tau\Bigr)\:ds 
\end{equation}
By \eqref{5.1}, \eqref{5.14}, and \eqref{5.14b}, we have that
\begin{multline*}
\Bigl(\sum_{\substack{|\alpha|+\mu\le
208\\\mu\le N-2}} |L^\mu \partial^\alpha u(\tau,\cd)|\Bigr)^2
\sum_{\substack{|\alpha|+ \mu\le 208\\\mu\le N-1}} |L^\mu
\partial^\alpha
u(\tau,\cd)| \\\le C \varepsilon^3
(1+t)^{2b_{N-1}\varepsilon+c'_{N-1}\varepsilon+} (1+t)^{-1} (1+|x|)^{-2}.
\end{multline*}
Since the norm is taken over $|x|\approx (t-s)$, it follows that
\eqref{D.4} is bounded by $C\varepsilon^3
(1+t)^{2b_{N-1}\varepsilon+c'_{N-1}\varepsilon+}$.  Since we may take
$\tilde{a}_N > 2b_{N-1}+c'_{N-1}$, this will be sufficient for $N\ge
2$.  
When $N=1$, there are no occurences of $L$ in \eqref{D.4}, and
appropriate bounds follow simply from \eqref{5.1}, \eqref{5.12}, 
and \eqref{5.14b}.

Plugging this into \eqref{D.2}, we see that
\begin{multline}\label{D.5}
\sum_{\substack{|\alpha|+\mu\le M\\\mu\le N}} \|L^\mu \partial^\alpha
u'(t,\cd)\|_2 \le C\varepsilon (1+t)^{C\varepsilon+\sigma}
\\+C(1+t)^{C\varepsilon}\sum_{\substack{|\alpha|+\mu\le M-200\\\mu\le
N}} \|\langle x\rangle^{-1/2} L^\mu \partial^\alpha u'\|^2_{L^2(S_t)}
\end{multline}
which yields \eqref{D.1} for $M\le 200$.

For $M>200$, \eqref{D.1} will follow from a simple induction argument
using the following lemma.

\begin{lemma}\label{lemma7.1}
Under the above assumptions, if $M\le 205$, $1\le N\le 4$, and
\begin{multline}\label{D.6}
\sum_{\substack{|\alpha|+\nu\le M\\\nu\le N}} \|L^\nu \partial^\alpha
u'(t,\cd)\|_2 + \sum_{\substack{|\alpha|+\nu\le M-3\\\nu\le N}}
\|\langle x\rangle^{-1/2} L^\nu \partial^\alpha u'\|_{L^2(S_t)}
\\+\sum_{\substack{|\alpha|+\nu\le M-4\\\nu\le N}} \|L^\nu Z^\alpha
u'(t,\cd)\|_2 + \sum_{\substack{|\alpha|+\nu\le M-6\\\nu\le N}}
\|\langle x\rangle^{-1/2} L^\nu Z^\alpha u'\|_{L^2(S_t)} \le C\varepsilon
(1+t)^{C\varepsilon+\sigma} 
\end{multline}
with $\sigma>0$, then there is a constant $C'$ so that
\begin{multline}\label{D.7}
\sum_{\substack{|\alpha|+\nu\le M-2\\\nu\le N}} \|\langle
x\rangle^{-1/2} L^\nu \partial^\alpha u'\|_{L^2(S_t)} +
\sum_{\substack{|\alpha|+\nu\le M-3\\\nu\le N}} \|L^\nu Z^\alpha
u'(t,\cd)\|_2
\\+\sum_{\substack{|\alpha|+\nu\le M-5\\\nu\le N}} \|\langle
x\rangle^{-1/2} L^\nu Z^\alpha u'\|_{L^2(S_t)}\le C'\varepsilon
(1+t)^{C'\varepsilon +C'\sigma}.
\end{multline}
\end{lemma}

\noindent{\em Proof of Lemma \ref{lemma7.1}:}  Let us start with the
first term on the left side of \eqref{D.7}.  Using \eqref{4.11},
\eqref{5.1}, and \eqref{5.5} as in \eqref{C.2}, we see that 
$$(\log(2+t))^{-1/2} \sum_{\substack{|\alpha|+\nu\le M-2\\\nu\le N}}
\|\langle x\rangle^{-1/2} L^\nu \partial^\alpha u'\|_{L^2(S_t)}$$
 is
controlled by $C\varepsilon (\log(2+t))^{1/2}$ plus 
\begin{equation}\label{D.8}
C\sum_{\substack{|\alpha|+\nu\le M-1\\\nu\le N}} \int_0^t \|L^\nu
\partial^\alpha \Box u(s,\cd)\|_2\:ds +
C\sum_{\substack{|\alpha|+\nu\le M-2\\\nu\le N}} \|L^\nu
\partial^\alpha \Box u\|_{L^2(S_t)}.
\end{equation}

When $\Box u$ is replaced by the quadratic terms $B(du)+Q(du,d^2u)$,
the first term in \eqref{D.8} is controlled by
\begin{multline*}
C\int_0^t \Bigl\|\sum_{|\alpha|\le 200} |\partial^\alpha u'(s,\cd)|
\sum_{\substack{|\alpha|+\nu\le M\\\nu\le N}} |L^\nu \partial^\alpha
u'(s,\cd)| \Bigr\|_2\:ds
\\+C\int_0^t\Bigl\|\sum_{|\alpha|\le M} |\partial^\alpha
u'(s,\cd)|\sum_{\substack{|\alpha|+\nu\le M-200\\\nu\le N}} |L^\nu
\partial^\alpha u'(s,\cd)|\Bigr\|_2\:ds
\\+C\int_0^t \Bigl\|\Bigl(\sum_{\substack{|\alpha|+\nu\le M\\\nu\le
N-1}} |L^\nu \partial^\alpha u'(s,\cd)|\Bigr)^2\Bigr\|_2\:ds.
\end{multline*}

Notice that by \eqref{5.16} and \eqref{D.6}, the desired bound holds
for the first term.  By \eqref{2.26}, the last term is bounded by
$$C\sum_{\substack{|\alpha|+\nu\le M+2\\\nu\le N-1}} \|\langle
x\rangle^{-1/2}L^\nu Z^\alpha u'\|^2_{L^2(S_t)},$$
and the appropriate bound follows from the \eqref{5.20} with $N$
replaced by $N-1$.  This is sufficient to show that the result holds
for this case when $M\le 200$.  When $M\ge 201$, we must also handle
the second term above.  By \eqref{2.26}, this is controlled by
$$C\sum_{|\alpha|\le M} \|\langle x\rangle^{-1/2} \partial^\alpha
u'\|^2_{L^2(S_t)} + C\sum_{\substack{|\alpha|+\nu\le M-198\\\nu\le N}}
\|\langle x\rangle^{-1/2} L^\nu Z^\alpha u'\|^2_{L^2(S_t)},$$
and the bounds follow from \eqref{C.1} and \eqref{D.6}.  When $\Box u$
is quadratic, the desired estimates for the last term in \eqref{D.8}
follow from very similar arguments.

It remains to bound \eqref{D.8} when $\Box u$ is replaced by the
higher order terms $P(u,du)+R(u,du,d^2u)$.  Since $M+1\le 206$, by \eqref{5.1},
\eqref{5.14} and \eqref{5.14b}, we have
$$\Bigl\|\Bigl(\sum_{\substack{|\alpha|+\nu\le M+1\\\nu\le N-1}}
|L^\nu \partial^\alpha u(s,\cd)|\Bigr)^2
\sum_{\substack{|\alpha|+\nu\le M+1\\\nu\le N}} |L^\nu \partial^\alpha
u(s,\cd)|\Bigr\|_2 \le C\varepsilon^2 (1+t)^{-1-}$$
for $\varepsilon>0$ sufficiently small.
Upon integration, it is easy to see that these terms satisfy the
desired bounds, and this finishes the proof that the first term on the
left side of \eqref{D.7} is dominated by its right side.

To control the second term in \eqref{D.7}, we will use \eqref{4.10}.
This means that we must estimate the first term in its right side,
which satisfies
\begin{multline*}
\sum_{\substack{|\alpha|+\nu\le M-3\\\nu\le N}} \|\Box_\gamma L^\nu
Z^\alpha u(s,\cd)\|_2 \le C\Bigl\|\sum_{|\alpha|\le 200} |Z^\alpha
u'(s,\cd)|
\sum_{\substack{|\alpha|+\nu\le M-3\\\nu\le N}} |L^\nu Z^\alpha
u'(s,\cd)|\Bigr\|_2
\\+C \Bigl\|\sum_{|\alpha|\le M-3}|Z^\alpha u'(s,\cd)|
\sum_{\substack{|\alpha|+\nu \le M-203\\\nu\le N}} |L^\nu Z^\alpha
u'(s,\cd)|\Bigr\|_2
\\+C\Bigl\|\Bigl(\sum_{\substack{|\alpha|+\nu\le M-3\\\nu\le N-1}}
|L^\nu Z^\alpha u'(s,\cd)|\Bigr)^2\Bigr\|_2
\\+C\Bigl\|\Bigl(\sum_{\substack{|\alpha|+\nu\le M-2\\\nu\le N-1}}
|L^\nu Z^\alpha u(s,\cd)|\bigr)^2 \sum_{\substack{|\alpha|+\nu\le
M-2\\\nu\le N}} |L^\nu Z^\alpha u(s,\cd)|\Bigr\|_2.
\end{multline*}
Applying \eqref{5.1} and \eqref{5.16} to the first term, \eqref{2.26}
to the second and third terms, and \eqref{5.1}, \eqref{5.14}, and
\eqref{5.14b} to the cubic term, we have that this is dominated by
\begin{multline}\label{D.9}
\frac{C\varepsilon}{1+s}Y^{1/2}_{M-3-N,N}(s) +
C\sum_{\substack{|\alpha|+\nu\le M-1\\\nu\le N-1}} \|\langle
x\rangle^{-1/2} L^\nu Z^\alpha u'(s,\cd)\|_2^2
\\+C\sum_{\substack{|\alpha|+\nu\le M-203\\\nu\le N}} \|\langle
x\rangle^{-1/2} L^\nu Z^\alpha u'(s,\cd)\|_2^2 + C\varepsilon^3 (1+s)^{-1-}
\end{multline}
with $Y_{M-3-N,N}(t)$ as in \eqref{4.9}.

Plugging this into \eqref{4.10}, applying the inductive hypothesis
(\eqref{5.20} with $N$ replaced by $N-1$), using Gronwall's
inequality, and arguing as in the proof of \eqref{C.3}, we see that
the bound for the second term in \eqref{D.7}  
follows for $M\le 203$.  When $M\ge 203$, we must also deal
with the third term in \eqref{D.9}, but this is done trivially by
applying \eqref{D.6}.

Using \eqref{4.12} and the arguments that procede, this in turn
implies that the third term in \eqref{D.7} is bounded by the right
side, which completes the proof.\qed

From \eqref{D.5} and the lemma, one gets \eqref{D.1} and
\eqref{5.19}.  

In order to complete the proof of Theorem \ref{theorem1.1}, one must
show that \eqref{5.20} follows from \eqref{5.12}-\eqref{5.19} and
\eqref{5.20} with $N$ replaced by $N-1$.  The first step is to show
that
\begin{equation}\label{D.10}
\sum_{\substack{|\alpha|+\nu\le M\\\nu\le N}} \|L^\nu \partial^\alpha
u'(t,\cd)\|_2 \le C\varepsilon (1+t)^{A_N\varepsilon}
\end{equation}
for some $A_N$.  For this, as in the proof of \eqref{D.1}, we will use
\eqref{4.8} once we are able to establish an appropriate version of
\eqref{4.7} for $N_0+\nu_0\le 300-8N$, $\nu_0\le N$.  Notice that for
$M\le 300-8N$, we have
\begin{multline*}
\sum_{\substack{j+\mu\le M\\\mu\le N}} \Bigl(|\tilde{L}^\mu
\partial_t^j \Box_\gamma u|+|[\tilde{L}^\mu \partial_t^j,
\Box-\Box_\gamma]u|\Bigr) 
\\\le C\Bigl(\sum_{\substack{j+\mu\le M\\\mu\le N}} |\tilde{L}^\mu
\partial_t^j \partial u| + \sum_{\substack{j+\mu\le M-1\\\mu\le N}}
|\tilde{L}^N \partial^j_t \partial^2 u|\Bigr)\Bigl(\sum_{|\alpha|\le
200} |\partial^\alpha u'|+\Bigl(\sum_{\substack{|\alpha|+\mu\le 190\\\mu\le
N-1}}|L^\mu \partial^\alpha u|\Bigr)^2\Bigr)
\\+C\sum_{\substack{|\alpha|+\mu\le M-200\\\mu\le N}} |L^\mu
\partial^\alpha u'| \sum_{\substack{|\alpha|+\mu\le M\\\mu\le N-1}}
|L^\mu \partial^\alpha u'|
\\+C\sum_{\substack{|\alpha|+\mu\le M\\\mu\le N-1}} |L^\mu
\partial^\alpha u'| \sum_{\substack{|\alpha|+\mu\le
\max(M/2,M-200)\\\mu\le N-1}}|L^\mu \partial^\alpha u'|
\\+C\Bigl(\sum_{\substack{|\alpha|+\mu\le 255\\\mu\le N-1}} |L^\mu
\partial^\alpha u|\Bigr)^2 \sum_{\substack{|\alpha|+\mu\le 180\\\mu\le
N}} |L^\mu \partial^\alpha u|
\\+C\sum_{\substack{|\alpha|+\mu\le M\\\mu\le N-1}} |L^\mu
\partial^\alpha u'| \sum_{\substack{|\alpha|+\mu\le 190\\\mu\le N-1}}
|L^\mu \partial^\alpha u| \sum_{\substack{|\alpha|+\mu\le 180\\\mu\le
N}} |L^\mu \partial^\alpha u|.
\end{multline*}
By this, \eqref{2.26}, \eqref{5.1},  \eqref{5.14}, and \eqref{5.16},
and elliptic regularity, we get that for $M\le 300-8N$
\begin{multline*}
\sum_{\substack{j+\mu\le M\\\mu\le N}} \Bigl(\|\tilde{L}^\mu
\partial_t^j \Box_\gamma
u(t,\cd)\|_2+\|[\tilde{L}^\mu\partial^j_t,\Box-\Box_\gamma]
u(t,\cd)\|_2\Bigr) \le \frac{C\varepsilon}{1+t}
\sum_{\substack{j+\mu\le M\\\mu\le N}} \|\tilde{L}^\mu \partial^j_t
u'(t,\cd)\|_2 
\\+C\sum_{\substack{|\alpha|+\mu\le M-200\\\mu\le N}} \|\langle
x\rangle^{-1/2} L^\mu \partial^\alpha u'(t,\cd)\|^2_2
+C\sum_{\substack{|\alpha|+\mu\le M+2\\\mu\le N-1}} \|\langle
x\rangle^{-1/2} L^\mu Z^\alpha u'(t,\cd)\|^2_2
\\+C\Bigl\|\Bigl(\sum_{\substack{|\alpha|+\mu\le 255\\\mu\le N-1}} |L^\mu
\partial^\alpha u|\Bigr)^2 \sum_{\substack{|\alpha|+\mu\le 180\\\mu\le
N}} |L^\mu \partial^\alpha u|\Bigr\|_2
\\+C\Bigl\|\sum_{\substack{|\alpha|+\mu\le M\\\mu\le N-1}} |L^\mu
\partial^\alpha u'| \sum_{\substack{|\alpha|+\mu\le 190\\\mu\le N-1}}
|L^\mu \partial^\alpha u| \sum_{\substack{|\alpha|+\mu\le 180\\\mu\le
N}} |L^\mu \partial^\alpha u|\Bigr\|_2.
\end{multline*}
Based on this, \eqref{4.7} holds with
\begin{multline}\label{D.11}
H_{N,M-N}(t)=C\sum_{\substack{|\alpha|+\mu\le M-200\\\mu\le N}} \|\langle
x\rangle^{-1/2} L^\mu \partial^\alpha u'(t,\cd)\|^2_2
\\+C\sum_{\substack{|\alpha|+\mu\le M+2\\\mu\le N-1}} \|\langle
x\rangle^{-1/2} L^\mu Z^\alpha u'(t,\cd)\|^2_2
+C\Bigl\|\Bigl(\sum_{\substack{|\alpha|+\mu\le 255\\\mu\le N-1}} |L^\mu
\partial^\alpha u|\Bigr)^2 \sum_{\substack{|\alpha|+\mu\le 180\\\mu\le
N}} |L^\mu \partial^\alpha u|\Bigr\|_2
\\+C\Bigl\|\sum_{\substack{|\alpha|+\mu\le M\\\mu\le N-1}} |L^\mu
\partial^\alpha u'| \sum_{\substack{|\alpha|+\mu\le 190\\\mu\le N-1}}
|L^\mu \partial^\alpha u| \sum_{\substack{|\alpha|+\mu\le 180\\\mu\le
N}} |L^\mu \partial^\alpha u|\Bigr\|_2.
\end{multline}

Notice that by \eqref{5.1}, \eqref{5.14}, and \eqref{5.15}, 
the third term on the right of \eqref{D.11} is controlled
by $C\varepsilon^3 (1+t)^{-1-}$.  Also notice that by \eqref{5.1},
\eqref{5.13}, and \eqref{5.15}, the last term in \eqref{D.11} is
dominated by
$$\frac{C\varepsilon^2
(1+t)^{\tilde{b}_N\varepsilon+\tilde{c}_N\varepsilon+}}{1+t}
\sum_{\substack{|\alpha|+\mu\le M \\\mu\le N-1}}\|\langle 
x\rangle^{-1} L^\mu \partial^\alpha u'(t,\cd)\|_2.$$  
From this, we see that
\begin{multline}\label{D.12}
\int_0^t H_{N,M-N}(s)\:ds \le C\varepsilon^3 +
C\sum_{\substack{|\alpha|+\mu\le M-200\\\mu\le N}} \|\langle
x\rangle^{-1/2} L^\mu \partial^\alpha u'\|_{L^2(S_t)}^2
\\+C\sum_{\substack{|\alpha|+\mu\le M+2\\\mu\le N-1}} \|\langle
x\rangle^{-1/2} L^\mu Z^\alpha u'\|^2_{L^2(S_t)}
\\+C\varepsilon^2 \int_0^t
(1+s)^{-1+\tilde{b}_N\varepsilon+\tilde{c}_N\varepsilon+}
\sum_{\substack{|\alpha|+\mu\le M\\\mu\le N-1}} \|\langle
x\rangle^{-1} L^\mu \partial^\alpha u'(s,\cd)\|_2\:ds.
\end{multline}
If one applies the Schwarz inequality and
uses \eqref{5.20} (with $N$ replaced by $N-1$), the last term above is
$O(\varepsilon^3)$ for sufficiently small $\varepsilon$.

If we use this in \eqref{4.8} and apply the inductive hypothesis to
handle terms that involve $N-1$ or fewer occurences of $L$, we see
that
\begin{multline}\label{D.13}
\sum_{\substack{|\alpha|+\mu\le M\\\mu\le N}} \|L^\mu \partial^\alpha
u'(t,\cd)\|_2 \le C\varepsilon (1+t)^{C\varepsilon+\sigma}
\\+C(1+t)^{C\varepsilon} \sum_{\substack{|\alpha|+\mu\le
M-200\\\mu\le N}} \|\langle x\rangle^{-1/2} L^\mu \partial^\alpha
u'\|^2_{L^2(S_t)}
\\+C (1+t)^{C\varepsilon} \int_0^t \sum_{\substack{|\alpha|+\mu\le
M\\\mu\le N-1}} \|L^\mu \partial^\alpha u'(s,\cd)\|_{L^2(|x|<1)} \:ds
\end{multline}
since the conditions on the data give $\int e_0(\tilde{L}^v
\partial_t^j u)(0,x)\:ds\le C\varepsilon^2$ if $\nu+j\le 300$.

As before, if we apply \eqref{4.14}, \eqref{5.1}, and \eqref{5.5}, the last
integral is dominated by $C\varepsilon \log(2+t)$ plus
\begin{equation}\label{D.14}
\sum_{\substack{|\alpha|+\mu\le M+1\\\mu\le N-1}} \int_0^t
\Bigl(\int_0^s \|L^\mu \partial^\alpha \Box
u(\tau,\cd)\|_{L^2(||x|-(s-\tau)|<10)} \:d\tau\Bigr)\:ds.
\end{equation}
When $\Box u$ is replaced by $B(du)+Q(du,d^2u)$, as in the proof of
\eqref{D.1}, we can apply \eqref{2.26} and finite overlap of the sets
$\{(\tau,x)\,:\, ||x|-(j-\tau)|<20\}$, $j=0,1,2,\dots$ to see that this
is bounded by
$$C\sum_{\substack{|\alpha|+\mu\le M+3\\\mu\le N-1}} \|\langle
x\rangle^{-1/2} L^\mu Z^\alpha u'\|_{L^2(S_t)}^2 \le C\varepsilon^2
(1+t)^{2a_{N-1}\varepsilon}.$$
The last inequality follows from the inductive hypothesis
\eqref{5.20}.

We must now examine the case when  
$\Box u$ in \eqref{D.14} is replaced by the cubic terms
$P(u,du)+R(u,du,d^2u)$.  Here, we see that \eqref{D.14} is bounded by
\begin{multline}\label{D.15}
C\int_0^t\Bigl(\int_0^s \Bigl\|\Bigl(\sum_{\substack{|\alpha|+\mu\le
190\\\mu\le N-1}} |L^\mu \partial^\alpha u(\tau,\cd)|\Bigr)^2
\\\times\sum_{\substack{|\alpha|+\mu\le
M+2\\\mu\le N-1}} |L^\mu \partial^\alpha
u'(\tau,\cd)|\Bigr\|_{L^2(||x|-(s-\tau)|<10)} \:d\tau\Bigr)\:ds
\\+C\int_0^t\Bigl(\int_0^s \Bigl\|\Bigl(\sum_{\mu\le N-1} |L^\mu
u(\tau,\cd)|\Bigr)^3\Bigr\|_{L^2(||x|-(s-\tau)|<10)} \:d\tau\Bigr)\:ds.
\end{multline}
Since the norm is taken over $|x|\approx (s-\tau)$, we can apply
\eqref{5.1} and \eqref{5.13} to bound the first term by
$$C\varepsilon^2 \int_0^t \Bigl(\int_0^s
\frac{(1+\tau)^{2\tilde{b}_N\varepsilon+}}{(1+\tau) (1+(s-\tau))}
\sum_{\substack{|\alpha|+\mu\le M+2\\\mu\le N-1}} \|L^\mu
\partial^\alpha u'(\tau,\cd)\|_{L^2(||x|-(s-\tau)|<10)}\:d\tau\Bigr)\:ds.$$
By the inductive hypothesis \eqref{5.20}, it follows that this term is
dominated by $C\varepsilon^3 (1+t)^{2\tilde{b}_N\varepsilon +
a_{N-1}\varepsilon +}$.  
Using three applications of \eqref{5.1} and \eqref{5.13} we see that
$$\Bigl(\sum_{\mu\le N-1} |L^\mu u(\tau,x)|\Bigr)^3 \le C\varepsilon^3
(1+\tau)^{3\tilde{b}_N\varepsilon +} (1+\tau)^{-1} (1+|x|)^{-2},$$
and thus it follows that the last term in \eqref{D.15} is controlled
by $C\varepsilon^3 (1+t)^{3\tilde{b}_N\varepsilon +}$.

Plugging these bounds in \eqref{D.13}, it follows that
\begin{multline}\label{D.16}
\sum_{\substack{|\alpha|+\mu\le M\\\mu\le N}} \|L^\mu \partial^\alpha
u'(t,\cd)\|_2 \le C\varepsilon (1+t)^{C\varepsilon+\sigma}
\\+C(1+t)^{C\varepsilon} \sum_{\substack{|\alpha|+\mu\le M-200\\\mu\le
N}} \|\langle x\rangle^{-1/2}L^\mu \partial^\alpha u'\|^2_{L^2(S_t)}
\end{multline}
which yields \eqref{D.10} for $M\le 200$.  
For $M>200$, similar to what we have seen previously, \eqref{D.10}
will follow from a simple induction argument using the following
lemma.

\begin{lemma}\label{lemma7.2}
Under the above assumptions, if $M\le 300-8N$ and
\begin{multline}\label{D.17}
\sum_{\substack{|\alpha|+\nu\le M\\\nu\le N}}\|L^\nu \partial^\alpha
u'(t,\cd)\|_2 + \sum_{\substack{|\alpha|+\nu\le M-3\\\nu\le N}}
\|\langle x\rangle^{-1/2} L^\nu \partial^\alpha u'\|_{L^2(S_t)}
\\+\sum_{\substack{|\alpha|+\nu\le M-4\\\nu\le N}} \|L^\nu Z^\alpha
u'(t,\cd)\|_2 + \sum_{\substack{|\alpha|+\nu\le M-6\\\nu\le N}}
\|\langle x\rangle^{-1/2}L^\nu Z^\alpha u'\|_{L^2(S_t)}\le
C\varepsilon (1+t)^{C\varepsilon+\sigma}
\end{multline}
with $\sigma>0$, then there is a constant $C'$ so that
\begin{multline}\label{D.18}
\sum_{\substack{|\alpha|+\nu\le M-2\\\nu\le N}} \|\langle
x\rangle^{-1/2} L^\nu \partial^\alpha u'\|_{L^2(S_t)} +
\sum_{\substack{|\alpha|+\nu\le M-3\\\nu\le N}} \|L^\nu Z^\alpha
u'(t,\cd)\|_2
\\+\sum_{\substack{|\alpha|+\nu\le M-5\\\nu\le N}} \|\langle
x\rangle^{-1/2} L^\nu Z^\alpha u'\|_{L^2(S_t)} \le C'\varepsilon
(1+t)^{C'\varepsilon +C'\sigma}.
\end{multline}
\end{lemma}

\noindent{\it Proof of Lemma \ref{lemma7.2}:} Here we use arguments
similar to those applied to prove Lemma \ref{lemma7.1}.  We begin by
showing that the first term on the left side of \eqref{D.18} satisfies
the bound.  Using \eqref{4.11}, \eqref{5.1}, and \eqref{5.5} as in
\eqref{C.2}, we see that
\begin{multline}\label{D.19}
(\log(2+t))^{-1/2}\sum_{\substack{|\alpha|+\nu\le M-2\\\nu\le N}}
\|\langle x\rangle^{-1/2}L^\nu\partial^\alpha u'\|_{L^2(S_t)} \le
C\varepsilon (\log(2+t))^{1/2} 
\\+C\sum_{\substack{|\alpha|+\nu\le M-1\\\nu\le N}} \int_0^t \|L^\nu
\partial^\alpha \Box u(s,\cd)\|_2\:ds +
C\sum_{\substack{|\alpha|+\nu\le M-2\\\nu\le N}} \|L^\nu
\partial^\alpha \Box u\|_{L^2(S_t)}.
\end{multline}
Notice that the second term in the right side of \eqref{D.19} is
\begin{multline}\label{D.20}
\le C\int_0^t
\Bigl\|\sum_{\substack{|\alpha|+\nu\le 190\\\nu\le N-1}}
|L^\nu \partial^\alpha u'(s,\cd)|\sum_{\substack{|\alpha|+\nu\le M\\\nu\le
N}} 
|L^\nu \partial^\alpha u'(s,\cd)|\Bigr\|_2\:ds
\\+C\int_0^t \Bigl\|\sum_{\substack{|\alpha|+\nu\le M\\\nu\le N-1}}
|L^\nu \partial^\alpha u'(s,\cd)|\sum_{\substack{|\alpha|+\nu\le
M-190\\\nu\le N}} |L^\nu \partial^\alpha u'(s,\cd)|\Bigr\|_2\:ds
\\+C\int_0^t \Bigl\|\sum_{\substack{|\alpha|+\nu\le 190\\\nu\le N-1}}
|L^\nu \partial^\alpha u(s,\cd)| \sum_{\substack{|\alpha|+\nu\le 180\\\nu\le
N}} 
|L^\nu \partial^\alpha u(s,\cd)| \sum_{\substack{|\alpha|+\nu\le
M\\\nu\le N}} 
|L^\nu \partial^\alpha u'(s,\cd)|\Bigr\|_2\:ds
\\+C\int_0^t \Bigl\|\Bigl(\sum_{\nu\le N-1} |L^\nu u(s,\cd)|\Bigr)^2
\sum_{\nu\le N} |L^\nu u(s,\cd)|\Bigr\|_2\:ds.
\end{multline}
By \eqref{5.1}, \eqref{5.13}, and \eqref{5.15}, it follows that the
last term is $O(\varepsilon^3)$.  Applying \eqref{5.1}, \eqref{5.13},
and \eqref{5.15} to the first and third terms and using \eqref{2.26}
and the Schwarz inequality on the second, we see that \eqref{D.20} is
\begin{multline*}
\le C\varepsilon^3+ 
C\varepsilon \int_0^t (1+s)^{-1+\tilde{b}_N\varepsilon+\tilde{c}_N\varepsilon+}
\sum_{\substack{|\alpha|+\nu\le M \\\nu\le N}} \|L^\nu \partial^\alpha
u'\|_2\:ds
\\+ C\sum_{\substack{|\alpha|+\nu\le M+2\\\nu\le N-1}} \|\langle
x\rangle^{-1/2} L^\nu Z^\alpha u'(s,\cd)\|_{L^2(S_t)}
\sum_{\substack{|\alpha|+\nu\le M-190\\\nu\le N}} \|\langle
x\rangle^{-1/2} L^\nu \partial^\alpha u'\|_{L^2(S_t)}.
\end{multline*}
When $N\le 190$, the last term above is unnecessary, and the bound
follows from \eqref{D.17}.  For $N>190$, one uses \eqref{D.17} and the
inductive hypothesis \eqref{5.20} to bound the additional term.  Since
the same arguments can be employed to bound the last term in
\eqref{D.19}, this finishes the proof of the bound for the first term
in the left side of \eqref{D.18}.

To control the second term on the left side of \eqref{D.18}, we will
use \eqref{4.10}.  The main step is to estimate the first term on its
right.  Here, we have
\begin{multline}\label{D.21}
\sum_{\substack{|\alpha|+\nu\le M-3\\\nu\le N}} \|\Box_\gamma L^\nu
Z^\alpha u(t,\cd)\|_2 \le C\Bigl\|\sum_{|\alpha|\le 200} |Z^\alpha
u'(t,\cd)| \sum_{\substack{|\alpha|+\nu\le M-3\\\nu\le N}} |L^\nu
Z^\alpha u'(t,\cd)|\Bigr\|_2
\\+C\Bigl\|\sum_{|\alpha|\le M-3} |Z^\alpha
u'(t,\cd)|\sum_{\substack{|\alpha|+\nu\le M-203\\\nu\le N}} |L^\nu
Z^\alpha u'(t,\cd)|\Bigr\|_2
\\+C\Bigl\|\Bigl(\sum_{\substack{|\alpha|+\nu\le M-3\\\nu\le N-1}}
|L^\nu Z^\alpha u'(t,\cd)|\Bigr)^2\Bigr\|_2
\\+C \Bigl\|\Bigl(\sum_{\substack{|\alpha|+\nu\le 190\\\nu\le
N-1}} |L^\nu Z^\alpha u(t,\cd)|\Bigr)^2
\sum_{\substack{|\alpha|+\nu\le M-3\\\nu\le N}} |L^\nu Z^\alpha
u'(t,\cd)|\Bigr\|_2
\\+C \Bigl\|\sum_{\substack{|\alpha|+\nu\le 190\\\nu\le N-1}}
|L^\nu Z^\alpha u(t,\cd)| \sum_{\substack{|\alpha|+\nu\le 180\\\nu\le
N}} |L^\nu Z^\alpha u(t,\cd)| \sum_{\substack{|\alpha|+\nu\le
M-3\\|\alpha|\ge 1, \nu\le N}} |L^\nu Z^\alpha u(t,\cd)|\Bigr\|_2
\\+C\Bigl\|\sum_{\substack{|\alpha|+\nu\le 190\\\nu\le N-1}}
|L^\nu Z^\alpha u(t,\cd)| \sum_{\substack{|\alpha|+\nu\le 180\\\nu\le
N}} |L^\nu Z^\alpha u(t,\cd)|\sum_{\substack{|\alpha|+\nu\le
M-3\\\nu\le N-1}} |L^\nu Z^\alpha u'(t,\cd)|\Bigr\|_2
\\+C\Bigl\|\Bigl(\sum_{\mu\le N-1} |L^\mu u(t,\cd)|\Bigr)^2
\sum_{\mu\le N} |L^\mu u(t,\cd)|\Bigr\|_2.
\end{multline}
With $Y_{M-N-3,N}(t)$ as in \eqref{4.9}, we can apply \eqref{5.1} and
\eqref{5.16} to bound the first term in the right by
$$\frac{C\varepsilon}{1+t}Y^{1/2}_{M-N-3,N}(t).$$
By applying \eqref{5.1} and \eqref{5.13}, the same bound holds for the
fourth term in the right side of \eqref{D.21}.  Using \eqref{2.26},
the second and third terms in the right of \eqref{D.21} are controlled
by
$$C\sum_{\substack{|\alpha|+\nu\le M-1\\\nu\le N-1}} \|\langle
x\rangle^{-1/2} L^\nu Z^\alpha u'(t,\cd)\|^2_2 +
C\sum_{\substack{|\alpha|+\nu\le M-203\\\nu\le N}} \|\langle
x\rangle^{-1/2} L^\nu Z^\alpha u'(t,\cd)\|^2_2.$$
Since the coefficients of $Z$ are $O(|x|)$, one may use \eqref{5.1}, \eqref{5.13}
and \eqref{5.15} to bound the fifth and sixth terms in the right side of
\eqref{D.21} by
$$
C\varepsilon^2 (1+t)^{-1+\tilde{b}_N\varepsilon+
\tilde{c}_N\varepsilon+} \Bigl(\sum_{\substack{|\alpha|+\nu\le M-4\\\nu\le
N}} \|L^\nu Z^\alpha u'(t,\cd)\|_2+\sum_{\substack{|\alpha|+\nu\le
M-3\\\nu\le N-1}} \|L^\nu Z^\alpha u'(t,\cd)\|_2\Bigr)
$$
Finally, the last term in \eqref{D.21} is easily seen to be $\le
C\varepsilon^3 (1+t)^{-1-}$ by \eqref{5.1}, \eqref{5.13}, and
\eqref{5.15}.

If we use these bounds for \eqref{D.21} in \eqref{4.10}, we see that
\begin{multline}\label{D.22}
\partial_t Y_{M-N-3,N}(t)\le \frac{C\varepsilon}{1+t} Y_{M-N-3,N}(t)
\\+CY^{1/2}_{M-N-3,N}(t)\Bigl[\sum_{\substack{|\alpha|+\nu\le
M-1\\\nu\le N-1}} \|\langle x\rangle^{-1/2} L^\nu Z^\alpha
u'(t,\cd)\|_2^2 \\+ \sum_{\substack{|\alpha|+\nu\le M-203\\\nu\le N}}
\|\langle x\rangle^{-1/2} L^\nu Z^\alpha u'(t,\cd)\|_2^2
\\+\varepsilon^2 (1+t)^{-1+\tilde{b}_N\varepsilon
+\tilde{c}_N\varepsilon+} \Bigl(\sum_{\substack{|\alpha|+\nu\le
M-4\\\nu\le N}} \|L^\nu Z^\alpha u'(t,\cd)\|_2 +
\sum_{\substack{|\alpha|+\nu\le M-3\\\nu\le N-1}} \|L^\nu Z^\alpha
u'(t,\cd)\|_2\Bigr)
\\+\varepsilon^3 (1+t)^{-1-}\Bigr] + C\sum_{\substack{|\alpha|+\nu\le
M-2\\\nu\le N}} \|L^\nu \partial^\alpha u'(t,\cd)\|^2_{L^2(|x|<1)}. 
\end{multline}
Thus, by Gronwall's inequality, we have
\begin{multline}\label{D.23}
Y_{M-N-3,N}(t)\le C(1+t)^{2C\varepsilon} \Bigl[\varepsilon +
\sum_{\substack{|\alpha|+\nu \le M-1\\\nu\le N-1}} \|\langle
x\rangle^{-1/2} L^\nu Z^\alpha u'\|^2_{L^2(S_t)}
\\+\sum_{\substack{|\alpha|+\nu\le M-203\\\nu\le N}} \|\langle
x\rangle^{-1/2} L^\nu Z^\alpha u'\|^2_{L^2(S_t)}
\\+\varepsilon^2
\int_0^t
(1+s)^{-1+\tilde{b}_N\varepsilon+\tilde{c}_N\varepsilon+}
\sum_{\substack{|\alpha|+\nu\le M-4\\\nu\le N}} \|L^\nu Z^\alpha
u'(s,\cd)\|_2 \:ds
\\+ \varepsilon^2 \int_0^t
(1+s)^{-1+\tilde{b}_N\varepsilon+\tilde{c}_N\varepsilon+} 
 \sum_{\substack{|\alpha|+\nu\le M-3\\\nu\le N-1}} \|L^\nu Z^\alpha
u'(s,\cd)\|_2\:ds\Bigr]^2 
\\+ C(1+t)^{C\varepsilon}
\sum_{\substack{|\alpha|+\nu\le M-2\\\nu\le N}} \|\langle
x\rangle^{-1/2} L^\nu \partial^\alpha u'\|^2_{L^2(S_t)}
\end{multline}
since $Y_{M-N-3,N}(0)\le C\varepsilon^2$ by \eqref{1.11}.

For $M\le 203$, the third term on the right does not appear, and the
proof of the bound
$$\sum_{\substack{|\alpha|+\nu\le M-3\\\nu\le N}} \|L^\nu Z^\alpha
u'(t,\cd)\|_2^2 \le CY_{M-N-3,N}(t)\le C\varepsilon^2
(1+t)^{C'\varepsilon+C'\sigma}$$
 is completed by applying the inductive hypothesis
\eqref{5.20} (with $N$ replaced by $N-1$) to the second and fifth
terms on the right, applying \eqref{D.17} to the fourth term on the
right, and using the bound for the first term on the left of
\eqref{D.18} to control the last term in \eqref{D.23}.  For $M>203$, a
subsequent application of \eqref{D.17} to the third term in
\eqref{D.23} completes the proof of the bound for the second term in
the left of \eqref{D.18}.  

Using \eqref{4.12} and the arguments above, this in turn implies that
the third term in \eqref{D.18} is controlled by its right side, which
completes the proof of the lemma. \qed

By \eqref{D.16} and the lemma, we get \eqref{D.10}.  The inductive
argument using the lemma also yields \eqref{5.20} which completes the
proof of Theorem \ref{theorem1.1}.


\end{document}